\documentclass{amsart}

\newtheorem%
{thm}{Theorem}[section]
\newtheorem%
{proposition}[thm]{Proposition}
\newtheorem%
{lemma}[thm]{Lemma}
\newtheorem%
{lemmadef}[thm]{Lemma-Definition}
\newtheorem%
{corollary}[thm]{Corollary}
\newtheorem%
{conjecture}[thm]{Conjecture}
\newtheorem%
{pre-thm}[thm]{Pre-theorem}
\newcommand{\dontprint}[1]{\relax}
\newcommand{\SL}{{\mathit{SL}}}
\newcommand{\ISL}{{\mathit{ISL}}}

\usepackage{latexsym, amsmath, amssymb, amsfonts, amscd}
\usepackage{amsthm}
\usepackage{graphicx, pb-diagram}

\theoremstyle{definition}
\newtheorem{definition}[thm]{Definition}

\theoremstyle{remark}
\newtheorem{remark}[thm]{Remark}
\newtheorem{ep}[thm]{Example}

\newcommand{\lra}{\longrightarrow}
\newcommand{\rra}{\rightrightarrows}

\newcommand{\Z}{\ensuremath{\mathbb Z}}
\newcommand{\Q}{\ensuremath{\mathbb Q}}
\newcommand{\ZZ}{\ensuremath{\mathbb Z}}
\newcommand{\C}{\ensuremath{\mathbb C}}
\newcommand{\R}{\ensuremath{\mathbb R}}

\newcommand{\CP}{\mathbb{C}P}
\newcommand{\RP}{\mathbb{R}P}

\newcommand{\cL}{\mathcal {L}}
\newcommand{\cC}{\mathcal{C}}            
\newcommand{\cE}{\mathcal{E}}
\newcommand{\cF}{\mathcal{F}}
\newcommand{\cM}{\mathcal{M}}
\newcommand{\cO}{\mathcal{O}}
\newcommand{\cX}{\mathcal{X}}
\newcommand{\cY}{\mathcal{Y}}

\newcommand{\cG}{\mathcal{G}}

\newcommand{\cU}{\mathcal{U}}
\newcommand{\cV}{\mathcal{V}}

\newcommand{\cd}{\check{\delta}}

\DeclareMathOperator{\id}{id}

\DeclareMathOperator{\Hom}{Hom}

\newcommand{\tU}{\tilde{U}}

\newcommand{\bt}{\mathbf{t}}                  
\newcommand{\bs}{\mathbf{s}}                  

\usepackage{latexsym, amsmath, amssymb, amsfonts, amscd}
\usepackage{amsthm}
\usepackage{graphicx, pb-diagram}
\usepackage[all,poly,web,knot]{xy}

\renewcommand{\latticebody}{\drop@{ }}

\DeclareMathOperator{\colim}{colim}

\newcommand{\Oanx}{\cO^\times_{\!an}}
\newcommand{\Oan}{\cO_{\!an}}
\newcommand{\Opol}{\cO_{\!pol}}
\newcommand{\Ole}[1]{\cO_{\!\leq #1}}

\newcommand{\Tr}{{\mathcal T\! r}}
\newcommand{\Homm}{\mathcal{H}\hspace{-.5mm}om}

\def\undertilde#1{\put(0,-4){$\scriptscriptstyle\sim$}{#1}}

\def\rrarrow{\hspace{.05cm}\mbox{\put(0,-2){$\rightarrow$}\put(0,2){$\rightarrow$}\hspace{.5cm}}}
\def\rrrarrow{\hspace{.05cm}\mbox{\,\put(0,-3){$\rightarrow$}\put(0,1){$\rightarrow$}\put(0,5){$\rightarrow$}\hspace{.5cm}}}
\def\rrrrarrow{\hspace{.05cm}\mbox{\,\put(0,-3.5){$\rightarrow$}\put(0,0){$\rightarrow$}\put(0,3.5){$\rightarrow$}\put(0,7){$\rightarrow$}
               \hspace{.5cm}}}
\def\hookright{\ar[r]|<{\put(-3,5.2){$\scriptstyle \subset$}}}

\title{A gerbe for the elliptic gamma function}

\author[G. Felder, A. Henriques, C. A. Rossi, C. Zhu]{Giovanni Felder, Andr\'e Henriques, Carlo A. Rossi and Chenchang Zhu}
\address{Department of Mathematics, ETH Zurich, 8092 Zurich, Switzerland}
\address{Mathematical Institute, Westf\" alische Wilhelms-Universit\"at, 48149 M\" unster, Germany}
\address{Department of Mathematics, ETH Zurich, 8092 Zurich, Switzerland}
\address{Department of Mathematics, ETH Zurich, 8092 Zurich, Switzerland}
\begin{document}
\rightline{NSF-KITP-05-77}

\bigskip

\begin{abstract} The identities for elliptic gamma functions discovered by
A. Var\-chenko and one of us are generalized to an infinite set of identities for
elliptic gamma functions associated to pairs of planes in 3-dimensional space.
The language of stacks and gerbes gives a natural framework for
a systematic description of these identities and their domain of validity.
A triptic curve is the quotient of the complex plane by
a subgroup of rank three (it is a stack).
Our identities can be summarized by saying that elliptic gamma functions
form a meromorphic section of a hermitian
holomorphic abelian gerbe over the universal oriented triptic curve.
\end{abstract}
\maketitle
\setcounter{tocdepth}{1}
\tableofcontents
\section{Introduction}\label{s-i}
The elliptic gamma function \cite{R} is a
solution of an elliptic version of the
Euler functional equation $\Gamma(z+1)=z\,\Gamma(z)$, in which the rational function $z$ is replaced
by a theta function. In the conventions used in \cite{FV1},
\begin{gather*}
\Gamma(z+\sigma,\tau,\sigma)=
\theta_0(z,\tau)\Gamma(z,\tau,\sigma),\\
\theta_0(z,\tau)=\prod_{j=0}^\infty
(1-e^{2\pi i((j+1)\tau-z)})
(1-e^{2\pi i(j\tau+z)}).
\end{gather*}
As in this difference equation three periods
$1,\tau,\sigma$ are involved, one should expect
that studying modular properties of this function
should involve $\ISL_3(\mathbb Z)=\SL_3(\mathbb Z)\ltimes\mathbb Z^3$, just as theta functions have transformation properties under $\ISL_2(\mathbb Z)=\SL_2(\mathbb Z)\ltimes\mathbb Z^2$.
Indeed, in \cite{FV1} three-term functional relations involving $\Gamma$ at points differing by a
$\ISL_3(\mathbb Z)$-action by fractional linear transformations were discovered. An interpretation
of these identities was given in terms of group cohomology.

In this paper we complete the picture and show that the elliptic gamma function
is a special case of a family of infinite products defining
a {\em meromorphic section of a holomorphic gerbe} on the complex
stack $\cX_3=[X_3/\ISL_3(\mathbb Z)]$, or, equivalently, an equivariant section
of a holomorphic equivariant gerbe on $X_3$. Here $X_3$
denotes the total space of the
restriction to $\CP^2-\RP^2$ of the dual tautological bundle $O(1)\to\CP^2$.
This statement, when translated into concrete terms, implies an infinite
set of identities between infinite products generalizing the identities
found in \cite{FV1}.

Our result is a gerbe version of the fact that the theta function
$\theta_0(z,\tau)$ is
a holomorphic section of a line bundle on the universal elliptic curve.
The analogy is clearer if we extend the definition of $\theta_0(z,\tau)$ from $\mathrm{Im}\,\tau>0$ to
$\mathrm{Im}\,\tau\neq0$ by the rule $\theta_0(z,-\tau)=\theta_0(-z,\tau)^{-1}$.
Remarkably, this extended theta function has transformation properties under $\ISL_2(\mathbb Z)$
given by the same formulae, i.e., the multipliers are given by analytic continuation
from the domain $\mathrm{Im}\,\tau>0$. Passing to homogeneous coordinates,
 these transformation properties can
be rephrased by saying that $\theta_0(w/x_2,x_1/x_2)$
is a meromorphic section of a
line bundle on the stack $\cX_2=[X_2/\ISL_2(\mathbb Z)]$, where $X_2$ is the total space
of $O(1)\to(\CP^1-\RP^1)$. Equivalently, $\theta_0$ is an
equivariant meromorphic section of an equivariant bundle on $X_2$.

Equivariant holomorphic (abelian) gerbes can be understood as a
generalisation of equivariant line bundles. Here is a naive account of
this story, based on local trivializations and cocycles, as in \cite{Chatterjee}
and \cite{hitchin:gerbe}.
A more intrinsic approach is used in Section \ref{s-pseudodivisors} and
the Appendix. Recall that a holomorphic line bundle on
a complex manifold $X$ can be described in terms of an open
covering $\mathcal U=(V_a)_{a\in I}$ of $X$ by holomorphic nowhere vanishing
transition functions $\phi_{a,b}\in\mathcal O^\times(V_a\cap V_b)$
such that $\phi_{b,a}=\phi_{a,b}^{-1}$ and obeying the cocycle
condition \[\phi_{a,b}(x)\phi_{b,c}(x)=\phi_{a,c}(x),\qquad x\in V_a\cap
V_b\cap V_c\]
on triple intersections.
Now suppose that a group $G$ acts on  $X$ by holomorphic diffeomorphisms, and
that $\mathcal U$ is an invariant open covering, namely that there is a $G$ action on the index set $I$ such that
$V_{ga}=gV_a$ for all $g\in G$, $a\in I$. An equivariant structure on a line
bundle $L$ is given by a lift of the action, namely the data of linear isomorphisms $\phi(g,x)\colon L_{g^{-1}x}\to L_x$ between fibres for each
$g\in G$ depending holomorphically on $x$ and such that $\phi(gh;x)=\phi(g;x)\phi(h;g^{-1}x)$.
In terms of local trivialisations
on the charts $V_a$ the lift of the action of $g\in G$ from $V_{g^{-1}a}$ to $V_{a}$
is given by holomorphic functions $\phi_a(g;\cdot)\in\mathcal O^\times(V_{a})$
such that
\begin{gather*}
\phi_{g^{-1}a,g^{-1}b}(g^{-1}x)\phi_b(g;x)=\phi_a(g; x)\phi_{a,b}(x),\qquad x\in V_a\cap V_b,
\\
\phi_a(gh;x)=\phi_a(g;x)\phi_{g^{-1}a}(h;g^{-1}x),\qquad x\in V_a.
\end{gather*}
The first equation expresses the fact that the $\phi_a(g;x)$
define a global bundle map and the second is the action property.
Let $\mathcal M^\times$ denote the sheaf of non-zero meromorphic functions.
A (non-zero) meromorphic section of $L$ is then a collection $s_a\in\mathcal M^\times(V_a)$,
$a\in I$
such that $s_b=s_a\phi_{a,b}$ and an equivariant meromorphic section obeys
additionally $s_a(x)=\phi_a(g;x)s_{g^{-1}a}(g^{-1}x)$.


These equations mean that $\phi_{a,b},\phi_a$ are components of a group-\v Cech 1-cocycle
$\phi\in C_G^{0,1}(\mathcal{U},\mathcal O^\times)\oplus C_G^{1,0}(\mathcal{U},\mathcal O^\times)$
in the double complex $\oplus_{p,q}C_G^{p,q}(\mathcal{U},\mathcal O^\times)$. For any
$G$-equivariant sheaf $\mathcal F$ of abelian groups on $X$, such as $\mathcal
O^\times,\mathcal M^\times$,
the \v Cech complex $\check C^\cdot(\mathcal{U},\mathcal F)$ is a complex of $G$-modules
and thus we have a double complex
\[
C_G^{p,q}(\mathcal U,\mathcal F)=C^p(G,\check C^q(\mathcal{U},\mathcal{F})),
\]
with group cohomology differential and \v Cech differential.
The total cohomology of this double complex
maps to the equivariant cohomology of the sheaf $\mathcal F$.
If $D$ denotes the total differential of this complex, a meromorphic
section of a holomorphic $G$-equivariant line bundle associated to the 1-cocycle $\phi$ is then
given by a 0-chain $s\in C_G^{0}(\mathcal U,\mathcal M^\times)$ such
that $Ds=\phi\in C_G^1(\mathcal U,\mathcal O^\times)\subset
C_G^1(\mathcal U,\mathcal M^\times)$.

These notions can be easily extended to the next level: one can define an
equivariant holomorphic gerbe (with structure group $\mathbb C^\times$ and
trivial band) on a complex $G$-manifold $X$ to be a cocycle
$\phi\in C_G^2(\mathcal U,\mathcal O^\times)$ and an equivariant meromorphic
section to be a cochain $s\in C^1_G(\mathcal U,\mathcal M^\times)$ such that
$Ds=\phi$.
{}From a more geometric  point of view, gerbes are obtained from line
bundles by ``replacing functions by line bundles''. Thus a {\em holomorphic
 gerbe} is given in terms of the covering $\mathcal U$ by a collection
$L_{a,b}$ of line bundles on $V_a\cap V_b$, the {\em transition bundles}, and isomorphisms
\begin{equation}\label{e-iso1}
\phi_{a,b,c}\colon L_{a,c}\to L_{a,b}\otimes L_{b,c}\quad \text{on} \;V_a\cap V_b\cap V_c.
\end{equation}
The bundles are such
that $L_{a,a}$ is the trivial bundle on $V_a$ and $L_{b,a}=L_{a,b}^{-1}$, the dual bundle. The
isomorphisms $\phi_{a,b,c}$ are skew-symmetric in $a,b,c$, when viewed as invertible sections
of $L_{a,b}\otimes L_{b,c}\otimes L_{c,a}$ on $V_a\cap V_b\cap V_c$, and are
coassociative on fourfold intersections:
\[
(\phi_{a,b,c}\otimes\mathrm{id})\circ\phi_{a,c,d}=
(\mathrm{id}\otimes\phi_{b,c,d})\circ\phi_{a,b,d}\quad \text{on}\;V_a\cap V_b\cap V_c\cap V_d.
\]
Suppose now that $G$ acts on $X$ and $\mathcal U$ is a $G$-invariant open covering.
A {\em $G$-equivariant holomorphic gerbe} is a holomorphic gerbe on $X$ with the additional data of line bundles
$L_a(g)\to V_a$,  $g\in G$, $a\in I$ and isomorphisms
\begin{gather}
\phi_{a,b}(g)\colon(g^{-1})^*L_{g^{-1}a,g^{-1}b}\otimes L_b(g)\to
 L_a(g)\otimes L_{a,b}\quad\text{on}\; V_a\cap V_b\label{e-iso2}
\\
 \phi_a(g,h)\colon L_a(gh)\to L_a(g)\otimes (g^{-1})^*L_{g^{-1}a}(h)\quad\text{on}\; V_a,\label{e-iso3}
\end{gather}
obeying natural coherence conditions.
If $L_{a,b}$ is the trivial bundle $V_a\cap V_b\times\mathbb C$, then $\phi_{a,b,c},\phi_{a,b}(g),\phi_a(g,h)$ are
invertible
holomorphic functions and both associative and coherence condition mean that they
form  a 2-cocycle in $C_G^2(\mathcal U,\mathcal O^\times)=C_G^{0,2}(\mathcal U,\mathcal O^\times)\oplus C_G^{1,1}(\mathcal U,\mathcal O^\times)\oplus C_G^{2,0}(\mathcal U,\mathcal O^\times)$.

In the language of transition bundles, a meromorphic section of the gerbe defined by the data
$(\mathcal U,\phi)$ is just a collection of meromorphic sections $s_{a,b}$ of $L_{a,b}$
such that $s_{a,b}=s_{b,a}^{-1}$ and
$\phi_{a,b,c}\circ s_{a,c}=s_{a,b}\otimes s_{b,c}$. Note that if $s_{a,b}$
are  {\em holomorphic} sections then they form what
is called a {\em trivialisation} of the gerbe.
 An {\em equivariant meromorphic
section} $s$ is a
a collection of meromorphic sections $s_{a,b}$ of $L_{a,b}$, forming a meromorphic section, and meromorphic sections $s_a(g)$ of $L_a(g)$ compatible with
the maps $\phi$. In terms of local trivialisations, these conditions are equivalent
to $Ds=\phi$.

\dontprint{
Thus to every 2-cocycle we can associate a bundle gerbe with trivial $L_{a,b}$ and
conversely to every bundle gerbe with trivialisation (choice of isomorphism
$L_{a,b}\to (V_{a}\cap V_{b})\times \mathbb C$ for all $a,b$) we associate a 2-cocycle $\phi$, so that
different trivialisations give rise to equivalent cocycles.
}
\dontprint{The coherence conditions are best understood in terms of sheaves of categories.
In that language a holomorphic gerbe on a complex manifold
is a sheaf of categories locally equivalent to the category of holomorphic
line bundles and obeying certain axioms.
In particular, given a good open covering $\mathcal U$, we have a collection
of categories $\mathcal C(V_{a_0}\cap\cdots\cap V_{a_q})$ and restriction functors
$\mathcal C(V')\to\mathcal C(V'')$ whenever the multiple intersection $V''$
is contained in $V'$........................................................
}

Let now $X=X_3$ be the total space of
the restriction of the line bundle
$O(1)\to\mathbb C P^2$ to $\mathbb C P^2
-\mathbb R P^2$. A point of $X$ is
an equivalence class of pairs
$(w,x)\in\mathbb
C\times (\mathbb C^3-\mathbb C\cdot\mathbb R^3)$ modulo $(w,x)\cong (\lambda w,\lambda x)$, $\lambda\in\mathbb C^\times$. The
group $G=\ISL_3(\mathbb Z)=\SL_3(\mathbb Z)\ltimes \mathbb Z^3$ acts on $X_3$ as follows: $n\in\mathbb Z^3$ acts by $(w,x)\mapsto (w-n\cdot x,x)$
and $g\in \SL_3(\mathbb Z)$ acts by $(w,x)\mapsto (x,gx)$. There is an
invariant open covering $\mathcal U=(V_a)$ indexed by the set of
oriented 2-dimensional subspaces of
$\mathbb Q^3$: $V_a$ is the set of
$x\in\mathbb C^3$ such that $\mathrm{Im}
(\alpha\cdot x\overline{\beta\cdot x)}>0$
for some (and thus any) oriented basis
$(\alpha,\beta)$ of
$a$.

\begin{thm}\label{t-Quiverful} The elliptic gamma function
$\Gamma(w/x_3,x_1/x_3,x_2/x_3)$ is the
component associated to the x-z and
y-z planes of a meromorphic section
\[
 s=(\Gamma,\Delta)\in C_G^1(\mathcal U,\mathcal M^\times)=C_G^{0,1}
(\mathcal U,\mathcal M^\times)\oplus
C_G^{1,0}(\mathcal U,\mathcal M^\times)
\]
of an equivariant holomorphic
gerbe $\mathcal G$ on $X$ given by a 2-cocycle
$\phi\in C_G^2(\mathcal U,\mathcal
O^\times)$. In explicit terms, we have
identities
\begin{gather*}
\Gamma_{a,b}(y)\Gamma_{b,a}(y)=1,\quad y\in V_a\cap V_b,
\\
\phi_{a,b,c}(y)\Gamma_{a,c}(y)=\Gamma_{a,b}(y)\Gamma_{b,c}(y),\quad y\in V_a\cap V_b\cap V_c,
\\
{\phi_{a,b}(g;y)}{\Gamma_{g^{-1}a,g^{-1}b}(g^{-1}y)}{\Delta_b(g;y)}={\Delta_a(g; y)}{\Gamma_{a,b}(y)}
,\qquad
y\in V_a\cap V_b,
\\
\phi_a(g,h;y)\Delta_a(gh;y)=\Delta_a(g;y)\Delta_{g^{-1}a}(h;g^{-1}y),
\qquad y\in V_a,
\end{gather*}
for  all $a,b,c\in I, g,h\in G$.
\end{thm}

The functions $\Gamma, \Delta,\phi$ are defined in Section \ref{s-iii} and this theorem
is proved in Section \ref{s-iv} (see Theorem \ref{t-main}). We call {\em gamma gerbe}
the gerbe defined by the cocycle $\phi$.

In Section \ref{s-pseudodivisors} we give a more geometric description of this gerbe.
First of all, just like $\cX_2$ is the universal curve over the moduli stack
$[(\CP^1-\RP^1)/\SL_2(\mathbb Z)]$ of elliptic curves,
$\cX_3$ is the universal curve over the moduli stack $[(\CP^2-\RP^2)/\SL_3(\mathbb Z)]$ of {\em oriented triptic curves}: an oriented triptic curve is a stack of
the form $\cE=[\mathbb C/\iota(\mathbb Z^3)]$ for some $\mathbb R$-epimorphism
$\iota\colon \mathbb R^3\to\mathbb C$, equipped with a generator of $H^3(\cE,\mathbb Z)\simeq \mathbb Z$.
On this universal curve we introduce geometrically the {\em pseudodivisor
gerbe}; it is given as a sheaf of categories, see Section \ref{s-pseudodivisors}.
We then show that the gamma gerbe $\cG^\Gamma$ and the pseudodivisor gerbe
$\cG^{div}$ are isomorphic by giving a holomorphic trivialisation of
$\cG^\Gamma\otimes(\cG^{div})^*$.

Holomorphic $G$-equivariant gerbes on $X$ are classified by
the second equivariant cohomology group $H_G^2(X,\mathcal O^\times)$ \cite{Gir71, Bry},
the Brauer group $\mathrm{Br}(\cX)$ of the stack $\cX=[X/G]$.
The section $s$ defines a class
$[s]\in H^1_{G}(X,\mathcal M^\times/\mathcal O^\times)$ mapping by the connecting
homomorphism to the  class $[\phi]\in H^2_G(X,\mathcal O^\times)$ of the
gerbe. Thus, in analogy with the case of holomorphic line bundles, we can say that
$s$ is a {\em Cartier divisor} for the gerbe $\mathcal G$.
The exponential short exact sequence gives then
the {\em Dixmier--Douady class}  $c(\mathcal G)$ of the gerbe, which is the image of $[\phi]$ by the connecting
homomorphism $H^2_G(X,\mathcal O^\times)\to H^3_G(X,\mathbb Z)$.
In Section \ref{s-Brauer} we compute the Brauer group $Br(\cE)=\cE\times \Z$
of a triptic curve $\cE$ and identify the class of the restriction of the
gamma gerbe with $(0,1)$. For the universal triptic curve, we compute the
$H^3_G(X_3,\Z)$ modulo torsion and obtain in Section \ref{sec-charclass}
the following result.
\begin{thm}\label{t-Aquarius}
We have
\[
H^3_G(X,\mathbb Z)/\mathrm{torsion}\cong \mathbb Z\oplus\mathbb Z
\]
and the image of $c(\cG)$ is a primitive vector in the lattice $\mathbb Z\oplus\mathbb Z$.
\end{thm}

In particular the gerbe is non-trivial.
We prove a more precise version of this statement
below (Theorem \ref{t-sonostufo}).

Finally, we show in Section \ref{s-v}
that the gamma gerbe is a {\em hermitian} holomorphic gerbe.
In terms of local trivialisations, this means that the line bundles $L_{a,b}$,
$L_a(g)$ have a hermitian structure and the structure maps \eqref{e-iso1},
\eqref{e-iso2}, \eqref{e-iso3} are unitary.
This result implies that, in spite of the fact that the quotient space
is very singular, the gamma gerbe admits an {\em equivariant connective structure}, i.e.,
a collection of connections on $L_{a,b}$, $L_a(g)$ compatible with the structure maps. Indeed, we prove (Theorem \ref{thm-eqgerbesconnstr}) that every hermitian equivariant
holomorphic gerbe on a complex manifold has a unique equivariant connective structure that is at the same time
 compatible with the complex and the hermitian structure, generalizing
a result of Chatterjee \cite{Chatterjee}. It
is a gerbe analogue of the theorem that holomorphic hermitian line bundles admit
a unique compatible connection. In the case of line bundles, the curvature of
a connection gives a de Rham representative for the first Chern class and it is
a closed form of type (1,1) if the connection comes from a hermitian structure
on a holomorphic line bundle. Similarly, hermitian forms on holomorphic gerbes
on a complex $G$-manifold $X$
have an invariant, the {\em (1,1)-curvature} in
$H^1_G(X,\underline\Omega^{1,1}_{\mathit{cl}})$, the first equivariant cohomology
group
with values in the sheaf of closed (1,1)-forms, see Section \ref{s-v}. An equivariant
 \v Cech representative for this class is given by the collection of the
 curvatures of the connections on $L_{a,b}$, $L_a(g)$. The natural map
$H^1_G(X,\underline\Omega^{1,1}_{\mathit{cl}})\to H^3_{G\, dR}(X,\mathbb C)$ gives
then an equivariant de Rham realization of the Dixmier--Douady class $c(\cG)\in H^3_G(X,\mathbb Z)$ (up to a factor $i/2\pi$).

The Appendix contains a description of gerbes on stacks as central extension of
groupoids, as well as results on the cohomology on stacks that we use in the
paper.

Both the cocycle defining the gamma gerbe and the hermitian structure are
constructed using multiple Bernoulli polynomials. The importance of these polynomials
in this setting was
first recognized by Narukawa \cite{Na}, who extended the modular relations
of \cite{FV1} to the multiple elliptic gamma functions of Nishizawa \cite{Ni}.
Most of the algebraic identities we prove can be generalized to the
multiple gamma function setting, related to $\ISL_n(\mathbb Z)$, and our proofs extend in a straightforward way.
However, the geometric setting presented here is for the moment limited
to $n\leq 3$.

Note that the identities proven here imply
(by taking logarithms) new modular
identities for the
generalized Eisenstein series considered in
\cite{FV3}. Also, the study of the q-deformed KZB equation on elliptic
curve in \cite{FVKZB} suggests that there is a non-abelian version of this
construction for each representation of a simple Lie algebra.

\noindent{\bf Acknowledgments}
We are grateful to Damien Calaque, Pavel Etingof,
Nigel Hitchin,
Christophe Soul\'e and Ping Xu for useful discussions,
comments and suggestions on the subject of
this paper. G. F. wishes to thank particularly
Alexander Varchenko: some of the basic ideas for this work
originated in early discussions with him. G. F. thanks the organizers
of the program ``Mathematical Structures in String Theory'' and KITP, where
part of this work was done, for hospitality.
Research at KITP was supported in part by the NSF through Grant Nr.~PHY99-07949.
A. H. thanks ETH for hospitality. This work has been partially supported by the European Union through the FP6
Marie Curie RTN {\em ENIGMA} (Contract number MRTN-CT-2004-5652).

\section{Theta function and
 elliptic gamma function}\label{s-ii}
\subsection{The theta function}
For $z\in\mathbb C$, $\mathrm{Im}\,\tau>0$, let
\[
\theta_0(z,\tau)=\prod_{j=0}^\infty(1-e^{2\pi i((j+1)\tau-z))})(1-e^{2\pi i(j\tau+z)}).
\]
Up to a factor $ie^{-\pi i(z-\tau/6)}\eta(\tau)^{-1}$, where $\eta$ is
the Dedekind eta function, $\theta_0$ is the first (odd) Jacobi theta function.
It is holomorphic for $w\in\mathbb C, \mathrm{Im}\,\tau>0$.
We extend $\theta_0$ to a meromorphic function
on the domain $\mathrm{Im}\,\tau\neq 0$ by setting
$\theta_0(z,-\tau)=\theta_0(z+\tau,\tau)^{-1}$
(which is also equal to $\theta_0(-z,\tau)^{-1}$).

This function is periodic with period 1 in both
arguments. It obeys the functional equation
\[
\theta_0(z+\tau,\tau)=-e^{-2\pi i z}\theta_0(z,\tau),
\]
 and the modular relation
\[
\theta_0\left(\frac w{x_2},\frac{x_1}{x_2}\right)
\theta_0\left(\frac w{x_1},\frac{x_2}{x_1}\right)
=\exp\left(-\pi i P_2(w,x)\right)
\]
where $P_2(w,x)=
({w^2}-(x_1+x_2)w+(x_1^2+x_2^2+3x_1x_2)/6)/x_1x_2$.
It follows that $\theta_0$ transforms under $\SL_2(\mathbb Z)$ by
\[
\theta_0\left(\frac z{c\tau+d},\frac{a\tau+b}{c\tau+d}\right)
=e^{\pi i Q(g;z,\tau)}\theta_0(z,\tau),\quad g^{-1}=
\left(\begin{array}{cc}a&b\\ c&d
  \end{array}
  \right)\in \SL_2(\mathbb Z),
\]
for some polynomial $Q(g;z,\tau)\in\mathbb Q[\tau,\tau^{-1},z]$ of degree at most 2
in $w$. An explicit formula for this polynomial is
given below.
\subsection{The elliptic gamma function}\label{ss-tegf}
The elliptic gamma function is a meromorphic function $\Gamma(z,\tau,\sigma)$
defined on $\mathbb C\times(\mathbb C-\mathbb R)\times(\mathbb C-\mathbb
R)$. For $\mathrm{Im}\,\tau, \mathrm{Im}\,\sigma>0$
it is defined as a convergent infinite product
\[
\Gamma(z,\tau,\sigma)=
\prod_{j,k=0}^\infty\frac
{1-e^{2\pi i((j+1)\tau+(k+1)\sigma-z)}}
{1-e^{2\pi i(j\tau+k\sigma+z)}}\,,
\]
and is extended by the rules
\begin{equation}\label{e-Cleophile}
\Gamma(z,\tau,\sigma)=\frac1{\Gamma(z-\tau,-\tau,\sigma)}
=\frac1{\Gamma(z-\sigma,\tau,-\sigma)}\,.
\end{equation}
For the purpose of this paper, the domain where $\mathrm{Im}\,\tau<0$ and $\mathrm{Im}\,\sigma>0$
is particularly relevant. There the gamma function is given by the product
\begin{equation}\label{e-expelliatur}
\Gamma(z,\tau,\sigma)=
\prod_{j,k=0}^\infty\frac
{1-e^{2\pi i(z-(j+1)\tau+k\sigma)}}
{1-e^{2\pi i(-z-j\tau+(k+1)\sigma)}}\,,
\end{equation}
Clearly $\Gamma$ is a periodic function of period 1 in all its arguments.
Moreover, it obeys an interesting set of identities.
First of all, $\Gamma$ obeys the difference equations
\begin{eqnarray*}
\Gamma(z+\tau,\tau,\sigma)&=&\theta_0(z,\sigma)\Gamma(z,\tau,\sigma)\\
\Gamma(z+\sigma,\tau,\sigma)&=&\theta_0(z,\tau)\Gamma(z,\tau,\sigma).
\end{eqnarray*}
Moreover it obeys the ``modular'' three-term relations \cite{FV1}
\begin{equation}\label{e-sectumsempra}
\Gamma(z,\tau,\sigma)=\Gamma(z,\tau,\tau+\sigma)\Gamma(z+\sigma,\tau+\sigma,\sigma).
\end{equation}
\begin{equation}\label{e-levicorpus}
\Gamma\left(\frac w{x_3},\frac{x_1}{x_3},\frac{x_2}{x_3}\right)
\Gamma\left(\frac w{x_1},\frac{x_2}{x_1},\frac{x_3}{x_1}\right)
\Gamma\left(\frac w{x_2},\frac{x_3}{x_2},\frac{x_1}{x_2}\right)
=e^{-\pi iP_3(w,x)/3},
\end{equation}
\begin{eqnarray}\label{e-bern}
P_3(w,x)&=&\frac{w^3}{x_1x_2x_3}-3\,\frac{x_1+x_2+x_3}{2x_1x_2x_3}\,w^2
\\ &&
 +\frac{x_1^2+x_2^2+x_3^2+3x_1x_2+3x_2x_3+3x_1x_3}{2x_1x_2x_3}\,w\notag \\
&&-\frac1{4}(x_1+x_2+x_3)\left(\frac1{x_1}+\frac1{x_2}+\frac1{x_3}\right).\notag
\end{eqnarray}
\subsection{Multiple gamma functions and Bernoulli polynomials}
 The functions $G_0(z,\tau)=\theta_0(z,\tau)$, $G_1(z,\tau,\sigma)=\Gamma(z,\tau,\sigma)$ are the first two of a sequence
introduced by Nishizawa \cite{Ni}
of multiple elliptic gamma functions $G_n$ of
$n+2$ variables obeying
\[
G_n(z+\tau_i,\tau_0,\ldots,\tau_n)
=G_{n-1}(z,\tau_0,\ldots,\hat\tau_i,\cdots,\tau_n)G_{n}(z,\tau_0,\ldots,\tau_n).
\]
As shown by Narukawa \cite{Na}, these functions
obey functional relations generalizing
the modular relations of $\theta_0,\Gamma$:
\[
\prod_{k=1}^r G_{r-2}\left(
\frac w{x_k},\frac {x_1}{x_k},\ldots,
\widehat{\frac{x_k}{x_k}},\ldots, \frac{x_r}{x_k}\right)
=\exp\left(-\frac{2\pi i}{r!} B_{r,r}(w,x)\right),
\]
where $B_{r,n}$ are multiple Bernoulli polynomials
defined by the generating function
\[
e^{wt}\prod_{j=1}^r\frac t{e^{x_jt}-1}=
\sum_{n=0}^\infty B_{r,n}(w,x_1,\dots,x_r)\frac {t^n}{n!}.
\]
In particular, $P_2=B_{2,2}$, $P_3=B_{3,3}$. These polynomials obey a number
of interesting identities. One of them is the difference relation
\begin{equation}\label{e-difference}
B_{r,n}(w+x_i,x)-B_{r,n}(w,x)=B_{r-1,n-1}(w,x_1,\dots,\hat x_i,\dots,x_r),
\end{equation}
which can readily be deduced using the generating function.

\section{Gamma functions associated to wedges}\label{s-iii}

Let $\Lambda=\mathbb Z^3\subset V=\mathbb Q^3$. We consider $V$ as an oriented
vector space with the standard volume form $\det\colon\wedge^3V\to\mathbb Q$. The dual lattice
$\Lambda^\vee=\mathrm{Hom}_\mathbb Z(\Lambda,\mathbb Z)$ can be viewed as a lattice in the dual
vector space $V^*$ and also comes with an
orientation form $\mathrm{det}\colon \wedge^3\Lambda^\vee\to\mathbb Z$.
To each pair of oriented 2-dimensional subspaces of
$V^*$ we associate a function
that reduces to the elliptic gamma function for coordinate planes.

\subsection{Wedges}\label{s-wedges}
An element $v$ of a free $\mathbb Z$-module $M$
is called {\em primitive} if
$\lambda\in \mathbb Z$, $v\in\lambda M$ implies $\lambda=\pm 1$
(in particular $0$ is not primitive).
We denote by $M_{\mathit{prim}}\subset M$ the set of primitive elements of $M$.
Primitive elements in $\Lambda$ are in one-to-one correspondence with oriented planes
through the origin in $V^*$. The plane $H(a)$ corresponding to
$a\in\Lambda_{\mathit{prim}}$ has equation $\alpha(a)=0$ and divides $\Lambda^\vee$ into two
subsets:
\[
\Lambda^\vee_+(a)=\{\alpha\in\Lambda^\vee\,|\,\alpha(a)>0\},\qquad
\Lambda^\vee_-(a)=\{\alpha\in\Lambda^\vee\,|\,\alpha(a)\leq 0\}.
\]
\begin{definition} A {\em wedge} is an ordered pair of oriented planes
through $0$ in $V^*$ or, equivalently, a pair $(a,b)\in\Lambda_{\mathit{prim}}^2$.
A wedge is {\em in general position} if
its planes intersect in a line, i.e., if $a,b$ are linearly independent.
\end{definition}

The intersection line of a wedge $(a,b)$ in general position is spanned by the
linear form $\det(a,b,\cdot)\in \Lambda^\vee$. This form is not primitive in
general, but there exists a unique $\gamma=\gamma_{a,b}\in\Lambda^\vee_p$ and
positive integer $s$ such that
\[
\det(a,b,x)=s\gamma(x),\quad\forall x\in V.
\]
The lattice points in the intersection are then the integer multiples of
$\gamma$. We call $\gamma_{a,b}$ the {\em direction vector} of the wedge $(a,b)$
and $s=\mathrm{mod}(a,b)$ its {\em modulus}.
We set $\mathrm{mod}(a,b)=0$ if $a,b$ are linearly
dependent.

Let $\mathrm{Aut}(\Lambda)= \SL_3(\mathbb Z)$ be the group of linear automorphisms
of $V$ preserving the lattice and the volume form. This group acts on $\Lambda$ preserving
$\Lambda_{\mathit{prim}}$. The modulus is an invariant of a wedge. There are two orbits with
modulus zero and the number of orbits with modulus $s>0$ is given by the Euler $\phi$-function
$\phi(s)$. Indeed, if $e_1,e_2,e_3$ is a basis of $\Lambda$ with $\det(e_1,e_2,e_3)=1$,
then it is easy to see that any wedge $(a,b)$ can be brought by an $\SL_3(\mathbb Z)$
transformation to the normal form $(e_1,\pm e_1)$ or $(e_1,r e_1+se_2)$ with $s=1,2,\dots$,
$0\leq r<s$, $\gcd(r,s)=1$, and that each orbit contains exactly one wedge in normal
form.

\subsection{Gamma functions}
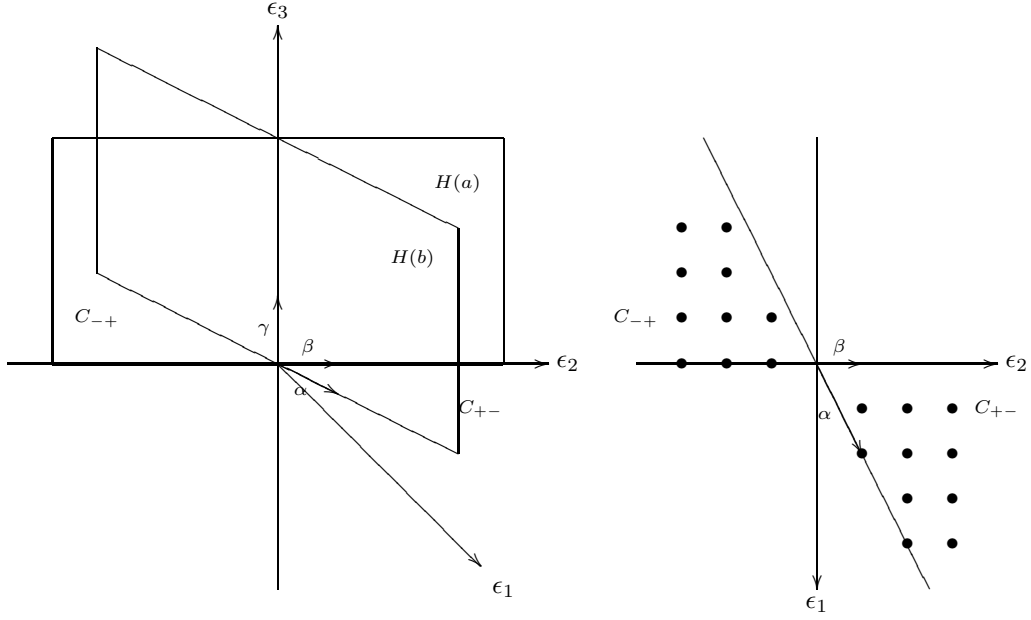
\begin{figure}
   \[
\begin{xy}
  *\xybox{(0,0);<3mm,0mm>:<0mm,3mm>::
  ,0
  ,{\xylattice{-12}{12}{-10}{15}}}="S"
  ,{(30,0)="a" 
     ,(-24,42)="b" 
   ,(-30,30)."a"*\frm{-}="L"}
   ,{ (-24, 42) \ar@{-} (24, 18)}
   ,{ (0, 0) \ar (0, 9)^{\gamma} }
   ,{ (-24, 42) \ar@{-} (-24, 30)}
,{ (-24, 30) \ar@{-} (-24, 12)},{ (-24, 12) \ar@{-} (0, 0)}
   ,{ (24, 18) \ar@{-} (24, -12)}
   ,{ (0, 0) \ar@{-} (24, -12)}
   ,{ (0, 0) \ar (8, -4)_{\alpha}}
   ,{ (0, 0) \ar (8, 0)^{\beta}}
,{ (0, 0) \ar@{->} (27, -27)}, {(30,-30)*{\epsilon_1}}
  ,{"S"+L \ar "S"+R*+!L{{\epsilon_2}}}
 ,{"S"+D \ar "S"+U*+!D{{\epsilon_3}}}
,{(24, 24)*{_{H(a)}}}, {(18, 14)*{_{H(b)}}}, {(27, -6)*{_{C_{+-}}}},  {(-24, 6)*{_{C_{-+}}}}
    \end{xy}
\quad \begin{xy}
  *\xybox{(0,0);<3mm,0mm>:<0mm,3mm>::
  ,0
  ,{\xylattice{-8}{8}{-10}{10}}}="S"
   ,{ (-15, 30) \ar@{-} (15, -30)}, {(0, 0) \ar (6,-12)_{\alpha} }
   ,{ (0, 0) \ar (6, 0)^{\beta}}
  ,{"S"+L \ar "S"+R*+!L{{\epsilon_2}}}
 ,{"S"+U \ar "S"+D*+!U{{\epsilon_1}}}
, {(24, -6)*{_{C_{+-}}}},  {(-24, 6)*{_{C_{-+}}}}
,{(-6, 0)*{\bullet}},{(-12, 0)*{\bullet}},{(-18, 0)*{\bullet}}
,{(-6, 6)*{\bullet}},{(-12, 6)*{\bullet}},{(-18, 6)*{\bullet}}
,{(-12, 12)*{\bullet}},{(-18, 12)*{\bullet}}
,{(-12,18)*{\bullet}},{(-18, 18)*{\bullet}}
,{(6,-6)*{\bullet}},{(12, -6)*{\bullet}},{(18,-6)*{\bullet}}
,{(6,-12)*{\bullet}},{(12, -12)*{\bullet}},{(18,-12)*{\bullet}}
,{(12, -18)*{\bullet}},{(18,-18)*{\bullet}}
,{(12,-24)*{\bullet}},{(18,-24)*{\bullet}}
    \end{xy}
   \]
\caption{
These pictures show the situation when $a=e_1$, $b=e_1-2e_2$.
Then $H(a)=span\{\epsilon_2, \epsilon_3\}$ and
$H(b)=span\{2\epsilon_1+\epsilon_2, \epsilon_3\}$, where $(\epsilon_i)$ is
the basis dual to the standard basis $(e_i)$ of $V$. Then we take
$\alpha=2\epsilon_1+\epsilon_2$, $\beta=\epsilon_2$ and
$\gamma=\epsilon_3$. The left picture shows the wedge of the ordered
pair $(a, b)$ and the right one shows the distribution of
$\delta$'s.}\label{f-1}
\end{figure}

Let $V_\mathbb C=V\otimes_\mathbb Q\mathbb C$ be
the complexification of $V$.
We associate to a wedge $(a,b)\in\Lambda_{\mathit{prim}}^2$ a function $\Gamma_{a,b}$
on an open subset of $\mathbb C\times V_\mathbb C$.
If $a=b$ we set
\[
\Gamma_{a,a}(w;x)=1.
\]
If $a,b$ are linearly independent, the direction vector
$\gamma=\gamma_{a,b}$ spans the intersection line $H(a)\cap H(b)\cap \Lambda^\vee$, so
the group $\mathbb Z\gamma$ acts on
\[
C_{\theta\eta}(a,b)=\Lambda^\vee_\theta(a)\cap \Lambda^\vee_\eta(b),\quad \theta,\eta\in
\{+,-\}
\]
by translation. We then set
\begin{equation}\label{e-imperius}
\Gamma_{a,b}(w,x)=\frac
{\prod_{\delta\in C_{+-}(a,b)/\mathbb Z\gamma}(1-e^{-2\pi i(\delta(x)-w)/\gamma(x)})}
{\prod_{\delta\in C_{-+}(a,b)/\mathbb Z\gamma}(1-e^{2\pi i(\delta(x)-w)/\gamma(x)})}\,.
\end{equation}
The geometry of this construction is depicted for an example in Fig.~\ref{f-1}.

\subsection{Domain of definition}\label{ss-Domain}
Let $a\in\Lambda_{\mathit{prim}}$. The plane $H(a)$ has a natural orientation: an ordered basis
$\alpha,\beta$ of $H(a)$ is oriented if $\det(\alpha,\beta,\delta)>0$
whenever $\delta(a)>0$.
\begin{lemma}\label{l-cruciatus} Let $x\in V_\mathbb C$. The following are equivalent:
\begin{enumerate}
\item[(i)] $\mathrm{Im}(\alpha(x)\overline{\beta(x)})>0$
for some oriented basis $\alpha,\beta$ of $H(a)$.
\item[(ii)] $\mathrm{Im}(\alpha(x)\overline{\beta(x)})>0$
for all oriented bases $\alpha,\beta$ of $H(a)$.
\item[(iii)] $\beta(x)\neq 0$ and $\mathrm{Im}(\alpha(x)/\beta(x))>0$
for some oriented basis of $H(a)$.
\item[(iv)] $\alpha(x)\neq 0$ and $\mathrm{Im}(\beta(x)/\alpha(x))<0$
for some oriented basis of $H(a)$.
\end{enumerate}
\end{lemma}

\begin{proof} The equivalence between (i) and (ii) follows from the fact
that the action of $\SL_2(\mathbb Q)$ by fractional linear transformations
preserves the upper half-plane. The equivalence with the other statements
is clear.
\end{proof}

Let $U^+_a$ be the open subset of $x\in V_\mathbb C$ obeying any of the equivalent
conditions of Lemma \ref{l-cruciatus}.

\begin{proposition}\label{p-Goodenough}
Let $a,b\in\Lambda_{\mathit{prim}}$ be linearly independent.
Then the ratio of infinite products \eqref{e-imperius} converges to a meromorphic function
on $U^+_a\cap U^+_b$.
 It has simple zeros at
$w=\delta(x)+n\gamma(x)$, with $\delta\in C_{+-}(a,b)$, $n\in\mathbb Z$, and simple
poles at
$w=\delta(x)+n\gamma(x)$, with $\delta\in C_{-+}(a,b)$, $n\in\mathbb Z$.
\end{proposition}

\begin{proof}
Let us first consider the denominator. Notice that
$C_{-+}(a,b)=
(\mathbb Q_{\leq 0}\alpha+\mathbb Q_{>0}\beta+\mathbb Q\gamma)\cap\Lambda^\vee$,
for any $\alpha\in H(b)$, $\beta\in H(a)$ such that $\alpha(a)>0$, $\beta(b)>0$.
Thus the condition for convergence of the infinite product in the denominator
is $\mathrm{Im}(\alpha(x)/\gamma(x))<0$ (i.e., $\mathrm{Im}(\gamma(x)/\alpha(x))>0$),
$\mathrm{Im}(\beta(x)/\gamma(x))>0$. Now, by construction, $\beta,\gamma$ is
an oriented basis of $H(a)$ and $\gamma,\alpha$ an oriented basis of $H(b)$,
so the convergence condition is $x\in U_a^+\cap U_b^+$.
A similar argument applies to the numerator.
 In this domain we obtain a meromorphic
function whose divisor can be read off the zeros of the factors of the
infinite products.
\end{proof}

\begin{remark} We did not define $\Gamma_{a,b}$ for $a=-b$.
There is no need to define it as its domain of definition is empty.
\end{remark}

\subsection{Expression in terms of ordinary elliptic gamma functions}
The gamma function associated to a wedge may be expressed (in many ways)
as a finite product of ordinary elliptic gamma functions:

\begin{proposition}\label{p-horcrux}
Let $a,b\in\Lambda_{\mathit{prim}}$ be linearly independent and let $\gamma\in\Lambda_{\mathit{prim}}$
be the direction
vector of the wedge $(a,b)$.
Let $\alpha,\beta\in\Lambda^\vee$ be such that
$\alpha(b)=\beta(a)=0$ and $\alpha(a)>0$, $\beta(b)>0$.
Then
\[
\Gamma_{a,b}(w,x)=\prod_{\delta\in F/\mathbb Z\gamma}
\Gamma\left(\frac{w+\delta(x)}{\gamma(x)},\frac{\alpha(x)}{\gamma(x)},
\frac{\beta(x)}{\gamma(x)}\right),
\]
where $F$ is the set of $\delta\in\Lambda^\vee$ such that
\[
0\leq\delta(a)<\alpha(a),\qquad 0\leq\delta(b)<\beta(b).
\]
\end{proposition}

\begin{proof} Since $C_{+-}(a,b)$ is the subset
of $\Lambda$ defined by the inequalities $\delta(a)>0,\delta(b)\leq0$, it follows
that each
$\delta\in C_{+-}(a,b)$
can be uniquely written as
$\delta=-\bar\delta+(j+1)\alpha-k\beta$,
with
$\bar\delta\in F$
and
$j,k\in\mathbb Z_{\geq0}$ (modulo $\mathbb{Z}\gamma$).
Similarly if
$\delta\in C_{-+}(a,b)$
then
$\delta=-\bar\delta-j\alpha+(k+1)\beta$
with
$\bar\delta\in F$
and
$j,k\in\mathbb Z_{\geq0}$.
\end{proof}

\begin{ep} If $a=e_1$, $b=e_2$ are the first two standard basis vectors,
$\Gamma_{a,b}$ reduces to the ordinary elliptic gamma function in the
domain $\mathrm{Im}(x_2/x_3)>0$, $\mathrm{Im}(x_3/x_1)>0$, $w\in\mathbb C$:
\[
\Gamma_{e_1,e_2}(w,x)=
\Gamma\left(\frac w{x_3},\frac{x_1}{x_3},\frac{x_2}{x_3}\right).
\]
In this domain $\Gamma$
is given by the convergent
product \eqref{e-expelliatur}
\end{ep}

\begin{remark} One could use Proposition \ref{p-horcrux} to define gamma
functions of wedges. The fact that this definition is independent of choices is a consequence
of the multiplication formulae \cite{FV2}.
\end{remark}

\subsection{Inversion and three-term relations}
The three-term relations \eqref{e-sectumsempra}, \eqref{e-levicorpus}
satisfied by the elliptic gamma function are  special cases
of a collection of three-term relations parametrized by a triple $a,b,c$
of primitive elements of $\Lambda$.

\begin{thm}\label{t-avedakedavra1}
Let $a,b\in \Lambda_{\mathit{prim}}$. Then, for $x\in U_a^+\cap U_b^+$ and $w\in\mathbb C$
we have the
inversion relation
\[
\Gamma_{a,b}(w,x)\Gamma_{b,a}(w,x)=1.
\]
More generally,
let $a,b,c\in\Lambda_{\mathit{prim}}$. Then, for $x\in U_a^+\cap U_b^+\cap U_c^+$ and $w\in \mathbb C$
\[
\Gamma_{a,b}(w,x)\Gamma_{b,c}(w,x)\Gamma_{c,a}(w,x)=\exp\left(-\frac{\pi i}{3} P_{a,b,c}(w,x)\right),
\]
for some polynomial $P_{a,b,c}(w,x)\in \mathbb Q(x)[w]$ of degree at
most 3 in $w$ with rational coefficients.
\end{thm}

Before giving the explicit description of the polynomial $P_{a,b,c}$ we need
to discuss for which triples the domain of validity of the three-term relation
is non-empty.

\begin{lemma}\label{l-muffliato}
Let $a,b,c\in\Lambda_{\mathit{prim}}$. Then $U^+_a\cap U^+_b\cap U^+_c$ is non-empty
if and only if
$a,b,c$ are on the same side of some plane through the origin
\end{lemma}
\begin{proof}
Suppose that $a,b,c$ all lie on the same line. Then either $a=b=c$ and $U^+_a\neq\emptyset$
or the sum of two is zero, in which case the intersection is empty since
$U^+_a$, $U^+_{-a}$ are disjoint.

Suppose that $a,b,c$ lie on a plane. Then the corresponding planes $H(a)$, $H(b)$, $H(c)$
intersect in a line $\mathbb Q\sigma$ and we may choose $\alpha,\beta,\gamma\in V^*$
so that $(\alpha,\sigma)$, $(\beta,\sigma)$, $(\gamma,\sigma)$
are oriented bases of these subspaces.
Now $a,b,c$ are on the same side of a plane if and only if
$\alpha,\beta,\gamma$ have the same property. In this case there is a $t\in V$
such that $\alpha(t),\beta(t),\gamma(t)>0$ and $\sigma(t)=0$. Then for any
$s\in V$ such that $\sigma(s)>0$, we have $x=s+it\in U^+_a\cap U^+_b\cap U^+_c$
so the intersection is non-empty. If $\alpha,\beta,\gamma$ are not on the same
side of a plane, we may rescale them so that $\alpha+\beta+\gamma=0$. Then
for any $x$ the three complex numbers $\alpha(x)\overline{\sigma(x)}$,
$\beta(x)\overline{\sigma(x)}$, $\gamma(x)\overline{\sigma(x)}$ add up to
zero, so they cannot all have positive imaginary part and the intersection
is empty.

Finally, assume that $a,b,c$ are a basis of $V$, which we may assume to
be oriented, and let $\alpha,\beta,\gamma$
be the dual basis of $V^*$. Both bases are always on the same side of some plane
through the origin. Moreover
$(\alpha,\beta)$, $(\beta,\gamma)$, $(\gamma,\alpha)$
are oriented bases of $H(c),H(a),H(b)$ respectively. Then for example
$x\in V_\mathbb C$ so that
$\alpha(x)=1,\beta(x)=(-1-i\sqrt3)/2,\gamma(x)=(-1+i\sqrt3)/2$ belongs
to the intersection $U^+_a\cap U^+_b\cap U^+_c$.
\end{proof}

\begin{thm}\label{t-avedakedavra2}
The polynomial $P_{a,b,c}$ of Theorem \ref{t-avedakedavra1}
can be chosen in the following way:
\begin{enumerate}
\item[(i)]
If $a,b,c$ are linearly dependent then $P_{a,b,c}=0$.
\item[(ii)]
If $a,b,c$ are linearly independent and $\det(a,b,c)>0$,
\begin{eqnarray*}
P_{a,b,c}(w,x)&=&\sum_{\delta\in F(a,b,c)}P_3(w+\delta(x),\alpha(x),\beta(x),\gamma(x))
\\ &=&\frac{\det(\alpha,\beta,\gamma)}
{\alpha(x)\beta(x)\gamma(x)}{w^3}+\mathrm{deg}\leq
2
\end{eqnarray*}
Here $\alpha,\beta,\gamma$ are the direction vectors
of the wedges $(b,c)$, $(c,a)$ and $(a,b)$, respectively, $P_3=B_{3,3}$ is the Bernoulli polynomial \eqref{e-bern} and $F(a,b,c)$
is the set of $\delta\in\Lambda^\vee$ such that
\[
0\leq \delta(a)< \alpha(a),\quad
0\leq \delta(b)< \beta(b),\quad
0\leq \delta(c)< \gamma(c).
\]
\item[(iii)]
For all permutations $\sigma\in S_3$,
\[
P_{a,b,c}(w,x)=\mathrm{sign}(\sigma)P_{\sigma(a),\sigma(b),\sigma(c)}(w,x).
\]
\end{enumerate}
Moreover, $P_{a,b,c}(w,x)$ is $\SL_3(\mathbb Z)$-equivariant:
\[
P_{ga,gb,gc}(w,x)=P_{a,b,c}(w,g^{-1}x),\qquad g\in \SL_3(\mathbb Z).
\]
\end{thm}

\noindent{\it Proof of Theorems \ref{t-avedakedavra1} and \ref{t-avedakedavra2}}
We first prove the inversion relation. The direction vector of
$(a,b)$ is the opposite of the direction vector of $(b,a)$. This implies
that $\Gamma_{b,a}$ is obtained from $\Gamma_{a,b}$ by exchanging the
numerator and the denominator, proving the claim.

Assume that $a,b,c$ lie on the same plane and let $\gamma$ be the direction vector of $(a,b)$. After a
cyclic permutation of $a,b,c$, it may be assumed
that $\gamma$ is also the direction vector of $(b,c)$ and $(a,c)$.
In this situation,
\[
C_{+-}(a,c)=C_{+-}(a,b)\sqcup C_{+-}(b,c),\quad
C_{-+}(a,c)=C_{-+}(a,b)\sqcup C_{-+}(b,c),
\]
implying the identity in the form $\Gamma_{a,b}\Gamma_{b,c}=\Gamma_{a,c}$.

Now let $a,b,c$ be linearly independent and assume that $\det(a,b,c)>0$.
We use the representation of Proposition \ref{p-horcrux}.
Consider first the factor $\Gamma_{a,b}$. The direction
vector $\alpha$ of $(b,c)$ obeys $\alpha(a)>0$ and $\alpha(b)=0$ (since
$\alpha$ is proportional by a positive integer to $\det(b,c,\cdot)$).
Similarly, $\beta(a)=0$ and $\beta(b)>0$. Thus we can write $\Gamma_{a,b}$ as
product of standard elliptic gamma functions as in Proposition \ref{p-horcrux}.
Let $F(a,b)$ denote the set $F$ appearing there.
Since $\gamma(a)=\gamma(b)=0$ and $\gamma(c)>0$, in each class
of $F(a,b)/\mathbb Z\gamma$
there is a unique representative $\delta$, such that $\delta(c)\in\{0,\dots,\gamma(c)-1\}$.
Therefore we can replace
the product in Proposition \ref{p-horcrux} by a product over $F(a,b,c)$.
Applying the same procedure to the remaining gamma functions in the product,
we obtain a product over $F(a,b,c)$ of triple products of standard elliptic gamma
functions. Identity \eqref{e-levicorpus} gives then the result.

The coefficient of $w^3$ in $P_{a,b,c}$ is then the coefficient of $P$
times the cardinality of $F(a,b,c)$. The latter is the volume of the parallelepiped
in $V^*$ generated by the vectors $\alpha,\beta,\gamma$.

It is clear that $P_{a,b,c}$, defined by (i), (ii) for $\det(a,b,c)\geq0$
is invariant under cyclic permutations of $a,b,c$. Applying the inversion
relation to all terms of the three term relation we get $\Gamma_{b,a}
\Gamma_{c,b}\Gamma_{a,c}=\exp(i\pi P_{a,b,c})$. Therefore we may extend the
definition of $P$ to general triples by setting $P_{c,b,a}=-P_{a,b,c}$. With
this definition the three-term relation is obeyed by all triples and (iii)
holds.

The equivariance can be easily checked directly from the definition of $P_{a,b,c}$.
\hfill$\square$

\medskip

\begin{remark}
 By Proposition~\ref{p-horcrux},
the identities of Theorems \ref{t-avedakedavra1} and \ref{t-avedakedavra2} can be written as
$n$-term 
relations for the ordinary elliptic gamma
function \eqref{e-expelliatur}.
The inversion relation for $a=e_1,b=e_2$ reads 
$\Gamma(z,\tau,\sigma)\Gamma(-z,-\sigma,-\tau)=1$.
The three-term relation 
\eqref{e-levicorpus}
is the special case of the second identity in Theorem \ref{t-avedakedavra1}
where $a=e_1,b=e_2,c=e_3$. If $a=e_1,b=e_1-e_2,c=e_2$ we obtain
the relation 
\[
\Gamma\left(-\frac w{x_3},-\frac {x_1+x_2}{x_3},\frac{x_2}{x_3}\right)
\Gamma\left(\frac w{x_3},\frac {x_1}{x_3},\frac{x_1+x_2}{x_3}\right)
\Gamma\left(-\frac w{x_3},-\frac {x_1}{x_3},-\frac{x_2}{x_3}\right)=1,
\]
which reduces to \eqref{e-sectumsempra} if we apply the inversion relation and
\eqref{e-Cleophile}.
For a simple ``new'' identity, consider the case
 $a=e_1$, $b=e_1+2e_2$,
$c=e_3$. Then, by using the direction vectors as
coordinates, we obtain the four-term relation
\begin{gather*}
\Gamma\left(
\frac{w}{y_3},
\frac{y_1}{y_3},
\frac{y_2}{y_3}
\right)
\Gamma\left(
\frac{w}{y_3}+\frac{y_1+y_2}{2\,y_3},
\frac{y_1}{y_3},
\frac{y_2}{y_3}
\right)
\Gamma\left(
\frac{w}{y_1},
\frac{{y_1+y_2}}{2\,y_1},
\frac{y_3}{y_1}
\right)
\Gamma\left(
\frac{w}{y_2},
\frac{y_3}{y_2},
\frac{{y_1+y_2}}{2\,y_2}
\right)
\\
=
\exp\left(-\frac{\pi i}3
\left(
P_3(w,y_1,y_2,y_3)+
P_3(w+{(y_1+y_2)}/2,y_1,y_2,y_3)\right)\right)
\end{gather*}

\dontprint{
{\em Details.}
a) $a=e_1,b=e_1+2e_2$. Direction vector $\gamma=\epsilon_3$.
$\alpha=2\epsilon_1-\epsilon_2$, $\beta=\epsilon_2$. Modulus=2.
\[
\Gamma_{e_1,e_1+2e_2}(w,x)
=
\Gamma\left(\frac{w}{x_3},\frac{2x_1-x_2}{x_3},\frac{x_2}{x_3}\right)
\Gamma\left(\frac{w+x_1}{x_3},\frac{2x_1-x_2}{x_3},\frac{x_2}{x_3}\right).
\]
b) $a=e_1+2e_2, b=e_3$. Direction vector $\gamma=2\epsilon_1-\epsilon_2$.
$\alpha=\epsilon_1,\beta=\epsilon_3$. Modulus=1.
\[
\Gamma_{e_1+2e_2,e_3}(w,x)
=
\Gamma\left(\frac{w}{2x_1-x_2},\frac{x_1}{2x_1-x_2},\frac{x_3}{2x_1-x_2}\right).
\]
c) $a=e_3, b=e_1$. Direction vector $\gamma=\epsilon_2$.
$\alpha=\epsilon_3,\beta=\epsilon_1$. Modulus=1.
\[
\Gamma_{e_3,e_1}(w,x)
=
\Gamma\left(\frac{w}{x_2},\frac{x_3}{x_2},\frac{x_1}{x_2}\right).
\]
Now set $y_1=2x_1-x_2,y_2=x_2,y_3=x_3$, so that $x_1=(y_1+y_2)/2$.
}
\end{remark}

\subsection{Action of $\ISL_3(\mathbb Z)$}\label{ss-Action}
The group $\SL_3(\mathbb Z)$ acts on $\Lambda$, $\Lambda^\vee$, $V$ and $V_\mathbb C$.
It acts transitively on $\Lambda_{\mathit{prim}}$.
If $a\in \Lambda_{\mathit{prim}}$ and $g\in \SL_3(\mathbb Z)$, $U^+_{ga}=gU^+_a$ and
\begin{equation}\label{e-OmicronPi}
\Gamma_{ga,gb}(w,x)=\Gamma_{a,b}(w,g^{-1}x), \qquad x\in U^+_a\cap U^+_b.
\end{equation}
The group $\ISL_3(\mathbb Z)=\SL_3(\mathbb Z)\ltimes \mathbb Z^3$
acts on $\mathbb C\times V_\mathbb C$ by
\[
(g,\alpha)\cdot(w,x)=(w-\alpha(x),g\cdot x),\quad
g\in \SL_3(\mathbb Z),\quad \alpha\in\Lambda^\vee=\mathbb Z^3.
\]
The subgroup $\Lambda^\vee$ maps $U^+_a$ to itself for all $a\in\Lambda$.
Here is a description of the action on gamma functions.

\begin{definition}\label{d-severus} A {\em framing} of
$\Lambda_{\mathit{prim}}$ assigns to each $a\in\Lambda_{\mathit{prim}}$
an oriented basis $f(a)=(\alpha_1,\alpha_2,\alpha_3)$ of the $\mathbb Z$-module
$\Lambda^\vee$ such that $\alpha_1(a)=1$ and $\alpha_2,\alpha_3\in H(a)$.
\end{definition}

For each framing $f$ we define a family of products of theta functions parametrized
by $a\in\Lambda_{\mathit{prim}}$ and $\mu\in \Lambda^\vee$:
\begin{equation}\label{e-DrStanhope}
\Delta_a(\mu;w,x)=\Delta_{a,f}(\mu;w,x)=\prod_{j=0}^{\mu(a)-1}
\theta_0
\left(
 \frac{w+j\alpha_1(x)}{\alpha_3(x)},\frac{\alpha_2(x)}{\alpha_3(x)}
\right).
\end{equation}
Here the product should be understood as $\prod_{-N\leq j<\mu(a)}/\prod_{-N\leq j<0}$
for any sufficiently large $N$ so it is defined also for negative $\mu(a)$.
Notice moreover that $\Delta_a(\mu;w,x)=1$, if $\mu$ belongs to the plane
$H(a)$.
\begin{proposition}\label{p-Lazare}
Let $a\in\Lambda_{\mathit{prim}}$, $\mu\in\Lambda^\vee$.
The function $(w,x)\to\Delta_a(\mu;w,x)$ is
meromorphic on $\mathbb C\times U_a^+$.
\end{proposition}

\begin{proof} The theta function $\theta(z,\tau)$ is holomorphic for $\mathrm{Im}\,\tau> 0$. The function $\Delta_a$, being a product of (inverses of) theta functions
with $\tau=\alpha_2(x)/\alpha_3(x)$ is meromorphic for $\mathrm{Im}(\alpha_2(x)/\alpha_3(x))>0$,
i.e., for $x\in U_a$.
\end{proof}

\begin{proposition}\label{p-felixfelicis}
Let $a,b\in\Lambda_{\mathit{prim}}$ be linearly independent and fix a framing of $\Lambda_{\mathit{prim}}$. Then
\[
\frac{\Gamma_{a,b}(w+\mu(x),x)}{\Gamma_{a,b}(w,x)}=
e^{\pi i P_{a,b}(\mu;x,w)}\frac{\Delta_a(\mu;w,x)}{\Delta_b(\mu;w,x)}
\]
for some (framing dependent) polynomial $P_{a,b}(\mu;x,w)\in \mathbb Q(x)[w]$ homogeneous of
degree $0$ in $x,w$ and of degree at most 2 in $w$.
\end{proposition}

Note that a change of framing multiplies $\Delta_a$ by the exponential of a polynomial
of degree at most 2 in $w$. Indeed a change of basis $\alpha_2,\alpha_3$ of $H(a)$
results in an $\SL_2(\mathbb Z)$ transformation of the arguments of $\theta_0$, producing
a multiplier given by the exponential of a
polynomial of degree $\leq 2$. Also, $\alpha_1$ is determined up to the addition
of an element of $H(a)\cap\Lambda^\vee$. Therefore, under a change of $\alpha_1$,
the theta functions in $\Delta_a$ change by a shift of $w$ by an integer linear
combination of $\alpha_2$ and $\alpha_3$. Under such shifts $\theta_0$ changes
by the exponential of a polynomial of degree $\leq 1$.

\medskip

\noindent{\it Proof of Proposition \ref{p-felixfelicis}\/.}
We consider exemplarily the case $\mu(a)>0$, $\mu(b)>0$; the other cases are treated in a similar way, whenever neither $\mu(a)$
nor $\mu(b)$ vanishes, with due modifications.
The cases where at least one of $\mu(a)$ and $\mu(b)$ vanishes are better treated separately, but the main computations
are similar to the following ones.

In the ratio of gamma
functions infinitely many terms cancel in the ratio of infinite products.
The remaining infinitely many terms can be collected as follows.
Let $\gamma$ be the direction vector of
$(a,b)$ and let $\alpha,\beta\in\Lambda_{\mathit{prim}}^\vee$
be such that $\alpha(b)=\beta(a)=0$, $\alpha(a),\beta(b)>0$.
Then
\[
\frac{\Gamma_{a,b}(w+\mu(x),x)}{\Gamma_{a,b}(w,x)}=
\prod_{\delta\in F_1}\left(- e^{-2\pi i(w+\delta(x))/\gamma(x)}\right)
\frac
{
\prod_{\delta\in F_2}\theta_0\left(\frac {w+\delta(x)}{\gamma(x)},\frac{\beta(x)}{\gamma(x)}
\right)
}
{
\prod_{\delta\in F_3}\theta_0\left(\frac {w+\delta(x)}{-\gamma(x)},\frac{\alpha(x)}{-\gamma(x)}
\right)
}\,.
\]
In these products $\delta$ runs over the finite subsets
of $\Lambda^\vee/\mathbb Z\gamma$ defined by the following inequalities:

\begin{gather*}
F_1\colon 0\leq\delta(a)<\mu(a),\quad 0\leq\delta(b)<\mu(b),
\\
F_2\colon 0\leq\delta(a)<\mu(a),\quad 0\leq\delta(b)<\beta(b),
\\
F_3\colon 0\leq\delta(a)<\alpha(a),\quad 0\leq\delta(b)<\mu(b).
\end{gather*}
We need to show that the right-hand side is $\Delta_a/\Delta_b$ up
to multiplication by an exponential functions of a polynomial.

Consider e.g.\ the product indexed by $F_2$: we notice that for each $0\leq j<\mu(a)$
there is exactly one $\delta\in\Lambda^\vee/\gamma\mathbb Z$
such that  $\delta(a)=j$ and $0\leq\delta(b)<\beta(b)$. This can be seen
most easily by going to the normal form, where $a=e_1$, $b=re_1+se_2$, $s>0,
0\leq r<s$, gcd$(r,s)=1$, see subsection \ref{s-wedges}: if $\epsilon_1,\epsilon_2,
\epsilon_3$ is the dual basis of $\Lambda^\vee$, $\beta=\epsilon_2$,
$\alpha=s\epsilon_1-r\epsilon_2$, $\gamma=\epsilon_3$ and the unique $\delta$ with
these properties is $\delta=j\epsilon_1+k\epsilon_2\mod\epsilon_3$ where
$k$ is the unique integer such that $0\leq jr+ks<s$.

If we are given a framing and $(\alpha_1,\alpha_2,\alpha_3)$ is the basis assigned to $a$,
this $\delta$ differs from $j\alpha_1$ by a lattice vector in $H(a)$, i.e., an
integer linear combination of $\beta$ and $\gamma$.
Thus we may replace $\delta$ by $j\alpha_1$ in the numerator up to the exponential of a polynomial.
A similar reasoning applies to the denominator.

The oriented bases $(\beta,\gamma)$ of $H(a)$ and $(\alpha,-\gamma)$ of $H(b)$
can then be replaced by the bases given by the framing up to the exponential
of a quadratic polynomial.

This proves the claim. \hfill$\square$

\begin{proposition}
\label{p-EleanorBold}
 Fix a framing of $\Lambda_{\mathit{prim}}$. Let $\mu,\nu\in\Lambda^\vee$, $a\in\Lambda_{\mathit{prim}}$.
Let $\mu=\sum m_i\alpha_i$, $\nu=\sum n_i\alpha_i$ be the
decomposition of $\mu,\nu$ in terms of the basis assigned to
$a$ by the framing. Then
\[
\Delta_a(\mu+\nu;w,x)=e^{2\pi i P_a(\mu,\nu;w,x)}
\Delta_a(\mu;w,x)\Delta_a(\nu;w+\mu(x),x)\]
where
\[
P_a(\mu,\nu;w,x)=\frac{n_1m_2}{\alpha_3(x)}
\left[w+
\frac{2m_1+n_1-1}{2}\alpha_1(x)
+\frac{m_2-1}{2}\alpha_2(x)+\frac{1}{2}\alpha_3(x)\right].
\]
\end{proposition}

\section{The gamma gerbe}\label{s-iv}

In this section we cast together the identities obtained in the previous sections and describe them
geometrically in the gerbe language. After fixing some sign conventions,
we begin by reviewing the well-known case of the Jacobi theta function,
based on the classical identities discovered by Jacobi and Dedekind.

\subsection{Conventions and notations}
Our convention for the differential $\check\delta$ of the \v Cech complex
$\check C(\mathcal U,\mathcal F)$ of a sheaf $F$
with open covering $\mathcal U=(V_a)_{a\in I}$ is
\[
\check\delta\phi_{a_0,\cdots,a_{p+1}}=\phi_{a_1,\dots,a_{p+1}}-\phi_{a_0,a_2,\dots,
a_{p+1}}+\cdots+(-1)^{p+1}\phi_{a_0,\dots,a_p},
\]
if $\phi=\oplus\phi_{a_0,\dots,a_p}\in \check C^p(\mathcal U,\mathcal F)=\oplus_{a_0,\cdots,a_p}\mathcal F(V_{a_0}\cap\cdots\cap V_{a_p})$. The differential $\delta$ of
the complex $C(G,M)=\oplus_p\mathrm{Maps}(G^p,M)$
of cochains of the group $G$ with values in the left $G$-module
$M$ is
\begin{eqnarray*}
\delta c(g_1,\dots,g_{p+1})&=&(-1)^{p+1}g_1\cdot c(g_2,\dots,g_{p+1})\\
&&+\sum_{j=1}^p(-1)^{p+1-j}
c(g_1,\dots,g_jg_{j+1},\dots,g_{p+1})\\
&&+c(g_1,\dots,g_p).
\end{eqnarray*}
In our case $G$ is a discrete group acting on a complex manifold and $M$ is
the \v Cech complex of an equivariant sheaf $\mathcal F$
corresponding to an invariant open covering.
Then the double complex $C_G(\mathcal U,\mathcal F)
=\oplus_{p,q}C^p(G,\check C^q(\mathcal U,\mathcal F))$ with total differential $D=\delta+(-1)^p\check \delta$ computes the equivariant cohomology $H_G(X,\mathcal F)$.
Many of the sheaves $\cF$ we consider, such as $\cO^\times$, $\cM^\times$,
are sheaves of multiplicative groups, so that the formulae for the differentials are written multiplicatively.

\subsection{The theta line bundle on $X_2$}
 Let $X_2$ be the total space of the dual tautological line bundle $O(1)\to(\mathbb C P^1-\mathbb R P^1)$. Thus $X_2$ is the quotient of $\mathbb C\times (\mathbb C^2-\{0\})$ by the action
 $(w,x)\mapsto(\lambda w,\lambda x)$ of $\mathbb C^\times$.
The group $G=\ISL_2(\mathbb Z)=SL_2(\mathbb Z)\ltimes \mathbb Z$ acts on $X_2$ by $(g,\mu)\cdot(w,x)=(w-\mu_1x_1-\mu_2x_2,g\cdot x)$.
As a complex manifold,
$X_2$ is isomorphic to $\mathbb C\times(\mathbb C-\mathbb R)$. It is covered by two invariant
contractible disjoint open sets $V_\pm=\{(w,x)\in\mathbb C\times (\mathbb C^2-\{0\})\,|\, \pm\mathrm{Im}(x_1/x_2)>0\}/\mathbb C^\times$. Let
\[
\theta_\pm(w,x)=\theta_0\left(\pm\frac w{x_2},\pm\frac{x_1}{x_2}\right)^{\pm1}, \qquad (w,x)\in V_\pm.
\]
Then, for every $h\in G$, $(w,x)\in V_\pm$,
\[
\theta_\pm(w,x)=\phi_\pm(g;w,x)\theta_\pm(h^{-1}(w,x)).
\]
The group 1-cocycle $\phi_\pm$ is given by the following formulae:
let $N:\SL_2(\mathbb Z)\to \mathbb Z/12\mathbb Z$ be the homomorphism whose value on generators
is
\[
N {1\;\phantom{-}1\choose 0\;\phantom{-}1}\equiv1\mod 12,\qquad
N {0\;-1\choose1\;\phantom{-}0}\equiv-3\mod 12.
\]
Then
if $g\in \SL_2(\mathbb Z)$
\[
\phi_\pm(g;w,x)=\exp\left(-\pi i Q\left(g;\frac w{x_2},\frac{x_1}{x_2}\right)\right),
\]
where for $g^{-1}={a\;b\choose c\;d}$,
\[
Q(g;z,\tau)=-\frac {c z^2}{c\tau+d}-\frac z{c \tau+d}+\frac z+\frac{1}{6}\frac{a\tau+b}{c\tau+d}-
\frac{\tau}6+\frac{1}6 N(g)
; 
\]
if $\mu\in\mathbb Z^2$,
\[
\phi_\pm(\mu;w,x)=\exp\left(2\pi i\mu_1\frac w{x_2}+\pi i\mu_1(\mu_1-1)\frac{x_1}{x_2}+\pi i\mu_1\right).
\]
This 1-cocycle extends in a trivial way to a  1-cocycle
\[(1,\phi_\pm)\in C^{0,1}_G(\mathcal{U},\mathcal O^\times)\oplus C^{1,0}_G
(\mathcal U,\mathcal O^\times),\] in the group-\v Cech
double complex for the covering $\mathcal U=\{V_+,V_-\}$
and defines a $G$-equivariant line bundle on $X_2$. The functions $\theta_{\pm}$ form
then a meromorphic equivariant section in $C^{0,0}_G(\mathcal U,\mathcal M^\times)$.
Clearly the \v Cech part of this story is trivial since $X_2$ is just the disjoint union of two
contractible sets, but this will change when we pass to $X_3$.

\subsection{The gamma gerbe on $X_3$}\label{sec:the gamma gerbe}
Let $X=X_3$ be the total space of the
line bundle $O(1)\to(\CP^2-\RP^2)$.
Geometrically we think of $O(1)$ as the dual bundle to the tautological line bundle and of
$\CP^2$ as the projectivization of $\Lambda_{\mathbb C}=\Lambda\otimes_\mathbb Z\mathbb C$, where $\Lambda=\mathbb Z^3$ is a free abelian group of rank 3 equipped
with a volume form $\det:\wedge^3\Lambda\to\mathbb Z$. The group
$\mathrm{Aut}(\Lambda)\cong \SL_3(\mathbb Z)$
of linear transformations of $\Lambda_\mathbb C$ mapping $\Lambda$ to itself and preserving the volume
form, acts naturally on $X$. The dual lattice $\Lambda^\vee=\mathrm{Hom}(\Lambda,\mathbb Z)
\cong \mathbb Z^3\subset V_\mathbb C^*$ acts on $O(1)$ fiberwise by translation and
we get an action of $\mathrm{Aut}(\Lambda)\ltimes\Lambda^\vee\cong \ISL_3(\mathbb Z)$.
More explicitly, this group acts linearly on $\mathbb C\times V_{\mathbb C}$
 via $(g,\mu)(w,x)=(w-\alpha(x),gx)$, and this action induces an action on $X=(\mathbb C\times (V_{\mathbb C}-\mathbb C\cdot V_\mathbb R))/
\mathbb C^\times$.

The complex manifold $X$ has a natural $\ISL_3(\mathbb Z)$-invariant open covering
$\mathcal U=(V_a)_{a\in\Lambda_{\mathit{prim}}}$ labeled
by the set of primitive vectors $\Lambda_{\mathit{prim}}$ in $\Lambda$: for $a\in\Lambda_{\mathit{prim}}$, $V_a=
\{ (w,x)\,|\, x\in U_a^+\}/\mathbb C^\times$, see subsection \ref{ss-Domain}.%
\footnote{One can show that this covering is good, i.e., all non-empty multiple intersections are contractible. Moreover the sets $V_a$ are domains of holomorphy. It
follows that the double complex $C(G,\check C(\mathcal U,\mathcal F))$
can be used to compute the equivariant cohomology with values in
equivariant constant or analytic coherent sheaves.}

Recall that a {\em framing} of $\Lambda_{\mathit{prim}}$ is a map $f:\Lambda_{\mathit{prim}}\to
\{$oriented bases of $\Lambda^\vee\}$ such the basis
 $f(a)=(\alpha_1,\alpha_2,\alpha_3)$ obeys $\alpha_1(a)=1$ and
$\alpha_2(a)=\alpha_3(a)=0$.

Let us fix a framing $f$. Then we have the gamma functions $\Gamma_{a,b}
\in\mathcal M^\times(V_a\cap V_b)$, associated to wedges $(a,b)\in\Lambda_{\mathit{prim}}^2$, and the products of theta functions $\Delta_a=\Delta_{a,f}\colon \mathbb Z^3\to
\mathcal M^\times(V_a)$, see Section \ref{s-iii}, Eqns.~\eqref{e-imperius}, \eqref{e-DrStanhope}.
Extend $\Delta_{a}$ to a map
$\ISL_3(\mathbb Z)\to\mathcal M^\times(V_a)$ by setting
$\Delta_{a}((g,\mu);w,x)=\Delta_{a}(\mu\circ g^{-1};w,x)$ for $g\in \SL_3(\mathbb Z)$, $\mu\in\Lambda^\vee=\mathbb Z^3$.

\begin{thm}\label{t-main}
 Let $G=\ISL_3(\mathbb Z)$ and $f$ be a framing of $\Lambda_{\mathit{prim}}$. View $\Gamma_{a,b}$, $\Delta_a$ as the components of a cochain
$
s=(\Gamma,\Delta)\in C^1_G(\mathcal U,\mathcal M^\times).
$
Then the functions forming the cocycle $\phi=D s$ are holomorphic and nowhere vanishing.
Thus $\phi\in C_G^2(\mathcal U,\mathcal O^\times)$ defines a holomorphic gerbe
$\mathcal G$ on $X_3$ and $s$ is a meromorphic section of this gerbe.
\end{thm}

This theorem is a rephrasing Theorem \ref{t-Quiverful}: the
identity $D s=\phi$, written explicitly, is the
set of equations appearing there.

\begin{proof} The first thing to check is that $s=(\Gamma,\Delta)\in C^1(\mathcal U,\mathcal M^\times)$,
namely that these functions are meromorphic on their domain of definition.
This is the content of
Propositions \ref{p-Goodenough}, \ref{p-Lazare}.

Now define $\phi$ by  $\phi=Ds$, namely by the identities of Theorem \ref{t-Quiverful}.
{\em A priori,} $\phi$ is a 2-cocycle with values in the sheaf $\mathcal M^\times$ of
invertible meromorphic functions and we
have to show that it takes value in the invertible {\em holomorphic} functions. It is sufficient
to show that it takes values in the exponentials of rational functions, since a meromorphic
function of this form is automatically holomorphic and nowhere vanishing.
The first identity in Theorem \ref{t-Quiverful} is treated in Theorems \ref{t-avedakedavra1}, \ref{t-avedakedavra2}, implying that $\phi_{a,b,c}$ is indeed holomorphic and invertible on
$V_a\cap V_b\cap V_c$. Proposition \ref{p-felixfelicis} implies that  $\phi_{a,b}(g;y)$
is the exponential of a rational function for $g$ in the translation subgroup of $G$. The
second identity of Theorem \ref{t-Quiverful} for general $g$ can be obtained from
Proposition \ref{p-felixfelicis} using the equivariance of $\Gamma$ under $\SL_3(\mathbb Z)$,
Equation \ref{e-OmicronPi}. Again the conclusion is that $\phi_{a,b}$ is the exponential of
a rational function. As for the third identity of Theorem \ref{t-Quiverful}, it appears in
Proposition \ref{p-EleanorBold} for $\mu,\nu$ in the translation subgroup.
Since $\Delta_{a,f}(g;w,x)$ for general $g$ is defined in terms of its restriction to
the translation subgroup, the identity of Proposition \ref{p-EleanorBold} gives an identity
for general group elements. However to reduce it to the third identity of Theorem \ref{t-Quiverful}
we need to relate $\Delta_{g^{-1}a,f}(\mu\circ g;w,g^{-1}x)$ to $\Delta_{a,f}(\mu;w,x)$ for $g\in \SL_3(\mathbb Z)$.
The relation involves the natural action of  $\SL_3(\mathbb Z)$ on the set of framings:
\[
\Delta_{g^{-1}a,f}(\mu\circ g;w,g^{-1}x)=\Delta_{a,gf}(\mu;w,x), \quad g\in \SL_3(\mathbb Z),\quad
\mu\in\Lambda^\vee.
\]
As discussed after Proposition \ref{p-felixfelicis}, a change of framing amounts to an
$\ISL_2(\mathbb Z)$-transformation of the theta functions and so we have
\begin{equation}\label{e-framing}
\Delta_{a,gf}(\mu;w,x)=\psi_{gf,f}(\mu;w,x)\Delta_{a,f}(\mu;w,x),\qquad (w,x)\in V_a,
\end{equation}
for some exponential of rational function $\psi$, which, by the same argument as above,
is an invertible holomorphic function on $V_a$ which enters in the third identity of
Theorem \ref{t-Quiverful} for general group elements.
 \end{proof}

\begin{remark} Equation \eqref{e-framing} in the proof shows that a change of framing
from $f$ to $f'$ amounts to replacing $\Delta$, by $\psi_{f',f}\Delta$ where $\psi_{f',f}\in C^{1,0}(\mathcal U,\mathcal O^\times)$ and thus $\phi$ by  a cocycle differing from $\phi$ by the exact
cocycle $D(1,\psi_{f',f})$.
\end{remark}

By construction, the components of the cocycle are exponentials of rational
functions. Their restriction of the components of the cocycle $\phi$ to the subgroup
$\Lambda^\vee=\mathbb Z^3$ is given by the explicit formulae calculated
in Section \ref{s-iii}:
for $a,b,c\in\Lambda_{\mathit{prim}}$, $\mu,\nu\in\Lambda^\vee\subset G$,
\begin{gather*}
\phi_{a,b,c}(w,x)=e^{-\frac{2\pi i}{3!}P_{a,b,c}(w,x)},
\qquad\phi_{a,b}(\mu;w,x)=e^{-\frac{2\pi i}{2!}P_{a,b}(\mu;w,x)},
\\
\phi_a(\mu,\nu;w,x)=e^{-2\pi i P_a(\mu,\nu;w,x)},
\end{gather*}
where $P_{a,b,c}(w,x),P_{a,b}(\mu;w,x),P_a(\mu,\nu;w,x)\in
\mathbb Q(x)[w]$ are the homogeneous rational functions
defined in Theorem \ref{t-avedakedavra2}, Proposition \ref{p-felixfelicis}
and Proposition \ref{p-EleanorBold}, respectively.

\begin{definition}\label{def: gamma gerbe}
The {\em gamma gerbe}
is the equivariant holomorphic gerbe $\mathcal G$ on $X_3$ defined by
the cocycle $\phi$.
\end{definition}

\subsection{The gamma gerbe as a groupoid central extension}
Now we also state our main result in the groupoid language
explained in the Appendix.

We take the invariant covering $V_a$ of $X=X_3$ with $a\in \Lambda_{\mathit{prim}}$. The
stack $\cX$$:=$$[X$ $/\ISL_3(\Z)]$ is presented by the action groupoid
$\ISL_3(\ZZ) \times X \rra X$.  Let $U_0=\sqcup V_a$ and
$U_1=U_0\times_{X}\times (\ISL_3(\ZZ) \times X) \times_{X} U_0$.
Then $U_1\rra U_0$ is a Lie groupoid that also presents $\cX$.
Notice that $U_1$ breaks into pieces: $U_1=\sqcup_{g, a, g^{-1}b}
W_{g, a, g^{-1}b}$ with $W_{g, a, g^{-1}b}=\{ (g, y), y\in V_a,
g^{-1} y \in V_{g^{-1}b}\}$. 
Then the gamma gerbe $\cG$
is a $\C^\times$-gerbe on the stack
$\cX$. Moreover, the elliptic gamma functions provide a
meromorphic section of that gerbe.

\begin{thm} \label{thm:gamma-gpd} The gamma gerbe $\cG$
corresponds to a central extension of groupoids,
\[
\xymatrix{
1 \ar[r] &\C^\times\times U_0  \ar[d] \ar@<-1ex>[d] \ar[r] & R \ar[d] \ar@<-1ex>[d] \ar[r] & U_1 \ar@<-1ex>[l]_{\Gamma_{a,b}\Delta_b^{-1}} \ar[d] \ar@<-1ex>[d] \ar[r] & 1 \\
& U_0 \ar[r]^{id} & U_0 \ar[r]^{id} &  U_0 }
\]
where $R$ is $L^\times_{a,b}\otimes L^\times_{b}(g)^{-1}$ on $W_{g, a,
g^{-1}b}$, and $\square^\times$ denotes the associated $\C^\times$-bundle 
of a line bundle.
The maps
$\Gamma_{a,b}\Delta_b^{-1}$ form a meromorphic groupoid homomorphism
from $U_1\rra U_0$ to $R\rra U_0$.
\end{thm}
\begin{proof}
As shown in Proposition
\ref{prop:relate}, $R$ is indeed a central extension of $U_1\rra
U_0$ and $R$ presents a $\C^\times$-gerbe $\cG$
corresponding to the cohomology class in $H^2(\cX, \cO^\times)$ represented by $(\phi_{a,b,c}, \phi_{a,b},
\phi_a)$. To see that
$\Gamma_{a,b}\Delta_b^{-1}$'s provide a meromorphic section of
$\cG$ in the above sense, we only have to notice
\[
\begin{split}
& (g, y, \Gamma_{a,b}(y)\Delta_b^{-1}(g; y) ) \cdot
 (h, g^{-1}y, \Gamma_{g^{-1}b, g^{-1} c}(g^{-1} y) \Delta_{g^{-1}c}^{-1}(h;
g^{-1} y) ) \\
=& (gh, y, \Gamma_{a,b}(y)\Delta_b^{-1}(g; y) \Gamma_{g^{-1}b,
g^{-1} c}(g^{-1} y) \Delta_{g^{-1}c}^{-1}(h; g^{-1}
y) \\
 &  \phi_{a,b,c}^{-1}(y) \phi_{b,c}(g;y)\phi_{c}(g,h;y)) \\
= & (gh, y, \frac{\Gamma_{a,b}(y)
\Gamma_{b,c}(y)}{\Gamma_{a,c}(y)} \frac{\Gamma_{g^{-1}b, g^{-1}
c}(g^{-1} y)\Delta_c (g, y)}{\Gamma_{b,c}(y) \Delta_{b}(g; y)}
\frac{\Delta_c(gh;
y)}{\Delta_{c}(g;y) \Delta_{g^{-1} c}(h; g^{-1}y)}\\
&\Gamma_{a,c}(y)\Delta_c^{-1}(gh; y) \phi_{a,b,c}^{-1}(y) \phi_{b,c}(g;y)\phi_{c}(g,h;y))\\
= &(gh, y, \Gamma_{a,c}(y)\Delta_c^{-1}(gh; y)).
\end{split}
\]  with the multiplication $\cdot$ on $R$ given
in \eqref{eq:mul-cent} and the identities
given in Theorem \ref{t-Quiverful}.
\end{proof}

Let $Y=\CP^2$$-$$\RP^2$. The fibres of the projection $[X/\ISL_3(\ZZ)]\to [Y/\SL_3(\ZZ)]$ have the
form $[\C/\iota(\ZZ^3)]$ for some linear map $\iota\colon \ZZ^3\to\C$ of rank 2 over $\R$.
When we restrict the gerbe to a fibre $\cE$$=$$[\C/\iota(\ZZ^3)]$, there is a nice groupoid presentation of the gerbe. We
may assume that $\iota(n_1,n_2,n_3)$=$n_1\tau+n_2\sigma+n_3$, with $\tau,\sigma\in\C$ and
$\mathrm{Im}\,\sigma\neq 0$. We regard $\cE$ as
the stack $[E_\sigma/\ZZ]$, where
$E_\sigma$ =$\C/(\ZZ+\sigma\ZZ)$ and $\ZZ$ acts
by shifts by $[\tau]$, the class of $\tau$ in
$E_\sigma$.
 Denote by $L([z])$ a
line bundle on $E_\sigma$ with Chern class 1 corresponding to $[z]\in E_\sigma$,
i.e., to the element $([z],1)$ of the Picard group $H^1(E_\sigma, \cO^\times)=E_\sigma \times \ZZ$. The
theta function $\theta_0(z, \sigma)$ is a \dontprint{famous} section of
$L([0])$, in the sense that it provides an equivariant section of
the trivial line bundle
$\C\times \C$ on $\C$ with an $\ZZ + \sigma \ZZ$ action. In fact $L([0])$ is the restriction of the
theta bundle on a fibre $E_\sigma$. Let
\begin{gather*}
L^{(n)}:=L([0])^*\otimes L([-\tau])^*\otimes...\otimes
L([-(n-1)\tau])^*, \quad n\geq 0, \\
L^{(n)}:=L([n\tau])\otimes ... \otimes L([-\tau]), \quad n<0,
\end{gather*}
and note that $L^{(n)}\otimes L^{(m)}=L^{(n+m)}$.
Then $\prod_{j=0}^{n-1} \theta_0^{-1}(z+j\tau, \sigma)$ is a
section of $L^{(n)}$ (where the product for $n-1<0$ is understood
as $\prod_{j=-N}^{n-1}/\prod_{j=-N}^0$). The disjoint union $\sqcup_n (L^{(n)})^\times$ of the associated $\C^\times$-bundles is a $\C^\times$-bundle over $E_\sigma \times \ZZ$.
Moreover, it is a groupoid over $E_\sigma$ with source and target maps
the composition of the projection $\sqcup_n L^{(n)}$ $\to$ $E_\sigma \times
\ZZ$ and the source and target maps of $E_\sigma \times
\ZZ\rra E_\sigma$. The multiplication is given by the tensor product
of line bundles
\[
L^{(n)} \times_{E_\sigma} L^{(m)}\to L^{(n)}\otimes L^{(m)}=L^{(n+m)}
\]
since $(\cdot (-n))^* L([k\tau])^* = L([(k+n)\tau])^*$. Here we
only demonstrate this for $n, m\geq 0$. The other  cases are similar.
Then $\sqcup_n L^{(n)}\rra E_\sigma$ is a $\C^\times$-central extension of $E_\sigma \times \ZZ \rra E_\sigma$ and
the section $\prod_{j=0}^{n-1} \theta_0^{-1}(z+j\tau,
\sigma)$ is a groupoid morphism.

\begin{proposition} The gamma gerbe $\cG$ restricted to $[\C/\iota(\ZZ^3)]$
is presented by $\sqcup_n L^{(n)}\rra E_\sigma$
and fits in the central extension of
groupoids,
\[
\xymatrix{
1 \ar[r] &\C^\times\times E_\sigma  \ar[d] \ar@<-1ex>[d] \ar[r] & \sqcup_n (L^{(n)})^\times \ar[d] \ar@<-1ex>[d] \ar[r] & E_\sigma \times \ZZ \ar@<-1ex>[l]_{\prod \theta_0^{-1}} \ar[d] \ar@<-1ex>[d] \ar[r] & 1 \\
& E_\sigma \ar[r]^{id} & E_\sigma \ar[r]^{id} &  E_\sigma, }
\]
where $\square^\times$ denotes the associated $\C^\times$-bundle 
of a line bundle.
\end{proposition}
\begin{proof}
A fibre at a point
$[(w, x_1, x_2, x_3)]$ $\in$ $\cX_3$ is $[\C/\iota(\ZZ^3)]$.
By Theorem \ref{thm:gamma-gpd}, the restriction of the gamma gerbe on this fibre can be
presented by $R|_{\C} \rra \C$, which is a central extension of $\C\times \ZZ^3\rra \C$. Suppose that
$[(w, x_1, x_2, x_3)] \in V_{e_1}$, then we can choose $\C$ to sit inside
$V_{e_1}$. Denote $[(w, x_1, x_2, x_3)]=[(z, \tau, \sigma, 1)]$. From the same theorem, $\Delta_{e_1}^{-1}$ provides a meromorphic section
of this central extension. To show the result, we only have to show that $R|_{\C}\rra \C$ is
Morita equivalent to $\sqcup_n (L^{(n)})^\times \rra E_\sigma$.  Since there is a surjective submersion
$\C \to E_\sigma$ by $z \mapsto [z]$, we only have to show that the
pull-back groupoid $\sqcup_n (L^{(n)})^\times \times_{E_\sigma\times E_\sigma}
\C \times \C$ is the same as $R|_{\C} \rra \C$. But what is
obvious is that the pull-back groupoid of $E_\sigma \times \ZZ \times_{E_\sigma\times E_\sigma}
\C \times \C$ is the action groupoid $\C \times \ZZ^2 \times \ZZ\rra
\C$, and $\C \times \ZZ^2 \times \ZZ$ maps to $E_\sigma\times \ZZ$ by
$(z, n_1, n_2, n_3)\mapsto ([z+n_3+n_2\sigma], n_1\tau)$. The
pull-back groupoid is again a $\C^\times$-bundle on $\C\times \ZZ^2 \times
\ZZ$, as the following diagram of pull-back groupoids shows:
\[
\xymatrix{
\sqcup_n(L^{(n)})^\times \ar[ddr] \ar@<-1ex>[ddr] \ar[dr] & &\sqcup_n (L^{(n)})^\times \times_{E_\sigma\times E_\sigma}
\C \times \C \ar[ll] \ar[dr] \ar[ddr]\ar@<-1ex>[ddr] &\\
& E_\sigma \times \ZZ  \ar@<-1ex>[ul]_{\prod\theta_0^{-1}} \ar[d]\ar@<-1ex>[d] & & \C\times \Z^2 \times \Z \ar[ll] \ar@<-1ex>[ul]_{\Delta_{e_1}^{-1}} \ar[d]\ar@<-1ex>[d]\\
& E_\sigma & & \C \ar[ll]}
\]
Therefore the pull-back groupoid is made up by disjoint union of
$\ZZ^2$ equivariant
$\C^\times$-bundles on $(\C\times \ZZ^2)\times \ZZ$ (hence they are all trivial). So the pull-back groupoid is the same as $R|_{\C}$ as a manifold. By
\eqref{e-DrStanhope}, the
section $\prod_{j=0}^{n-1} \theta_0^{-1}$ pull backs
to $\Delta_{e_1}^{-1}$. Since $\C^\times$ is one-dimensional, the multiplication on the $\C^\times$-central extension is determined by
the multiplication of one section. Since  $\prod_{j=0}^{n-1} \theta_0^{-1}$
is a groupoid morphism,  $\Delta_{e_1}^{-1}$ is a groupoid homomorphism for the pull-back groupoid, but it is also
a groupoid morphism for $R|_{\C}$. Therefore, the pull-back
groupoid is $R|_{\C}$.
\end{proof}

\begin{remark}
This direct fibrewise construction of the gamma gerbe via the
theta bundle provides a possibility to  build
higher gamma (n-) gerbes by the similar procedure.
\end{remark}

\section{Constructing gerbes via divisors}\label{s-pseudodivisors}

The purpose of this section is to give another, more geometrical construction of the gamma gerbe.

First of all, instead of talking about $\ISL_3(\Z)$-equivariant gerbes on $X=X_3$, we shall talk about gerbes on the quotient stack
$\cX=\cX_3:=[X/\ISL_3(\Z)]$.
A holomorphic gerbe on $\cX$ will then be a gadget (see \cite{Bry93} for a detailed discussion) that assigns to each \'etale map $U\to \cX$ a category
$\cG_U$ which, for $U$ small enough, is equivalent to the category of holomorphic line bundles on $U$.
As an example, the trivial gerbe is given by $$U\mapsto\{\text{holomorphic line bundles on }U\}.$$
If $\cL_1$, $\cL_2$ are two objects of $\cG_U$, then $\Hom(\cL_1,\cL_2)$ can be identified with the
set of sections of a line bundle over $U$.
We shall denote this line bundle by $\Homm(\cL_1,\cL_2)$.


We begin by an intrinsic description of $\cX$.

\begin{definition}\label{def triptic}
A {\em triptic curve} $\cE$ is a pointed holomorphic stack of the form $[\C/\iota(\Z^3)]$, where
$\iota:\Z^3\to\C$ is a map of (real) rank $2$.

An {\em orientation} of a triptic curve $\cE$ is given by the choice of a generator of $H^3(\cE,\Z)\simeq\Z$.
\end{definition}

\begin{proposition}
The stack $\Tr:=[(\C P^2-\R P^2)/\SL_3(\Z)]$ is the moduli space of oriented triptic curves.

The stack $\cX=[X/\ISL_3(\Z)]$ is the total space of the universal family of triptic curves over $\Tr$.
\end{proposition}

\begin{proof}
The map $[X/\Z^3]\to \C P^2-\R P^2$ is an $\SL_3(\Z)$-equivariant family of triptic curves.
Each fibre has a canonical orientation given by the generator of $H^3(\Z^3,\Z)=\Z$, and these orientations are compatible with the group action.
It follows that
\begin{equation}\label{universal curve}
[[X/\Z^3]/\SL_3(\Z)]=[X/\ISL_3(\Z)]\to [\C P^2-\R P^2/\SL_3(\Z)]
\end{equation}
is a family of oriented triptic curves.

To see that (\ref{universal curve}) is universal, consider an arbitrary family $\undertilde{\cE}\to S$ of oriented triptic curves.
We need to show that, locally, $\undertilde{\cE}$ is the pullback of $[X/\Z^3]$ via  some map $S\to \C P^2-\R P^2$,
and that this map is unique up to an element of $\SL_3(\Z)$.

Indeed, $\undertilde{\cE}$ can be written locally as $(S\times \C)/\Z^3$, where the action is given by $\mu(w,s)=(w-\iota_s(\mu),s)$ for some map
\begin{equation*}
\iota\mapsto\iota_s:S\to \hom(\Z^3,\C)=\C^3.
\end{equation*}
Modding out by rescaling we get a map to $\C P^2$, and
the condition on the rank of $\iota_s$ tells us that it doesn't hit $\R P^2$.
That map is then uniquely defined up to a change of basis of $\Z^3$.
\end{proof}

Before constructing the gamma gerbe on $\cX$, it is instructive to construct its restriction to the various triptic curves $\cE\subset \cX$.

\subsection{Gerbes on triptic curves}

Let $\cE$ be a triptic curve, let $\tilde \cE\simeq \C$ denote the universal cover of $\cE$, and let $\Lambda^\vee:=\tilde\cE\times_{\cE}\{0\}$.
Note that $\Lambda^\vee\simeq\Z^3$ comes with a canonical map $\iota$ to $\tilde \cE$ and that we have $[\tilde \cE/\iota(\Lambda^\vee)]=\cE$.
The action of $\Lambda^\vee$ on $\tilde\cE$ can also be identified with the action of $\pi_1(\cE)$ by deck transformations
and thus we get a canonical isomorphism $\Lambda^\vee=\pi_1(\cE)$.

Let $W:=\Lambda^\vee\otimes\R$ and define $K:=\ker(\iota_\R:W\to\tilde\cE)$.
Observe also that since $\mathrm{rk}(\iota)=2$, we always have $K\simeq\R$.
We shall also let $T:=W/\Lambda^\vee$ and $p:T\to\cE$ be the projection.
All this can then be assembled in the following short exact sequences (this being of course an abuse of language for what concerns the top and bottom rows):
\begin{equation}\label{lots of ses}
\begin{matrix}\xymatrix@R=.5cm{
&&\ar[d]0\\
&0\ar[r]&\Lambda^\vee\ar[r]^{\iota}\ar[d]&\tilde\cE\ar[r]&\cE\ar[r]&0\\
0\ar[r]&K\ar[r]&W\ar[r]^(.55){\iota_\R}\ar[d]&\tilde\cE\ar[r]&0\\
0\ar[r]&K\ar[r]&T\ar[rr]^p\ar[d]&&\cE\ar[r]&0\\
&&0
}\end{matrix}
\end{equation}

\begin{lemma}\label{lem:orient}
An orientation of $\cE$ in the sense of Definition \ref{def triptic} is equivalent to an orientation of $K$.
\end{lemma}

\begin{proof}
Since $\tilde\cE$ is contractible, the cohomology of $\cE$ can be identified with that of its fundamental group.
An orientation of $\cE$ is therefore the same thing as a choice of generator of $\Lambda^3(\pi_1\cE)\subset\Lambda^3(W)$,
which is the same thing as an orientation of $W$.

Since $\tilde\cE$ is a complex manifold, it comes with a canonical orientation.
We conclude by the middle row of (\ref{lots of ses})
that an orientation of $W$ is equivalent to an orientation of $K$.
\end{proof}

{}From now on, $\cE$ shall be an oriented triptic curve.
By virtue of Lemma \ref{lem:orient}, we shall then identify the ends of $K$ with $\{\infty,-\infty\}$.

Given an \'etale map $f:U\to \cE$, we let $T_U:=T\times_\cE U$, $Z_U:=\{0\}\times_\cE U$, and $p_U:T_U\to U$ be the projection.
\begin{equation}\label{two squares}
\begin{matrix}\xymatrix@R=.6cm@C=1cm{
\{0\}\hookright&T\ar[r]^p&\cE\\
Z_U\hookright\ar[u]&T_U\ar[r]^(.55){p_U}\ar[u]&U\ar[u]_f
}\end{matrix}
\end{equation}

\begin{lemma}
The space $T_U$ is a $K$-principal bundle over $U$ and $Z_U\to T_U$ is the inclusion of a discrete subset.
\end{lemma}

\begin{proof}
The first claim follows from the bottom row of (\ref{lots of ses}).
The second one holds because $\{0\}\hookrightarrow T$ is discrete and $T_U\to T$ is \'etale.
\end{proof}

We will say that a sequence of points in $Z_U$ {\em tends to $+\infty$} (respectively tends to $-\infty$) if its image in $U$ is relatively compact, and
if it tends to $+\infty$ (resp. $-\infty$) in the $K$-coordinate of $T_U$.
Let also ``$<$'' denote the corresponding total order on $K$.

\begin{definition}\label{def:pseudodivisor}
Let $f:U\to\cE$ be an \'etale map.
A {\em pseudodivisor} on $U$ is a function $D:Z_U\to \Z$ such that
for every sequence $y_n\in Z_U$ that tends to $+\infty$ we have $\lim D(y_n)=1$ and
for every sequence $z_n\in Z_U$ that tends to $-\infty$ we have $\lim D(z_n)=0$.
\end{definition}

We should note that pseudodivisors are not an empty notion.
Indeed, $K$-principal bundles over manifolds are always trivial.
So given an \'etale open $f:U\to \cX$, we may pick an isomorphism $\varphi:T_U\stackrel{\scriptscriptstyle\sim}{\to}U\times K$.
Pick $k_0\in K$ and let $\mathrm{pr}_2:U\times K\to K$ denote the projection.
The function $D:Z_U\to\Z$ given by
\begin{equation*}\label{D0}
D(y)=\begin{cases}1\quad\text{if } \mathrm{pr}_2(\varphi(y))\ge k_0\\
0\quad\text{if } \mathrm{pr}_2(\varphi(y))< k_0
\end{cases}
\end{equation*}
is then an example of an pseudodivisor on $U$.
Note also that pseudodivisors form a sheaf on $\cE$.

\begin{lemma}\label{pseudodivisors}
If $D_1$, $D_2$ are pseudodivisors, then $(p_U)_*(D_1-D_2)$ is a divisor on $U$, namely it takes finite values and has discrete support.
\end{lemma}

\begin{proof}
Let $U'\subset U$ be a subset with compact closure.
We need to show that the set of $y\in Z_{U'}$ such that $D_1(y)-D_2(y)\not=0$ is finite.

Since $T_U\to U$ is trivial, we may identify it with $U\times K$.
We shall thus write $y=(z,t)$, $z\in U$, $t\in K$, for a point $y\in Z_U$.
By definition of pseudodivisor, one may pick elements $a,b\in K$
such that $D_i(z,t)=0$ for all $(z,t)\in Z_U{'}$, $t<a$, and
such that $D_i(z,t)=1$ for all $(z,t)\in Z_{U'}$, $t>b$.
It follows that $D_1(z,t)-D_2(z,t)=0$ for all $(z,t)\in Z_{U'}$, $t\not\in [a,b]$.

Since $U'\times[a,b]$ is relatively compact in $T_U$ and since $Z_{U'}$ is discrete,
there is only a finite number of points $(z,t)\in Z_{U'}$ such that $t\in [a,b]$.
In particular, there is only a finite number of points such that $D_1(y)-D_2(y)\not=0$.
\end{proof}

Given two pseudodivisors, we shall henceforth abuse language and identify $D_1-D_2$ with its pushforward $(p_U)_*(D_1-D_2)$.

We are now in position to define the pseudodivisor gerbe on $\cE$.
It will assign to each \'etale map $f:U\to \cE$ a category
$\cG_U$, and to each pair of objects $\cL_1, \cL_2\in\cG_U$ a line bundle $\Homm(\cL_1,\cL_2)$ over $U$.

\begin{definition}\label{def:gerbe on fiber}
An object of $\cG_U$ is a pair $(L,D)$ consisting of a line bundle $L$ on $U$ and a pseudodivisor $D$ on $U$.
A morphism $(L_1,D_1) \to (L_2,D_2)$ is a section of the twisted line bundle $\big(\Homm(L_1,L_2)\big)(D_2-D_1)$,
where $\Homm(L_1,L_2)=L_1^*\otimes L_2$ denotes the usual internal hom of line bundles.
In other words
\begin{equation*}
\Homm_{\cG_U}\big((L_1,D_1),(L_2,D_2)\big):=\big(\Homm(L_1,L_2)\big)(D_2-D_1).
\end{equation*}

Note that this definition makes sense since by Lemma \ref{pseudodivisors}, $D_2-D_1$ is a divisor on $U$.
\end{definition}

We now globalize the above definition to the whole of $\cX$.

\subsection{The global construction}

Since all the constructions displayed in (\ref{lots of ses}) are canonical, they carry over to the universal triptic curve $\cX\to\Tr$.
The only difference is that now all the spaces will come with given maps to $\Tr$.
We shall call them $\underline{\Lambda}^\vee$, $\underline W$, $\underline T$, and $\underline K$ to indicate that fact.
We also have $\underline\cE=\cX=[X/\ISL_3(\Z)]$ and $\underline{\tilde\cE}=[X/\SL_3(\Z)]$.
For example, the last row of (\ref{lots of ses}) becomes the following short exact sequence of abelian group objects over $\Tr$:
\begin{equation*}
\xymatrix{
0\ar[r]&\underline{K}\ar[r]\ar[dr]&\underline{T}\ar[r]^p\ar[d]&\cX\ar[r]\ar[dl]&0\\&&\Tr
}
\end{equation*}
And since $\cX_3$ is a family of oriented triptic curves, we get by Lemma \ref{lem:orient} a canonical orientation of the line bundle $\underline{K}\to\Tr$.

Given an \'etale map $f:U\to \cX$, we now let $T_U:=\underline{T}\times_{\cX} U$, $Z_U:=\Tr\times_{\cX} U$, and let $p_U:T_U\to U$ denote the projection.
\begin{equation}\label{6 pull back}
\begin{matrix}\xymatrix@R=.6cm@C=1.4cm{
\Tr\hookright^{\substack{\text{zero}\\ \text{section}}}&\underline{T}\ar[r]^p&\cX\\
Z_U\hookright\ar[u]&T_U\ar[r]^(.55){p_U}\ar[u]&U\ar[u]_f
}\end{matrix}
\end{equation}
Once again, $T_U$ is a $\underline{K}$-principal bundle over $U$.
The only difference with (\ref{two squares}) is that $Z_U\hookrightarrow T_U$ is not discrete any more.
Instead, it is a union of submanifolds of (real) codimension 3.

The analogues of Definition \ref{def:pseudodivisor} and Lemma \ref{pseudodivisors} are now straightforward.

\begin{definition}\label{def:pseudodivisor global}
Let $f:U\to\cX$ be an \'etale map.
A pseudodivisor on $U$ is a continuous function $D:Z_U\to \Z$ such that
for every point $\{q\}\to\Tr$ with corresponding fibre $\cE=\{q\}\times_\Tr \cX$, the restriction of $D$ to
$\{q\}\times_\Tr Z_U$ is a pseudodivisor in the sense of Definition \ref{def:pseudodivisor}.
\end{definition}

\begin{lemma}
If $D_1$, $D_2$ are pseudodivisors, then $(p_U)_*(D_1-D_2)$ is a divisor on $U$.
\end{lemma}

\begin{proof}
The zero section $\Tr\to\cX$ is locally like the inclusion of a complex hypersurface.
Therefore the same holds for $p_U:Z_U\to U$.
So we just need to show that $(p_U)_*(D_1-D_2)$ is supported on locally finitely many pieces.

Let $U'\subset U$ be a subset with compact closure.
By the same argument as in Lemma \ref{pseudodivisors}, the support of $(D_1-D_2)|_{Z_{U'}}$ is
relatively compact in $T_U$.
The pushforward $(p_U)_*(D_1-D_2)$ is therefore a divisor.
\end{proof}

The pseudodivisor gerbe $\cG$ on $\cX$ (soon to be identified with the gamma gerbe) can now be defined in a way much similar to Definition \ref{def:gerbe on fiber}.
Once again, we identify $D_1-D_2$ with $(p_U)_*(D_1-D_2)$.

\begin{definition}\label{def: global gerbe}
To an \'etale map $U\to\cX$, the pseudodivisor gerbe $\cG$ assigns the category $\cG_U$ given as follows:
An object of $\cG_U$ is a pair $(L,D)$ consisting of a line bundle $L$ over $U$ and a pseudodivisor $D$ on $U$.
The morphisms are given by
\begin{equation*}
\Homm_{\cG_U}\big((L_1,D_1),(L_2,D_2)\big):=\big(\Homm(L_1,L_2)\big)(D_2-D_1).
\end{equation*}
\end{definition}

\subsection{Identifying the two gerbes}

Let $\cG^\Gamma$, $\cG^{div}$ denote the gamma gerbe and the pseudodivisor gerbe respectively, as given in Definitions \ref{def: gamma gerbe} and \ref{def: global gerbe}.
We shall show that $\cG^\Gamma\simeq\cG^{div}$ by producing a nowhere vanishing
holomorphic section of $\cG^\Gamma\otimes(\cG^{div})^*$.

To describe $(\cG^{div})^*$, we introduce the notion of a {\em dual pseudodivisor}.
It is obtained replacing the condition
``for every sequence $y_n\in Z_U$ that tends to $+\infty$ we have $\lim D(y_n)=1$'' by
``$\ldots=-1$''.
The dual gerbe $(\cG^{div})^*$ is then given by taking Definition \ref{def: global gerbe} and replacing pseudodivisors by dual pseudodivisors.

The holomorphic section of $\cG^\Gamma\otimes(\cG^{div})^*$
will be defined with respect to the same invariant covering $\{V_a\}$ that was used in section \ref{sec:the gamma gerbe}.
It consists of the following data:
\begin{itemize}
\item For each $a\in\Lambda_{\mathit{prim}}$, a pair $(L_a,D_a)$, where $L_a$ is a line bundle over $V_a$, and $D_a:Z_{V_a}\to \Z$ is a dual pseudodivisor.
\item For each $a,b\in\Lambda_{\mathit{prim}}$, a nowhere vanishing holomorphic section $s_{a,b}$ of the line bundle $\big(\Homm(L_b,L_a)\big)(D_a-D_b)$.
\item For each $a\in \Lambda_{\mathit{prim}}$, $g\in G$ a non-vanishing holomorphic section $s_{a;g}$ of the line bundle $\big(\Homm(L_a,g_*\,L_{g^{-1}a})\big)(g_*D_{g^{-1}a}-D_a)$.
\end{itemize}
satisfying the relations (we set $\phi_{a,b;g}(y)=\phi_{a,b}(g;y)$ etc.)
\begin{gather}
\phi_{a,b,c}\,s_{a,c}=s_{a,b}\circ s_{b,c},\nonumber\\
\phi_{a,b;g}\,g_*(s_{g^{-1}a,g^{-1}b})\circ s_{b;g}= s_{a;g}\circ s_{a,b},\label{3 cocycle formulas}\\
\phi_{a;g,h}\,s_{a;gh}=s_{a;g}\circ g_*(s_{g^{-1}a;h}).\nonumber
\end{gather}
The above relations are probably best understood by saying that the following diagrams commute
up to multiplication by the functions drawn in their middle.
\begin{gather*}
\qquad\begin{matrix}\xymatrix@C=.6cm{
&L_b\ar[dl]_{s_{a,b}}\ar@{}[d]|(.6){\phi_{a,b,c}}&\\
L_a&&L_c\ar[ll]^{s_{a,c}}\ar[ul]_{s_{b,c}}
}\end{matrix}
\qquad\quad
\begin{matrix}\xymatrix{
g_*L_{g^{-1}a}&L_a\ar[l]_(.4){s_{a;g}}\ar@{}[dl]|(.4){\phi_{a,b;g}}\\
g_*L_{g^{-1}b}\ar[u]|{g_*(s_{g^{-1}a,g^{-1}b})\hspace{1cm}}&L_b\ar[l]^(.4){s_{b;g}}\ar[u]_{s_{a,b}}
}\end{matrix}
\\
\begin{matrix}\xymatrix{
&g_*L_{g^{-1}a}\ar[dl]_{g_*(s_{g^{-1}a;h})}\ar@{}[d]|(.6){\phi_{a;g,h}\hspace{4mm}}&\\
(gh)_*L_{(gh)^{-1}a}&&L_a\ar[ll]^(.4){s_{a;gh}}\ar[ul]_{s_{a;g}}
}\end{matrix}
\end{gather*}

We now proceed to construct these pieces of data.
First of all, we let $L_a:=V_a\times \C$ be the trivial line bundle.
To write down $D_a$, we need some more explicit descriptions of $T_a:=T_{V_a}$, $Z_a:=Z_{V_a}$, and $p:= p_{V_a}$.

The space $\C P^2-\R P^2$ is the moduli space
of triptic curves $\cE$ equipped with a trivialisation of their fundamental groups $\Lambda^\vee=\pi_1(\cE)$.
So when restricted to $\C P^2-\R P^2$, the bundles $\underline{\Lambda}^\vee$, $\underline{W}$, $\underline{T}$ defined in (\ref{lots of ses}) have canonical trivialisations
\begin{equation*}
\begin{split}
\underline{\Lambda}^\vee\times_{\Tr}(\C P^2-\R P^2)&=\Z^3\times(\C P^2-\R P^2),\\
\underline{W}\times_{\Tr}(\C P^2-\R P^2)&=\R^3\times(\C P^2-\R P^2),\\
\underline{T}\times_{\Tr}(\C P^2-\R P^2)&=(\R^3/\Z^3)\times(\C P^2-\R P^2).\\
\end{split}
\end{equation*}
Recall that $V_a=(\C\times U_a^+)/\C^\times$ is an open subset of $X_3$.
We can then compute $T_a$ and $Z_a$ by means of the following pullback squares
\begin{equation*}
\xymatrix@C=.7cm{
\Tr\ar[r]^{\text{zero section}}&\underline{T}\ar[r]&\cX_3\ar[r]&\scriptstyle\Tr\\
\scriptstyle(\C P^2-\R P^2)\ar[r]\ar[u]&\scriptstyle(\R^3/\Z^3)\times(\C P^2-\R P^2)\ar[r]\ar[u]&\scriptstyle X_3/\Z^3\ar[u]\ar[r]&\scriptstyle(\C P^2-\R P^2)\ar[u]\\
\scriptstyle\Z^3\times(\C P^2-\R P^2)\ar[r]\ar[u]&\scriptstyle\R^3\times(\C P^2-\R P^2)\ar[r]\ar[u]&\scriptstyle X_3\ar[u]\\
Z_a=\Z^3\times(U_a^+/\C^\times)\ar[r]\ar[u]&T_a=\R^3\times(U_a^+/\C^\times)\ar[r]^(.75){p}\ar[u]&V_a\ar[u]
}
\end{equation*}
The projection $p$ can then be written explicitly as $p:(\mu,[x])\mapsto [(\mu(x),x)]$.
Given the above description of $Z_a$, we can now define $D_a$:
\[
D_a:Z_a\to\Z:=\begin{cases}
-1&\text{on }\Lambda_+^\vee(a)\times(U_a^+/\C^\times)\\
0&\text{on }\Lambda_-^\vee(a)\times(U_a^+/\C^\times)
\end{cases}
\]
The definition of $U_a^+$ is made up so that all sequences $y_n\in \Lambda^\vee$ tending to $+\infty$ eventually land in $\Lambda_+^\vee(a)$, and that all sequences tending to $-\infty$ eventually land in $\Lambda_-^\vee(a)$.
It follows that $D_a$ is indeed a dual pseudodivisor.

Since the $L_a$ are trivial, $s_{a,b}$ and $s_{a;g}$ are just functions.
We let $s_{a,b}:=\Gamma_{a,b}$ and $s_{a;g}:=\Delta_{a}(g;\cdot)$.
The equations (\ref{3 cocycle formulas}) are then identical to those in Theorem \ref{t-Quiverful}.

So we just need to show that when viewed as section of $\cO(D_a-D_b)$ and $\cO(g_*D_{g^{-1}a}-D_a)$ respectively,
$s_{a,b}$ and $s_{a;g}$ are non-vanishing and holomorphic.
In other words, we must show that the divisors of $s_{a,b}$ and $s_{a;g}$ are
given by $-p_*(D_a-D_b)$ and $-p_*(g_*D_{g^{-1}a}-D_a)$.
For this we rewrite (\ref{e-imperius}) formally as
\begin{equation}\label{it's just formal}
\Gamma_{a,b}(w,x)=\frac
{\prod_{\delta\in \Lambda_+^\vee(a)/\mathbb Z\gamma}(1-e^{-2\pi i(\delta(x)-w)/\gamma(x)})}
{\prod_{\delta\in \Lambda_+^\vee(b)/\mathbb Z\gamma}(1-e^{2\pi i(\delta(x)-w)/\gamma(x)})}\,.
\end{equation}
Neither the numerator nor the denominator converge, but (\ref{it's just formal}) is good enough to compute the zeroes and poles of $\Gamma_{a,b}$.
The numerator has zeroes at $p\big(\Lambda_+^\vee(a)\times(U_a^+/\C^\times)\big)$
and the denominator zeroes at $p\big(\Lambda_+^\vee(b)\times(U_a^+/\C^\times)\big)$.
It follows that the divisor of $s_{a,b}=\Gamma_{a,b}$ is $-p_*(D_a-D_b)$, which is what we wanted to show.
Similarly, we can rewrite (\ref{e-DrStanhope}) formally as
\begin{equation}\label{again it's just formal}
\Delta_{a}(g;w,x)=\frac
{\,\prod_{j=-\infty}^{\mu(a)-1}\,\theta_0\textstyle\left( \frac{w+j\alpha_1(x)}{\alpha_3(x)},\frac{\alpha_2(x)}{\alpha_3(x)}\right)\,}
{\prod_{j=-\infty}^{-1}\,\theta_0\textstyle\left( \frac{w+j\alpha_1(x)}{\alpha_3(x)},\frac{\alpha_2(x)}{\alpha_3(x)}\right)}
\end{equation}
where $\mu$ is the $\Z^3$ component of $g\in \ISL_3(\Z)$.
The numerator has zeroes at
$p\big((\mu+\Lambda_+^\vee(a))\times(U_a^+/\C^\times)\big)$
and the denominator has zeroes at $p\big(\Lambda_+^\vee(a)\times(U_a^+/\C^\times)\big)$.
It follows that the divisor of $s_{a;g}=\Delta_{a}(g;\cdot)$ is $-p_*(\mu_*D_a-D_a)=-p_*(g_*D_{g^{-1}a}-D_a)$, as desired.
Thus we have proven:
\begin{thm}
The  two gerbes over $\cX$ described in Definitions \ref{def: gamma gerbe} and \ref{def: global gerbe} are isomorphic.
\end{thm}

\section{The Brauer group of triptic curves}\label{s-Brauer}

In this section, we shall compute the Brauer group $Br(\cE):=H^2(\cE,\cO^\times)$ of a generic triptic curve $\cE$.
We shall then identify the class of the gamma gerbe inside this Brauer group.

\subsection{The cohomology of $\cE$ with coefficients in $\Oanx$}

Let $x_1,x_2,x_3$ be three complex numbers, and
let $\iota:\Z^3\to \C$ denote the linear map sending $e_i$ to  $x_i$.
We assume for the moment that the $x_i$ are $\Q$-linearly independent and that they span $\C$ as an $\R$-vector space.
So in particular $\iota$ is injective.
Let $\cE$ denote the triptic curve\footnote{Note that since $\iota$ is injective, the stack $[\C/\iota(\Z^3)]$ is actually a sheaf.} $[\C/\iota(\Z^3)]$.

Let $\Oanx=(\Oanx,\cdot)$ the sheaf of invertible analytic functions.
The goal of this section is to compute the cohomology of $\cE$ with coefficients in $\Oanx$,

In order to avoid redundant arguments, we put ourselves in a slightly more general situation.
So let $r\ge 2$ be an integer, let $x_j$, $j=1\ldots r$ be $\Q$-linearly independent complex numbers,
let $\iota$ be the map $\Z^r\to\C:e_j\mapsto \tau_j$, and let $\cE=[\C/\iota(\Z^r)]$.
We also assume that $x_1$ and $x_2$ span a lattice in $\C$, and call $E$ the elliptic cure $\C/\Z\{x_1,x_2\}$.

Let $\Oan=(\Oan,+)$ denote the sheaf of analytic functions, and
consider the exponential sequence
\begin{equation}\label{expSES}
0\longrightarrow\Z\longrightarrow\Oan\longrightarrow\Oanx\longrightarrow 0.
\end{equation}
We will compute $H^*(\cE,\Oanx)$ by first calculating $H^*(\cE,\Z)$ and $H^*(\cE,\Oan)$, and then using
the long exact sequence in cohomology.

We first compute $H^*(\cE,\Z)$.
Since $\Z$ is a constant sheaf, that's just ordinary integral cohomology applied to the \'etale homotopy type of $\cE$ \cite[Corollary 9.3]{AM69}.
Let $N$ denote the \v Cech nerve associated to the cover $\C\to \cE$.
It is the simplicial manifold given by all the
$n$-tuple fibered products $\C\times_\cE\ldots\times_\cE\C=(\Z^r)^{n-1}\times\C$.
\begin{equation}\label{it's a good cover}
N=\Big(\cdots\rrrrarrow(\Z^r)^2\times\C\rrrarrow\Z^r\times\C\rrarrow\C\Big)
\end{equation}
Since each degree of (\ref{it's a good cover}) is a disjoint union of contractible components, the cover
$\C\to \cX$ is a good cover.
We can therefore compute the \'etale homotopy type of $\cE$ by applying $\pi_0$ degreewise to $N$.
The resulting simplicial set is then the standard simplicial model for $B\Z^r$.
We deduce that
\[
H^*(\cE,\Z)=H^*(B\Z^r,\Z)=\Lambda^{\!*}(\Z^r).
\]
Morally, this argument can be summarized as follows.
The sheaf $\cE$ is a free quotient of a contractible space by the group $\Z^r$, it is therefore a model for $B\Z^r$.

Let $\Opol$ denote the locally constant sheaf of polynomial functions on $\C$, and let the same symbol also denote its pullback to $\cE$.
We will use $\Opol$ as a first approximation to $\Oan$.
Let $\Ole{n}$ denote the subsheaf of $\Opol$ of functions of degree at most $n$.
Clearly, $\Opol$ is the colimit of the sheaves $\Ole{n}$.

\begin{lemma}\label{compact}
The group $H^i(\cE,\Opol)$ is the colimit over $n$ of $H^i(\cE,\Ole{n})$.
\end{lemma}

\begin{proof}
The sheaves $\Ole{n}$ and $\Opol$ are locally constant and their stalks are given respectively by $\C\{1,\ldots,z^n\}$ and $\C[z]$.
It follows that
\[
H^i(\cE,\Opol)=H^i(B\Z^r,\C[z]),\qquad H^i(\cE,\Ole{n})=H^i(B\Z^r,\C\{1,\ldots,z^n\}),
\]
where the right hand side denotes cohomology with local coefficients.
The result follows since $B\Z^r$ admits a model with finitely many cells.
\end{proof}

\begin{remark}
If $\{\cF_n\}$ is an arbitrary directed system of sheaves on $\cE$, then it is also true that $H^i(\cE,\colim\cF_n)=\colim H^i(\cE,\cF_n)$.
This holds because $\cE$ can be resolved by a simplicial manifold which is compact in each degree.
\end{remark}

We now compute $H^*(\cE,\Opol)$ using the spectral sequence associated to the filtration by the $\Ole{n}$.
This spectral sequence converges by virtue of Lemma \ref{compact}.
The associated graded sheaves $\Ole{n}/\Ole{n-1}$ being isomorphic to the constant sheaf $\C$, we have
\begin{equation*}
H^*(\cE,\Ole{n}/\Ole{n-1})=H^*(\cE,\C)=H^*(B\Z^r,\C)=\Lambda^{\!*}(\C^r).
\end{equation*}
So our spectral sequence looks like
\begin{equation}\label{statement of SS}
E^{i,n}_1=\Lambda^{\!i}(\C^r)\Rightarrow H^i(\cE,\Opol),
\end{equation}
with differentials $d_k$ of bidegree $(1,-k)$.

Let $\alpha_j$, $j=1\ldots r$ denote the generators of $E_1^{1,0}=\C^r$.
We want to compute the image under $d_1$ of the standard generator $z\in E_1^{0,1}=\C$.
For this, we need to find an $(\Ole{1})$-valued \v Cech 0-cochain $c$ lifting the $(\Ole{1}/\Ole{0})$-valued 0-cocycle $z$, and take its coboundary.
We define $c$ with respect to the cover $\C\twoheadrightarrow \C/\iota(\Z^r)=\cE$.
It is then the usual coordinate $z$ on $\C$.
Its coboundary $\delta c$ is then a function on $\C\times \Z^r$ which can easily be computed to be $\sum a_i x_i$ on the component
$\C\times (a_1,\ldots,a_r)$.
This $(\Ole{0})$-valued 1-cocycle represents $\sum x_i\alpha_i$.
We have therefore computed
\begin{equation}\label{d1computation}
d_1(z)=\sum x_i\alpha_i.
\end{equation}
All the other $d_1$ differentials follow from (\ref{d1computation}) by multiplicativity.

More precisely, letting $  x=(x_1,\ldots,x_r)\in\C^r$, we can identify
$d_1:E_1^{i,n}\to E_1^{i+1,n-1}$ with $-\wedge n\!  x:\Lambda^{\!i}(\C^r)\to\Lambda^{\!i+1}(\C^r)$.
The sequence
\[
\Lambda^{\!i-1}(\C^r)
\stackrel{\wedge\! \textstyle  x}{\longrightarrow}
\Lambda^{\!i}(\C^r)
\stackrel{\wedge\! \textstyle  x}{\longrightarrow}
\Lambda^{\!i+1}(\C^r)
\]
being exact, we deduce that $E_2^{i,n}$ is zero for $n\ge 0$ and
$\Lambda^{\!i}(\C^r)\big/(\Lambda^{\!i-1}(\C^r)\!\wedge\!   x)$ for $n=0$.
The latter can also be identified with $\Lambda^{\!i}(\C^r/  x\C)$.
We have then proven
\begin{proposition}\label{Opol computation}
$
H^*(\cE,\Opol)=\Lambda^{\! *}(\C^r/  x\C).
$
\hfill$\square$
\end{proposition}

For $r=3$, this is what the spectral sequence (\ref{statement of SS}) looks like in coordinates:
\begin{equation}\label{OlenSS}
  \put(6,15){\vdots}
  \put(79,15){\vdots}
  \put(157,15){\vdots}
  \put(233,15){\vdots}
\xymatrix@C2cm{
\ar[dr]|{2\cdot\!\Big(
  \put(-1,8){$\scriptstyle x_1$}
  \put(-1,2){$\scriptstyle x_2$}
  \put(-1,-4){$\scriptstyle x_3$} \hspace{.25cm}
\Big)}
\C&
\ar[dr]|{2\cdot\!\Big(
  \put(-1,8){-$\scriptstyle x_2$}
  \put(-1,2){-$\scriptstyle x_3$}
  \put(3,-4){$\scriptscriptstyle 0$}
  \put(13,8){$\scriptstyle x_1$}
  \put(14,2){$\scriptscriptstyle 0$}
  \put(10,-4){-$\scriptstyle x_3$}
  \put(24,8){$\scriptscriptstyle 0$}
  \put(23,2){$\scriptstyle x_1$}
  \put(23,-4){$\scriptstyle x_2$} \hspace{1.05cm}
\Big)_{\phantom *}}
\C^3
&
\ar[dr]|{
2\cdot({
\scriptstyle x_3\,
    \text{-} x_2\:
            x_1})_{\phantom *}}
\C^3
&\C\\
\ar[dr]|{
\Big(
  \put(-1,8){$\scriptstyle x_1$}
  \put(-1,2){$\scriptstyle x_2$}
  \put(-1,-4){$\scriptstyle x_3$} \hspace{.25cm}
\Big)}
\C
&
\ar[dr]|{\Big(
  \put(-1,8){-$\scriptstyle x_2$}
  \put(-1,2){-$\scriptstyle x_3$}
  \put(3,-4){$\scriptscriptstyle 0$}
  \put(13,8){$\scriptstyle x_1$}
  \put(14,2){$\scriptscriptstyle 0$}
  \put(10,-4){-$\scriptstyle x_3$}
  \put(24,8){$\scriptscriptstyle 0$}
  \put(23,2){$\scriptstyle x_1$}
  \put(23,-4){$\scriptstyle x_2$} \hspace{1.05cm}
\Big)}
\C^3
&
\ar[dr]|{
({\scriptstyle x_3\,
      \text{-} x_2\:
               x_1})_{\phantom *}}
\C^3&\C\\
\C&\C^3&\C^3&\C
}
\end{equation}

For $r=2$, $\cX=E$ is an elliptic curve.
We have then computed that $H^0(E,\Opol)=H^1(E,\Opol)=\C$, and that the other cohomology groups vanish.
By Riemann-Roch and Serre's GAGA principle, we also have $H^0(E,\Oan)=H^1(E,\Oan)=\C$, and $H^i(E,\Oan)=0$ for $i\ge 2$.
This shows that $H^i(E,\Opol)$ and $H^i(E,\Oan)$ are abstractly isomorphic.

\begin{lemma}\label{Hpq}
The natural map $H^i(E,\Opol)\to H^i(E,\Oan)$ is an isomorphism.
\end{lemma}

\begin{proof}
Clearly, the only nontrivial case is $i=1$.
Since $H^1(E,\C)\to H^1(E,\Opol)$ is surjective and $\dim H^1(E,\Opol)=\dim H^1(E,\Oan)$,
it is enough to show that $H^1(E,\C)\to H^1(E,\Oan)$ is surjective.
Consider the short exact sequence
\[
0\longrightarrow\C\longrightarrow \Oan\stackrel{d}{\longrightarrow}\Omega_{\!an}^1\longrightarrow 0
\]
and the corresponding long exact sequence in cohomology
\[
H^1(E,\C)\to H^1(E,\Oan)\to H^1(E,\Omega_{\!an}^1)\to H^2(E,\C)\to H^2(E,\Oan)
\]
The map $H^1(E,\Omega_{\!an}^1)\to H^2(E,\C)$ is surjective because $H^2(E,\Oan)=0$.
But it is also the inclusion of the middle summand in the Hodge decomposition
$H^2=H^{2,0}\oplus H^{1,1}\oplus H^{0,2}$.
So it is an isomorphism.
It follows that $H^1(E,\C)\to H^1(E,\Oan)$ is surjective.
\end{proof}

As a corollary, we have the following result.

\begin{lemma}
The cohomology groups $H^i(\cE,\Oan/\Opol)$ are zero for all $i\ge 0$ and all $r\ge 2$.
\end{lemma}

\begin{proof}
For $r=2$, this follows from Lemma \ref{Hpq} and from the long exact sequence in cohomology.
For $r>2$, let us identify the cohomology of $\cE$ with the $(\Z^{r-2})$-equivariant
cohomology of $E$.
The latter can be computed as the cohomology of the simplicial manifold
\begin{equation}\label{actionnerve}
N=\left(\cdots\rrrrarrow\big(\Z^{r-2}\big)^2\times E\rrrarrow\Z^{r-2}\times E\rrarrow E\right).
\end{equation}
We then have a spectral sequence \cite[Prop. 3.2.]{Fri82}
\[
E_1^{i,j}=H^i(N_j,\Oan/\Opol)\Rightarrow H^{i+j}(N,\Oan/\Opol)=H^{i+j}(\cX,\Oan/\Opol)
\]
Since each $N_j$ is a disjoint union of $E$'s, we know from the $r=2$ case that $H^i(N_j,\Oan/\Opol)=0$.
The spectral sequence vanishes identically thus proving the result.
\end{proof}

\begin{corollary}
The natural map  $H^i(\cE,\Opol)\to H^i(\cE,\Oan)$ is an isomorphism.
\end{corollary}

We can now state the main theorem of this section.

\begin{thm}\label{thm:HOan}
The cohomology of $\cE$ with coefficients in $\Oanx$ reads
\begin{align}
H^i(\cE,\Oanx)\hspace{.43cm}&=\Lambda^{\!i}(\C^r/  x\C)\big/\Lambda^{\!i}(\Z^r)&
                                                                &\text{for}\quad 0\le i\le r-2\label{i le r-2}\\
H^{r-1}(\cE,\Oanx)&=\cE\times\Z\label{i = r-1}\\
H^i(\cE,\Oanx)\hspace{.43cm}&=0&                                &\text{for}\quad i\ge r,\label{i ge r}
\end{align}
where $x$ denotes the vector $(x_1,\ldots,x_r) \in \C^r$.

The map $\cE\to\cE:w\mapsto w+z$ induces the identity on $H^i(\cE,\Oanx)$ for $i\le r-2$, and induces
$(w,n)\mapsto(w+nz,n)$ on $H^{r-1}(\cE,\Oanx)=\cE\times\Z$.
\end{thm}

\begin{proof}
The cohomology $H^*(\cE,\Z)=\Lambda^{\! *}(\Z^r)$ lands in the first row of the spectral sequence (\ref{OlenSS}) and hits the standard basis of
$\Lambda^{\! *}(\C^r)$.
The map $H^i(\cE,\Z)\to H^i(\cE,\Oan)$ is therefore the composite
\begin{equation}\label{it's injective}
H^i(\cE,\Z)=\Lambda^{\! i}(\Z^r)
\longrightarrow    \Lambda^{\! i}(\C^r)
\longrightarrow
\Lambda^{\!i}(\C^r/  x\C)=H^i(\cE,\Oan).
\end{equation}
By Lemma \ref{Qbasis}, (\ref{it's injective}) is injective for $i\le r-1$.
The long exact sequence corresponding to (\ref{expSES}) therefore splits
into short exact sequences
\begin{equation}\label{i le r-2 SES}
0\longrightarrow H^i(\cE,\Z) \longrightarrow H^i(\cE,\Oan) \longrightarrow H^i(\cE,\Oanx) \longrightarrow 0
\end{equation}
for $i\le r-2$, and a four term exact sequence
\begin{equation}\label{i = r-1 4TES}
0\rightarrow H^{r-1}(\cE,\Z)\rightarrow H^{r-1}(\cE,\Oan)\rightarrow H^{r-1}(\cE,\Oanx)\rightarrow H^r(\cE,\Z)\rightarrow 0.
\end{equation}
Equation (\ref{i le r-2}) follows from (\ref{i le r-2 SES}) and our computation of $H^i(\cE,\Z)$ and $H^i(\cE,\Oan)$.

To get (\ref{i = r-1}), we note that $H^r(\cE,\Z)=\Z$.
The rightmost arrow of (\ref{i = r-1 4TES}) therefore splits, and we get
\[
H^{r-1}(\cE,\Oanx)=\Big(\Lambda^{\!r-1}(\C^r/  x\C)\big/\Lambda^{\!r-1}(\Z^r)\Big)\times \Z.
\]
To identify $\Lambda^{\!r-1}(\C^r/  x\C)\big/\Lambda^{\!r-1}(\Z^r)$ with $\cE$,
consider the isomorphism
\begin{equation}\label{both lines}
\begin{split}
\Lambda^{\!r-1}(\C^r/  x\C)\hspace{2cm}&\to\,\,\C\\
\overline e_j:=(-1)^j e_1\wedge\ldots\wedge e_{j-1}\wedge e_{j+1}\wedge\ldots \wedge e_r&\mapsto\,\, x_j
\end{split}
\end{equation}
It is an isomorphism because it is not zero and $\Lambda^{\!r-1}(\C^r/  x\C)$
is one dimensional.
It is well defined because the relation $x_k\overline e_j=x_j\overline e_k$ holds in
$\Lambda^{\!r-1}(\C^r/  x\C)$:
\[
x_k\overline e_j-x_j\overline e_k=
  x\wedge
\big(e_1\wedge\ldots\wedge e_{j-1}\wedge e_{j+1}\wedge\ldots\wedge e_{k-1}\wedge e_{k+1}\wedge\ldots \wedge e_r\big)=0.
\]
We now observe that
\[
\Lambda^{\!r-1}(\C^r/  x\C)\big/\Lambda^{\!r-1}(\Z^r)=
\Lambda^{\!r-1}(\C^r/  x\C)\big/\Z\{\overline e_j\}=
\C/\Z\{s_j\}=\cE,
\]
thus proving (\ref{i = r-1}).
At last, (\ref{i ge r}) holds because $H^i(\cE,\Oan)=H^{i+1}(\cE,\Z)=0$ for $i\ge r$.

We now examine the map $z_*:H^*(\cE,\Oanx)\to H^*(\cE,\Oanx)$ induced by $w\mapsto w+z$.
The latter being homotopic to the identity, it induces the identity on $H^*(\cE,\C)$.
As a consequence, $z_*$ is the identity on the spectral sequence (\ref{statement of SS}), and thus also on the short exact sequence (\ref{i le r-2 SES}).
This shows that $z_*$ is the identity on $H^{\le r-2}(\cE,\Oanx)$.

For similar reasons, the induced map on $H^{r-1}(\cE,\Oanx)=\cE\times\Z$ is the identity on the subgroup $\cE$ and on the quotient $\Z$.
To further determine $z_*$, we consider the quotient $\cE':=[\cE/z\Z]$.
Its cohomology is related to that of $\cE$ by means of the following long exact sequence:
\[
\ldots\longrightarrow H^i(\cE,\Oanx)\stackrel{1-z_*}{\longrightarrow}H^i(\cE,\Oanx)\longrightarrow H^{i+1}(\cE',\Oanx)\longrightarrow\ldots
\]
Letting $i=r-1$, we get
\[
\ldots\longrightarrow \cE\times\Z\stackrel{1-z_*}{\longrightarrow}\cE\times\Z\longrightarrow \cE'\times\Z\longrightarrow 0.
\]
It follows that (up to a choice of sign), the map $1-z_*$ is given by $(w,n)\mapsto (nz,0)$.
\end{proof}

As a special case of Theorem \ref{thm:HOan}, we have:

\begin{corollary}\label{thm:HOan r=3}
The cohomology of a generic triptic curve $\cE=\C/\iota(\Z^3)$ is given by
\begin{alignat*}{3}
H^0(\cE,\Oanx)&=\C^\times&
H^1(\cE,\Oanx)&=\big(\C^\times\big)^3\big/(e^{\lambda x_1}\!,e^{\lambda x_2}\!,e^{\lambda x_3}),\quad \lambda\in\C\\
H^2(\cE,\Oanx)&=\cE\times\Z&\qquad
H^i(\cE,\Oanx)&=0\quad\text{for}\quad i\ge 3.
\end{alignat*}
\end{corollary}

\begin{lemma}\label{Qbasis}
The map $\Lambda^{\!i}(\Z^r)\to\Lambda^{\!i}(\C^r/  x\C)$ is injective for all $i\le r-1$.
\end{lemma}

\begin{proof}
Let $b_1:=1$, $b_j:=- x_j/ x_1$, $j=2\ldots r$. We then have
\[
e_1=\sum_{j=2}^{r}b_je_j
\]
in $\C^r/ \tau\C$.
The $x_j$ were assumed $\Q$-linearly independent, so we may complete the set of $b$'s to a $\Q$-basis $\{b_\alpha\}$ of $\C$.
Since $\{e_j\}_{2\le j\le r}$ is a $\C$-basis of $\C^r/  x\C$,
the set $\{e_J\}$, $J\subset\{2,\ldots,r\}$, $|J|=i$ is a $\C$-basis of $\Lambda^{\!i}(\C^r/  x\C)$.
It follows that
\begin{equation}\label{horrible basis}
\big\{b_\alpha e_J\big\}_{J\subset\{2,\ldots,r\},\, |J|=i}
\end{equation}
is a $\Q$-basis of $\Lambda^{\!i}(\C^r/  x\C)$.
We now expand the $\Z$-basis $\{e_I\}$ of $\Lambda^{\!i}(\Z^r)$ in terms of (\ref{horrible basis}):
\begin{alignat*}{3}
e_I&=b_1e_I& &\text{if}\quad 1\not\in I\\
e_I&=e_1\wedge e_J=\sum_{j=2}^{r}b_je_j\wedge e_J=\sum_{\substack{j=2,\\j\not\in J}}^{r}b_j\,e_{J\cup\{j\}},& \qquad &\text{if}\quad I=\{1\}\cup J.
\end{alignat*}
Each element of (\ref{horrible basis}) appears in the expansion of at most one $e_I$.
Unless one of the $e_I$ maps to zero (which happens only if $i=r$),
this means that they are $\Q$-linearly independent in $\Lambda^{\!i}(\C^r/  x\C)$.
The map $\Lambda^{\!i}(\Z^r)\to\Lambda^{\!i}(\C^r/  x\C)$ is therefore injective.
\end{proof}

\subsection{The Dixmier--Douady class}

We now proceed to identify the class of the gamma gerbe inside the Brauer group $Br(\cE)=H^2(\cE,\cO^\times)=\cE\times\Z$.
We first determine its Dixmier--Douady class.

\begin{proposition}\label{t-hofame}
Let $\cG$ be the restriction of the gamma gerbe to a generic triptic curve $\cE$.
Then its Dixmier--Douady class $$c(\cG)\in H^3(\cE,\Z)=\Z$$ is a generator.
\end{proposition}

\begin{proof}
Pick an open $V_a$ covering $\cE$.
By Theorem \ref{t-main}, a cocycle representing $c(\mathcal G)$
is given by
\[
\psi=\frac1{2\pi i}\delta\log\phi_a\in C^3(\mathbb Z^3,\check C^0(V_a,\mathbb Z)).
\]
This can be computed from
Proposition \ref{p-EleanorBold}: in the
notation used there, let
$\lambda=\sum_i\ell_i\alpha_i$,
$\mu=\sum_i m_i\alpha_i$,
$\nu=\sum_i n_i\alpha_i$. Then
$\psi$ is the cocycle
\begin{eqnarray*}
\psi(\lambda,\mu,\nu)&=&
\delta P_a(\lambda,\mu,\nu)
\\
&=&
-P_a(\mu,\nu;w+\lambda(x),x)+
P_a(\lambda+\mu,\nu;w,x)\\
&&-
P_a(\lambda,\mu+\nu;w,x)+P_a(\lambda,\mu;w,x)
\\
&=&-m_1n_2\ell_3,
\end{eqnarray*}
representing a generator
of $H^3(\mathbb Z^3,\mathbb Z)=\mathbb Z$.
\end{proof}

It is now tempting to try to determine the $\cE$ coordinate of $c(\cG)$ in $Br(\cE)=\cE\times\Z$.
However, this question doesn't make any sense.
Indeed, we don't have yet a canonical isomorphism between $Br(\cE)$ and $\cX\times\Z$, we just have the short exact sequence
\begin{equation}\label{eq: splittable}
0\longrightarrow\cE\longrightarrow Br(\cE)\longrightarrow\Z\longrightarrow 0.
\end{equation}
coming from (\ref{i = r-1 4TES}).
So instead, we shall use Proposition \ref{t-hofame} to construct a splitting of (\ref{eq: splittable}) by sending the generator of $\Z$ to $c(\cG)\in Br(\cE)$.
This has the effect of fixing the isomorphism $Br(\cE)=\cE\times\Z$.

Now that we have fixed this isomorphism, we can construct gerbes representing each element of $Br(\cE)$.

\begin{proposition}
The class $(z,n)\in\cE\times\Z= Br(\cE)$ is represented
by the gerbe
\[
\cG^{\otimes n-1}\otimes z_*(\cG),
\]
where $\cG$ is the gamma gerbe, and $z_*(\cG)$ is its translate by $z$.
\end{proposition}

\begin{proof}
Let $c\in\cE\times\Z$ denote the class of $\cG$.
By Proposition \ref{t-hofame}, the $\Z$ component of $c$ is $1$,
and by the above convention, the $\cE$ component of $c$ is $0$.
By Theorem \ref{thm:HOan}, the $\Z$ component of $z_*(c)$ is then also $1$, and its $\cE$ component is $z$.
\end{proof}



\section{The characteristic class
of the gamma gerbe}\label{sec-charclass}
Let as above
$G=\ISL_3(\mathbb Z)=\SL_3(\mathbb Z)\ltimes \mathbb Z^3$ act on $X$, the
total space of the dual tautological bundle $O(1)$ on $Y=\CP^2-\RP^2$. The Dixmier--Douady class
$c(\mathcal G)$ of the
gamma gerbe $\mathcal G$
is a class in the equivariant cohomology $H_G^3(X,\mathbb Z)$. In this section we show
that $H_G^3(X,\mathbb Z)/\mathrm{torsion}=\mathbb Z\oplus\mathbb Z$ and that the image of
$c(\mathcal G)$ is a primitive vector in this group. In particular, this implies that the gerbe is topologically non-trivial.
We first compute the restriction
 of $c(\mathcal G)$ to a generic fibre. We then
compute $H^3_G(X,\mathbb Z)$ relating it to
group cohomology via a spectral sequence.

\subsection{The restriction to a generic fibre}
The bundle projection
 $p\colon X\to \CP^2-\RP^2$ is $\ISL_3(\mathbb Z)$ equivariant (with trivial action of
$\mathbb Z^3$ on the base). The isotropy group of the generic fibre $F=p^{-1}(x)\cong \mathbb C$ is $\mathbb Z^3$. The action
of $n\in Z^3$ is $z\to z+\sum n_ix_i$ and here
generic means that the $x_i$ are linearly independent over $\mathbb Q$.
Correspondingly, the
restriction to $F$ induces
a homomorphism $r\colon H^3_G(X,\mathbb Z)\to
H^3_{\mathbb Z^3}(F,\mathbb Z)$. It is the
composition of the natural homomorphisms
(restriction to a subgroup and pull-back by
the inclusion map)
\[
H^3_G(X,\mathbb Z)\to H^3_{\mathbb Z^3}(X,
\mathbb Z)\to H^3_{\mathbb Z^3}(F,\mathbb Z).
\]
Since $F$ is contractible, $H^\bullet_{\mathbb Z^3}(F,\mathbb Z)=H^\bullet(\mathbb Z^3,\mathbb Z)=
\wedge(\mathbb Z^3)$. Then Prop.~\ref{t-hofame} gives:

\begin{proposition}\label{prop-ho ancora fame}
The restriction to the generic fibre
\[
r\colon H^3_G(X,\mathbb Z)\to
H^3_{\mathbb Z^3}(F,\mathbb Z)=\mathbb Z,
\]
sends
the Dixmier--Douady class $c(\mathcal G)$
to a generator. In particular $c(\mathcal G)$
is non-trivial.
\end{proposition}

\subsection{The topology of $X$}

\begin{proposition}\label{Plinius}
The manifold $X$  retracts to the 2-sphere $S^2$, embedded as the rational
curve $x_1^2+x_2^2+x_3^2=0$, $w=0$. The retraction
$r\colon X\to S^2$ is
$\ISL_3(\mathbb Z)$-equivariant
for the following action on $S^2$: $\mathbb Z^3$ acts trivially; if $S^2$ is viewed as the space of row vectors of norm 1
in $\mathbb R^3$
then $g\in \SL_3(\mathbb Z)$ acts on $v\in S^2$ by $g\cdot v= vg^{-1}/|vg^{-1}|$.
\end{proposition}

\begin{proof} Clearly $X$ retracts to its zero section $Y$. Since the
subgroup $\mathbb Z^3$ acts on the fibres of
$X\to Y$, the retraction is $\ISL_3(\mathbb Z)$-equivariant for the trivial action of $\mathbb Z^3$ on the zero section.
It is thus sufficient to show that we have
an embedding $j\colon S^2\to Y$ and
an $\SL_3(\mathbb Z)$-equivariant
deformation retraction $r\colon Y\to S^2$.

The rational curve $x_1^2+x_2^2+x_3^2=0$ in
$\CP^2$ does not have any real points so it
is contained in $Y$. Thus we have
an embedding $j\colon S^2\cong \CP^1\to
Y$, induced by the $\mathbb C^\times$-equivariant map $\mathbb C^2-\{0\}\to \C^3-\{0\}$
\[
j(z_1,z_2)=(i(z_1^2-z_2^2),z_1^2+z_2^2,2iz_1z_2).
\]
To construct the retraction $r\colon Y\to S^2$,
notice that
to each point $[x]\in\CP^2-\RP^2$,
$x\in\mathbb C^3$ there corresponds an ordered
pair $(\mathrm{Re}\,x,\mathrm{Im}\,x)$ of
linearly independent vectors in $\mathbb R^3$, and all such pairs arise in this fashion. Moreover the oriented
plane defined by $[x]$ is independent of the choice of representative $x$. So we have a map $r$ from $Y$ to the manifold
of oriented planes through the origin
 in $\mathbb R^3$. This manifold is identified with $S^2$ via the normal vector map: in terms of the cross product,
\[
r([x])=
\frac
 {\mathrm{Re}\,x \,\times \mathrm{Im}\,x}
{|\mathrm{Re}\,x \,\times \mathrm{Im}\,x|}.
\]
The induced action of $\mathit{\SL}_3$ on
$S^2$ is obtained from the fact that
the cross product is an isomorphism of
$\SL_3(\mathbb Z)$-modules
$\wedge^2(\mathbb R^3)\to (\mathbb R^3)^*$.

A simple calculation shows that $r\circ j=\mathrm{Id}$ if we identify $S^2$ with $\CP^1$ via
the south pole stereographic projection
$(x,y,z)\to[x+iy,1+z]$. The map $r
\colon Y\to S^2$ is a fibre bundle with
contractible fibres: the fibre over $a$ is the
space of oriented bases of the plane normal to
$a$ modulo rotations and rescaling.
Thus $r\circ j$ is homotopic to the identity.
\end{proof}

\subsection{The cohomology of $\ISL_3(\mathbb Z)$}
\begin{proposition}\label{p-Ampliatus} The first few integral cohomology groups of $\ISL_3(\mathbb Z)=\SL_3(\mathbb Z)\ltimes\mathbb Z^3$ are
\[
H^j(\ISL_3(\mathbb Z),\mathbb Z)=\left\{
\begin{array}{cl}
\mathbb Z,&j=0,\\
0,&j=1,2,\\
\end{array}
\right.
\]
and there is an exact sequence
\[
0\to \mathbb Z\to
H^3(\ISL_3(\mathbb Z),\Z)/\mathrm{torsion}\stackrel{
\mathrm{res}}\longrightarrow H^3(\Z^3,\Z)\cong \Z,
\]
where $\mathrm{res}$ is induced by the
 restriction to the subgroup $\Z^3\subset \ISL_3(\Z)$ and the group $\Z$ on the left is
$H^2(\SL_3(\mathbb Z),\mathbb Z^3)/\mathrm{torsion}$. (We will
also show below that $\mathrm{res}$
is surjective).
\end{proposition}
\begin{proof}
The Lyndon--Hochschild--Serre spectral sequence of the exact sequence of groups
\[
0\to \mathbb Z^3\to \ISL_3(\mathbb Z)\to \SL_3(\mathbb Z)\to 1,
\]
converging to $H^\bullet(\ISL_3(\mathbb Z),\mathbb Z)$,
has $E_2$ term
\[
E_2^{p,q}=H^p(\SL_3(\mathbb Z), H^q(\mathbb Z^3,\mathbb Z)).
\]
 Now $H^q(\mathbb Z^3,\mathbb Z)=\wedge^q(\mathbb Z^3)$ with the natural
action of $\SL_3(\mathbb Z)$ on column vectors.

The relevant part of $E_2$ is ($n\mathbb Z_p$ denotes the direct sum of $n$ copies of $\mathbb Z/p\mathbb Z$)
\[
\begin{array}{cccccc}
q\geq 4 &0&0&0&0& \\
q=3 &\mathbb Z &         &             &           & \\
q=2 &     0    &   0     &\mathbb Z\oplus
\mathbb Z_{2^r}    &&            \\
q=1 &     0    &   0     &\mathbb Z\oplus
\mathbb Z_{2^r}    & &           \\
q=0 &\mathbb Z &   0     &    0        & 2\mathbb Z_2&2\mathbb Z_3\oplus2\mathbb Z_4
      \\
 &p=0       &p=1      &p=2          &p=3         &p=4
\end{array}
\]
for some integer $r$.
In the lowest row we have the groups
$E_2^{p,0}=H^p(\SL_3(\mathbb Z),\mathbb Z)$
computed by Soul\'e (see \cite{Sou}, Theorem 4).
The
$\SL_3(\mathbb Z)$-modules $\mathbb Z^3$ and $\wedge^2(\mathbb Z^3)\cong (\mathbb Z^3)^*$ are related
by an automorphism of $\SL_3(\mathbb Z)$ (inversion composed with transposition). Therefore $E_2^{p,1}\cong E_2^{p,2}$.
 As for $E_2^{p,1}$, we have
$H^0(\SL_3(\mathbb Z),\mathbb Z^3)=0$ (there are no invariant column vectors),
$H^1(\SL_3(\mathbb Z),\mathbb Z^3)=0$ (Sah,
\cite{Sah}, Theorem III.6)
and
$H^2(\SL_3(\mathbb Z),\mathbb Z^3)=\mathbb Z\oplus \mathbb Z_{2^r}$ for some $r$ (Hewitt, \cite{Hew}).
Finally $\wedge^3(\mathbb Z^3)\cong \mathbb Z$ so $E_2^{p,3}=E_2^{p,0}$.

{}From this description of the $E_2$ term
we obtain the groups
 $H^j(\ISL_3(\mathbb Z),\mathbb Z)$
for $j=0,1,2$. For $j=3$, higher differentials
should be considered. The group
 $E_\infty^{2,1}$ is the kernel of
$
d_2\colon E_2^{2,1}\to E_2^{4,0}
$.
Since $E_2^{4,0}$ is torsion, this kernel is
infinite cyclic plus torsion. Since $E_2^{3,0}$
is torsion, we have an embedding
$\Z=H^2(\SL_3(\Z),\Z)/\mathrm{torsion}\hookrightarrow H^3(\ISL_3(\mathbb Z),\mathbb Z)/\mathrm{torsion}$. The cokernel of this map is $E_\infty^{0,3}\subset E_2^{0,3}=H^3(\mathbb Z^3,\mathbb Z)\cong \mathbb Z$, which is vulnerable to $d_2,d_3,d_4$. The edge homomorphism $H^3(\ISL_3(\Z),\Z)\to
E_2^{0,3}=H^3(\Z^3,\Z)$ is induced by the restriction map.
\end{proof}
\subsection{The equivariant cohomology of $X$}
\dontprint{The equivariant cohomology of $X$ with integer
coefficients
is the total cohomology of the double complex $C^p(G,\check C^q(\mathcal U,\mathbb Z))$ where $\check C^\cdot(\mathcal U,\mathbb Z)$ is the \v Cech complex of the equivariant good covering $\mathcal U=(V_a)_{a\in\Lambda_{\mathit{prim}}}$.}
There is a spectral sequence converging to the
equivariant cohomology of $X$ with values in the constant sheaf $\mathbb Z$,
alias the integral cohomology of the stack $[X/\ISL_3(\mathbb Z)]$.
Its $E_2$ term is
$E_2^{p,q}=H^p(G,H^q(X,\mathbb Z))$. Since $X$ is connected, $H^0(X,\mathbb Z)=\mathbb Z$ (with trivial action of $G$) and
we have a natural edge homomorphism
\[
E_2^{p,0}=H^p(G,\mathbb Z)\to
H_G^p(X,\mathbb Z).
\]
\begin{proposition}\label{p-Corelia} Let $G=\ISL_3(\mathbb Z)$. Then $H^0_G(X,\mathbb Z)=H^2_G(X,\mathbb Z)=\mathbb Z$ and $H^1_G(X,\mathbb Z)=0$.
The edge homomorphism $H^3(G,\mathbb Z)\to
H^3_G(X,\mathbb Z)$ is surjective with finite kernel. In particular,
\[
H^3_G(X,\mathbb Z)/\mathrm{torsion}\cong H^3(G,\mathbb Z)/\mathrm{torsion}.
\]
\end{proposition}

\begin{proof}
Since $X$ is a deformation retract of $S^2$, we have
$H^0(X,\mathbb Z)= H^2(X,\mathbb Z)=\mathbb Z$ and
$H^q(X,\mathbb Z)=0$ for $q\neq 0,2$. The action
of $\ISL_3(\mathbb Z)$ on these groups is trivial since it is the restriction of an action of
$\ISL_3(\mathbb R)$ and any homomorphism
$G_\mathbb R\to\mathrm{Aut}(\mathbb Z)=\{\pm\mathrm{Id}\}$ from a connected Lie group $G_\mathbb R$ is necessarily trivial.
Thus $E_2^{p,0}=E_2^{p,2}=H^p(G,
\mathbb Z)$ and all other $E_2^{p,q}$ are trivial.
By Proposition \ref{p-Ampliatus}, $H^1(G,\mathbb Z)=H^2(G,\mathbb Z)=0$, so the $E_2$ term is
\[
\begin{array}{cccccc}
\vdots&\phantom{HHH}&&&\vdots&\\
0&0&0&0&0&\cdots \\
\mathbb Z&0&0&H^3(G,\mathbb Z)&H^4(G,\mathbb Z)&\cdots\\
0&0&0&0&0&\cdots \\
\mathbb Z&0&0&H^3(G,\mathbb Z)&H^4(G,\mathbb Z)&\cdots
\end{array}
\]
This implies the claim for $H^0_G$, $H^1_G$ and
gives an exact sequence
\[
0\to H^2_G(X,\mathbb Z)\to \mathbb Z
\stackrel{d_3}{\longrightarrow}H^3(G,\mathbb Z)\to H^3_G(X,\mathbb Z)\to 0.
\]
To complete the proof it is thus sufficient to show that the image of
$H^2_G(X,\mathbb Z)$ in $H^2(X,\mathbb Z)=\mathbb Z$ is non-zero.
We can use de Rham cohomology to prove this.
A 2-form representing the generator of $H^2(X,\mathbb Z)$ can be constructed with
the Fubini--Study K\"ahler form $\omega_2$ on
$\CP^2$. On $\CP^n$ $\omega_n$
 is represented by the $\mathbb C^\times$-basic 2-form on $\mathbb C^{n+1}-\{0\}$
\[
\omega_n=d\,d^c\log\,|x|^2,
\quad d^c=\frac1{4\pi i}(\partial-\bar\partial).
\]
\begin{lemma}\label{l-Corax}
Let $j:\CP^1\hookrightarrow\CP^2$ be the
parametrization of the rational curve $x_1^2+x_2^2+x_3^2=0$ induced by the map $\mathbb C^2\to
\mathbb C^3$
\[
 (z_1,z_2)\to(i(z_1^2-z_2^2),z_1^2+z_2^2,2iz_1z_2).
\]
Then
\[
j^*\omega_2=2\omega_1,
\]
\end{lemma}

\begin{proof}This follows from the identity
$|j(z)|^2=2|z|^4$.\end{proof}

Let us also denote by $\omega_2\in\Omega^2(X)$
the pull-back of $\omega_2\in\Omega^2(\CP^2)$ by the composition $X\to Y\hookrightarrow \CP^2$.
Since $\omega_1$ is normalized so that $\int_{\CP^1}\omega_1=1$, it follows that
the class of
$\frac12\omega_2\in\Omega^2(X)$ is a generator of $\mathbb Z=H^2(X,\mathbb Z)\subset H^2(X,\mathbb R)$.
Consider the double complex $C^{p,q}=C^p(G,\Omega^q(X))$ computing the equivariant cohomology with real coefficients. Denote the two differentials
by $d$, the de Rham
differential, and by $\delta$, the differential of
group cohomology.
 We must show that $\omega_2\in C^{0,2}$ extends to an equivariant cocycle in
$C^{0,2}\oplus C^{1,1}\oplus C^{2,0}$. We have
$d\omega_2=0$ and
\[
\delta\omega_2(g)=\omega_2-(g^{-1})^*
\omega_2=d\psi(g),
\]
where $\psi\in C^{1,1}$ is
\[
\psi(g)=
d^c\log\frac{|x|^2}{|g^{-1}x|^2}
\]
Note that the argument of the logarithm is $\mathbb C^\times$-invariant, so $\psi(g)$ is a well-defined form on $X$. Finally
\[
\delta\psi(g,h)=
d^c\left(
\log\frac{|x|^2}{|g^{-1}x|^2}
-\log\frac{|x|^2}{|(gh)^{-1}x|^2}
+\log\frac{|g^{-1}x|^2}{|h^{-1}g^{-1}x|^2}
\right)=0.
\]
It follows that $(\omega_2,\psi,0)\in C^{0,2}\oplus C^{1,1}\oplus C^{2,0}$ represents a class
in $H^2_G(X,\mathbb R)$ mapping to the non-trivial class $[\omega_2]\in H^2(X,\mathbb R)$.

\end{proof}
\begin{thm}\label{t-sonostufo}
The exact sequence of Proposition
\ref{p-Ampliatus} is also exact on the right and
thus we have an exact sequence
\[
0\to \mathbb Z\to
H^3_G(X,\mathbb Z)/\text{\rm torsion}\to
H^3(\mathbb Z^3,\mathbb Z)\cong \mathbb Z\to 0
\]
The image of the Dixmier--Douady class
$c(\mathcal G)\in H^3_G(X,\mathbb Z)$ of the gamma gerbe
is a generator of $H^3(\mathbb Z^3,\mathbb Z)$.
\end{thm}

\begin{proof}
This follows from the commutativity of the diagram of
restriction maps and edge homomorphisms
\[
\begin{array}{ccccc}
H^3(G,\Z)&\longrightarrow& H^3(\Z^3,\Z)&&
\\
\downarrow&&\downarrow&\searrow&
\\
H^3_G(X,\Z)&\longrightarrow& H^3_{\Z^3}(X,\Z)
&\longrightarrow& H^3_{\mathbb Z^3}(F,\Z)\cong \Z.
\end{array}
\]
The composition in the bottom row sends $c(\cG)$
to a generator of $\mathbb Z$ by
Proposition \ref{p-Ampliatus} and is thus
surjective. The left vertical arrow is surjective
by Proposition \ref{p-Corelia} and the map
$H^3(\Z^3,\Z)\to H^3_{\mathbb Z^3}(F,\Z)$ is an
isomorphism. It follows that the restriction
map $H^3(G,\mathbb Z)\to H^3(\Z^3,\Z)$ is also
surjective.

\end{proof}

\section{Hermitian structure and connective structure}\label{s-v}
In this section we introduce the notion of a hermitian structure
on a gerbe and its canonical connective structure. Then we construct
a hermitian structure on the gamma gerbe and compute its (1,1)-curvature
in $H^1(\cX,\underline\Omega^{1,1}_{cl})$.

\subsection{Hermitian structures on line bundles and gerbes}
Recall that a hermitian line bundle is a holomorphic line bundle with
a hermitian metric on each fibre depending smoothly on the base.
A hermitian line bundle has a canonical connection, the unique connection
preserving the hermitian metric and compatible with the holomorphic structure.
The curvature $F$ of the canonical connection is a closed $(1,1)$-form, and
$\frac i{2\pi}F$ is a de Rham representative of the first Chern class of the
bundle. Suppose that a discrete group $G$ acts on a complex manifold $X$. Then
a $G$-equivariant hermitian line bundle on $X$ is a hermitian line bundle
$L\to X$ with a lift of the action to $L$ preserving the hermitian metric.
In this case the canonical connection is $G$-equivariant, the curvature is
a $G$-invariant (1,1)-form $F$ and $\frac i{2\pi}F$ is an equivariant
de Rham representative of the first Chern class.
In terms of local trivializations, suppose that a holomorphic line bundle
is given by transition functions $\phi_{a,b}\in \cO^\times(V_a\cap V_b)$
with respect to an open covering $\cU=(V_a)_{a\in I}$.
Then a hermitian metric is given by positive smooth functions $h_a$ on $V_a$,
the norms squared of the trivializing sections,
such that $h_b/h_a=|\phi_{a,b}|^{-2}$. The canonical connection on $V_a$
is $d+\theta_a$ where $\theta_a=\partial\log h_a$ and the curvature is the global
(1,1)-form whose restriction to $V_a$ is $F=d\theta_a=\overline\partial\partial\log
h_a$. In the equivariant case, with an invariant open covering $\cU$,
the lift of the action to the line bundle is
given by functions $\phi_a$ on $G\times V_a$, as described in the introduction
and an equivariant hermitian structure obeys additionally $h_a(x)/h_{g^{-1}a}(g^{-1}x)
=|\phi_a(g;x)|^{-2}$ on $G\times V_a$. The curvature $F$ is then $G$-invariant and
$\frac i{2\pi} F$ represents the image in the equivariant
de Rham cohomology of first Chern class.

These notions have a straightforward generalization to (abelian) holomorphic gerbes:

\begin{definition} A {\em hermitian structure}
on a holomorphic gerbe given by line
bundles $L_{a,b}$ on $V_a\cap V_b$ is a collection of hermitian structures $\|\cdot\|_{a,b}$
on $L_{a,b}$ such that the maps $\phi_{a,b,c}$, see \eqref{e-iso1}, are unitary isomorphisms. A {\em hermitian gerbe} is a holomorphic gerbe with a hermitian structure.
\end{definition}
\begin{definition}
An {\em equivariant hermitian structure} on a $G$-equivariant holomorphic gerbe given
by line bundles $L_{a,b},$ $L_a(g)$ with respect to an invariant open covering
is the additional data of hermitian structures
$\|\cdot\|_{a;g}$ on $L_a(g)$ such that also \eqref{e-iso2},
\eqref{e-iso3} are unitary. A holomorphic gerbe
with an equivariant hermitian structure will
be called {\em equivariant hermitian gerbe}.
\end{definition}

The notion of equivalence of hermitian gerbes is patterned on the corresponding
notion for gerbes defined by transition bundles \cite{Chatterjee,hitchin:gerbe}:
an equivariant hermitian gerbe defined with respect to an
invariant open covering induces by restriction equivariant hermitian gerbes for all equivariant refinements of the covering.
Two equivariant hermitian gerbes $\cG,\cG'$ given by line bundles $(L_{a,b}, L_a(g))$ and
$(L'_{a,b},L'_a(g))$, respectively, are equivalent if, possibly after passing to a
common refinement $(V_a)$, there are hermitian line bundles $M_a$ on $V_a$
with unitary isomorphisms $L'_{a,b}\cong M_a\otimes L_{a,b}\otimes M_b^*$, $L'_a(g)\cong M^*_a\otimes L_a(g)\otimes (g^{-1})^*M_{g^{-1}a}$.

By definition, a $G$-equivariant hermitian gerbe on a complex manifold $X$ is the
same as a hermitian gerbe on the stack $\cX=[X/G]$. More generally,
in the language of groupoid presentations introduced in the Appendix,
a hermitian gerbe on a stack is a hermitian structure on the
central extension of a presentation groupoid, compatible with the groupoid
multiplication.

Recall that a connective structure on a holomorphic gerbe given in terms
of line bundles $L_{a,b}$ is given by smooth complex
connections $\nabla_{a,b}$ on $L_{a,b}$
compatible with the structure maps $\phi_{a,b,c}$, see \eqref{e-iso1},
meaning that
$\phi_{a,b,c}\circ\nabla_{a,c}= (\nabla_{a,b}\otimes{\id}+{\id}\otimes\nabla_{b,c})
\circ\phi_{a,b,c}$, on sections of $L_{a,c}$ restricted to $V_a\cap V_b\cap V_c$.
In the equivariant case, we have additional connections $\nabla_{a;g}$ on $L_a(g)$
compatible with the structure maps \eqref{e-iso2}, \eqref{e-iso3}. In the
language of groupoid presentations, a connective structure is a connection
on the $\C^\times$-bundle defining the groupoid central extension compatible
with the groupoid multiplication.

\begin{definition} A {\em compatible connective structure} on an equivariant hermitian gerbe
given by hermitian line bundles $L_{a,b},L_a(g)$
is a connective structure given by connections
$\nabla_{a,b},\nabla_{a;g}$  preserving the hermitian structure and compatible
with the holomorphic structure.
\end{definition}

\begin{thm}\label{thm-eqgerbesconnstr}
An equivariant hermitian gerbe possesses a unique compatible connective structure.
\end{thm}
\begin{proof}
Let the gerbe $\mathcal{G}$ be represented by line bundles $L_{a,b}$, $L_a(g)$.
Since $L_{a,b}$ and $L_a(g)$ are hermitian line bundles, they possess unique canonical compatible connections $\nabla_{a,b}$ and $\nabla_{a;g}$.
It remains to check compatibility of the connections with the isomorphism
\eqref{e-iso1}-\eqref{e-iso3}. This follows immediately from the fact that
these isomorphism  preserve the hermitian structure of $L_{a,b}$ and because tensor products and pull-backs preserve the corresponding canonical connections.
\end{proof}

Recall that equivalence classes of hermitian line bundles on a complex
manifold $\cX=X$ or stack $\cX=[X/G]$ form a group $\widehat{\mathrm{Pic}}(\cX)$,
an extension of the Picard group $\mathrm{Pic}(\cX)=H^1(\cX,\cO^\times)$ of $\cX$,
with respect to tensor product. The curvature $F=F(\nabla_h)$ of the canonical connection $\nabla_h$
of a hermitian structure $h$ on a line bundle $L$ gives a group homomorphism
\begin{eqnarray*}
&\widehat{\mathrm{Pic}}(\cX)\to H^0(\cX,\underline\Omega_{cl}^{1,1})&\\
&(L,h)\mapsto \frac i{2\pi}F
\end{eqnarray*}
to the group of global closed differential forms of type $(1,1)$. The normalization
is chosen so that the class of the image in de Rham cohomology is a representative
of the image of the first Chern class of $L$. In fancier terms, we have a commutative
diagram
\[
\begin{array}{ccc}
\widehat{\mathrm{Pic}}(\cX)        &\longrightarrow &\mathrm{Pic}(\cX)\\
        \downarrow                    &                &   \downarrow       \\
H^0(\cX,\underline\Omega^{1,1}_{cl})&\longrightarrow &H^2_{dR}(\cX,\C).
\end{array}
\]
The first arrow forgets the hermitian structure. The right vertical arrow sends
a line bundle to the image of the first Chern class in de Rham cohomology and
the lower arrow is induced by the inclusion of closed 2-forms in the de Rham
complex of sheaves $\underline\Omega^0\to\underline\Omega^1\to\cdots$ of smooth
complex-valued differential forms, whose hypercohomology is
$H^\bullet_{dR}(\cX,\C)$.

Similarly, the curvatures of the canonical connections on the hermitian
line bundles $L_{a,b}$, $L_a(g)$ form a 1-cocycle $F=(F_{a,b},F_{a;g})\in
C^{0,1}(G,\underline\Omega_{cl}^{1,1})\oplus C^{1,0}(G,\underline\Omega_{cl}^{1,1})$
with values in the sheaf of closed $(1,1)$-forms. We call the class of
$F$ in $H^1_G(X,\underline\Omega_{cl}^{1,1})$ $=$
$H^1(\cX,\underline\Omega_{cl}^{1,1})$ the {\em $(1,1)$-curvature} of the
hermitian gerbe.
Hermitian gerbes also form a group under tensor product. The group of equivalence
classes of hermitian gerbes $\widehat{\mathrm{Br}}(\cX)$ on $\cX$ is an extension
of the Brauer group $\mathrm{Br}(\cX)=H^2(\cX,\cO^\times)$ of equivalence classes
of holomorphic gerbes.  The $(1,1)$-curvature of a tensor product is the
sum of the $(1,1)$-curvatures of the factors. Moreover, the 1-cocycles
of equivalent hermitian gerbes differ by the differential of a 0-cochain,
consisting of the curvatures of the line bundles $M_a$ defining the equivalence.
Thus we have a group homomorphism
\begin{eqnarray*}
&\widehat{\mathrm{Br}}(\cX)\to H^1(\cX,\underline\Omega_{cl}^{1,1})&\\
&(\cG,h)\mapsto \frac i{2\pi}F.
\end{eqnarray*}
The normalization is again chosen so that we have a commutative diagram
\[
\begin{array}{ccc}
\widehat{\mathrm{Br}}(\cX)        &\longrightarrow &\mathrm{Br}(\cX)\\
        \downarrow                    &                &   \downarrow       \\
H^1(\cX,\underline\Omega^{1,1}_{cl})&\longrightarrow &H^3_{dR}(\cX,\C),
\end{array}
\]
with the right vertical arrow mapping a gerbe to the image in de Rham
cohomology of its Dixmier--Douady class.\footnote{Note that $F$ is given by
imaginary (1,1)-forms, so $\frac i{2\pi}F$  is real and we may take the
{\em real} de Rham complex.}

In explicit terms, suppose that a $G$-equivariant holomorphic gerbe $\cG$ is given on an invariant open covering $\cU=(V_a)$ with trivial bundles $L_{a,b}$, $L_{a}(g)$ as
in the case of the gamma gerbe. Then a hermitian structure is given by the
collection $h_{a,b}$, $h_{a;g}$ of the norms squared of the trivializing sections.
These functions form an equivariant cochain
\[
(h^{0,1},h^{1,0})\in C^1_G(\cU,\underline\R_+)=
 C^{0}(G,\check C^1(\cU,\underline \R_+))
\oplus
 C^{1}(G,\check C^0(\cU,\underline \R_+)),
\]
with values in the sheaf of positive smooth functions. The compatibility condition
with the structure maps is then
$
D h=|\phi|^{-2},
$
where $|\phi|^{-2}\in C^2(\cU,\underline \R_+)$ is the cocycle $(|\phi_{a,b,c}|^{-2},
|\phi_{a,b}|^{-2},|\phi_a|^{-2})$.
The canonical connective structure is given by the connections
$\nabla_{a,b}=d+\partial\log h_{a,b}$, $\nabla_{a;g}=d+\partial\log h_{a;g}$.
The $(1,1)$-curvature is then the class of the cocycle $F$ with
\[
F_{a,b}=\bar\partial\partial\log h_{a,b},\qquad
F_{a;g}=\bar\partial\partial\log h_{a;g}.
\]
Finally the complex de Rham cohomology of $\cX$ may be realized as the cohomology
of the complex $\Omega^n_G(X)=\oplus_{p+q=n}C^p(G,\Omega^q(X))$ of group cochains
with values in smooth complex differential forms.
A representative
of the image of the Dixmier--Douady class of the gerbe is obtained by unravelling
the definitions: a {\em curving}, see \cite{Bry93}, \cite{hitchin:gerbe}, is a collection of
$(1,1)$-forms $B_a$ on $V_a$ so that
$B_b-B_a=F_{a,b}$. A curving always exists and can be constructed by a
partition of unity. The {\em 3-curvature} $dB_a$ is then a globally defined
3-form. Also $F_{a;g}-B_a+(g^{-1})^*B_{g^{-1}a}$
is the restriction
to $V_a$ of a globally defined $(1,1)$-form and
\[
\frac i{2\pi}(d B_a,
F_{a;g}-B_a+(g^{-1})^*B_{g^{-1}a},
0,0
)
\in \Omega^3_G(X),
\]
is a cocycle representing the image in $H^3_{dR}(\cX,\C)$ of the Dixmier--Douady class.

\subsection{The case of the theta function}\label{ssec-thetabunherm}
It is a well-known fact that the theta bundle has an equivariant hermitian structure. The norm
squared of a local section $s(z,\tau)$ is $h_2(z,\tau)|s(z,\tau)|^2$. With our conventions,
 \[
h_2(z,\tau)=\exp\left(-2\pi\frac{(\mathrm{Im}\,z)^2}{\mathrm{Im}\,\tau}+2\pi\,\mathrm{Im}(z-\tau/6)\right).
\]
The statement that this function defines an equivariant hermitian structure on the theta bundle
is equivalent to the fact that the smooth function on $\mathbb C\times H_+$
\[
h_2(z,\tau)|\theta_0(z,\tau)|^2
\]
is invariant under $\ISL_2(\mathbb Z)$. Moreover $h_2$ is defined
on $\mathbb C\times(\mathbb C-\mathbb R)$ and obeys $h_2(-z,-\tau)=h_2(z,\tau)^{-1}$. As
a consequence the invariance property extends to the lower half-plane(recall that
we set $\theta_0(z,\tau)=\theta_0(-z,-\tau)^{-1}$).  In homogeneous coordinates we then
have the result:
\begin{proposition}\label{p-Bertie}
Let $p\colon X_2\times \mathbb C\to X_2$ be the $\ISL_2(\mathbb Z)$-equivariant
theta line bundle. Then
\[
\|\zeta\|^2=h_2\left(\frac w{x_2},\frac{x_1}{x_2}\right)|\zeta|^2,\quad \zeta\in p^{-1}(w,x)=\mathbb C,
\]
is an equivariant hermitian structure on the theta bundle. In other words, the function
\[
h_2\left(\frac w{x_2},\frac{x_1}{x_2}\right)\left|\theta_0\left(\frac w{x_2},\frac{x_1}{x_2}\right)\right|^2
\]
is invariant under the action of $\ISL_2(\mathbb Z)$.
\end{proposition}
To complete the picture, we write down the explicit formula for
the  $(1,1)$-form $c_1=\frac{i}{2\pi}\overline{\partial}\partial\log h_2$, which represents the first Chern class of the theta bundle. In
homogeneous coordinates $z=w/x_2$, $\tau=x_1/x_2$,
\begin{align*}
c_1=\frac{i}{2\pi}\overline{\partial}\partial\log h_2=\frac{i}{2\, \mathrm{Im}\,\tau}\left(d z-\frac{\mathrm{Im}\,z}{\mathrm{Im}\,\tau}d \tau\right)\wedge\overline{\left(d z-\frac{\mathrm{Im}\,z}{\mathrm{Im}\,\tau}d \tau\right)}.
\end{align*}
This 2-form is by construction $\ISL_2(\mathbb Z)$-invariant and restricts to the standard normalized volume form on each fibre $\C/(\ZZ+\tau\ZZ)$.

\subsection{A hermitian structure on the gamma gerbe}\label{ssec-hermgamma}
As in Subsection \ref{sec:the gamma gerbe}, let $G=\ISL_3(\Z)=\SL_3(\Z)\ltimes\Z^3$
act on $X=X_3$ and let $\cG$ be the gamma gerbe defined in terms of the
covering $(V_a)_{a\in\Lambda_{\mathit{prim}}}$ by Theorem \ref{t-main}.
By the preceding discussion, an equivariant hermitian structure on the gamma gerbe is
given by a set of
positive smooth functions $h_{a,b}$ on $V_a\cap V_b$ and $h_a$ on $G\times V_a$
(such that $h_{a;g}(y)=h_a(g;y)$)
obeying the conditions:
\begin{gather}
 h_{a,b}(y)h_{b,a}(y)=1,\quad y\in V_a\cap V_b,\label{e-hi}\\
h_{a,b}(y)h_{b,c}(y)h_{c,a}(y)=|\phi_{a,b,c}(y)|^{-2},\quad y\in V_a\cap V_b\cap V_c,\label{e-hii}
\\
{h_{g^{-1}a,g^{-1}b}(g^{-1}y)}{h_b(g;y)}=|\phi_{a,b}(g;y)|^{2}{h_a(g; y)}
,\,
y\in V_a\cap V_b,\label{e-hiii}
\\
h_a(g_1g_2;x)=|\phi_a(g_1,g_2;y)|^{2}h_a(g_1;y)h_{g_1^{-1}a}(g_2;g_1^{-1}y),
\quad y\in V_a,\label{e-hiv}
\end{gather}
for  all $a,b,c\in I, g,g_1,g_2\in G$.

Introduce the function of three complex variables
\[
h_3(z,\tau,\sigma)=\exp\left(
-\frac{2\pi}3 R_3(\mathrm{Im}\,z,\mathrm{Im}\,\tau,\mathrm{Im}\,\sigma)\right),
\]
\[
R_3(\zeta,t,s)=\frac{\zeta^3}{ts}-\frac32\left(\frac1 t+\frac1s\right)\zeta^2+\left(\frac t{2s}+\frac s{2t}+\frac32\right)\zeta-\frac{t+s}4\,.
\]
Note that both $h_2$ and $h_3$ (and their higher $n$ analogues)
can be expressed in terms of multiple Bernoulli
polynomials (see Subsection \ref{ss-tegf})
\begin{equation}\label{e-ObediahSlope}
h_n(z,\tau_1,\dots,\tau_{n-1})=
\exp\left(-(4\pi/n!)B_{n-1,n}(\zeta,t_1,\dots,t_{n-1})\right),
\end{equation}
where $\zeta=\mathrm{Im}\,z$, $t_j=\mathrm{Im}\,\tau_j$. In particular, $R_3=B_{2,3}$.

\begin{thm}\label{thm-bartycrouch}
Fix a framing of $\Lambda_{\mathit{prim}}$ and let $(V_a)_{a\in\Lambda_{\mathit{prim}}}$ be the open covering of $X_3$ of Section \ref{s-iv}.
For $a\neq b\in \Lambda_{\mathit{prim}}$, let
\[
h_{a,b}(w,x)=\prod_{\delta\in F/\mathbb Z\gamma}h_3\left(\frac{w+\delta(x)}{\gamma(x)},\frac{\alpha(x)}{\gamma(x)},\frac{\beta(x)}{\gamma(x)}\right), \quad (w,x)\in V_a\cap V_b,
\]
where $\alpha,\beta,\gamma$ and $F$ are as in Proposition \ref{p-horcrux}, and set $h_{a,a}=1$. Let $a\in\Lambda_{\mathit{prim}}$ and $(\alpha_1,\alpha_2,\alpha_3)$ the basis of $\Lambda^\vee$ associated
to $a$ by the framing and define for $\mu\in\Lambda^\vee=\mathbb Z^3$
\[
h_a(\mu;w,x)=\prod_{j=0}^{\mu(a)-1}h_2\left(\frac{w+j\alpha_1(x)}{\alpha_3(x)},\frac{\alpha_2(x)}{\alpha_3(x)}\right),\quad (w,x)\in V_a.
\]
Extend $h_a$ to $\ISL_3(\mathbb Z)$ by setting $h_a((g,\mu);w,x)=h_a(\mu\circ g^{-1};w,x)$.
Then the $1$-cochain $(h^{0,1},h^{1,0})$, specified by $h^{1,0}=(h_a)_{a\in\Lambda_{\mathit{prim}}}$ and $h^{0,1}=(h_{a,b})_{a,b\in\Lambda_{\mathit{prim}}}$, with values in the sheaf of positive smooth functions, defines a hermitian structure on the gamma gerbe.
\end{thm}

\begin{proof}
We first need to show that $h_{a,b}$ is well-defined, independently of the choice of $\alpha,\beta$. The general
$\alpha,\beta$ obeying the conditions of Proposition \ref{p-horcrux}  are of the form
\[
\alpha=n\alpha_0+m\gamma,\quad \beta=n'\beta_0+m'\gamma,\quad n,n'\in\mathbb Z_{>0}, \quad m,m'\in\mathbb Z,
\]
where $\alpha_0(a)$, $\beta_0(b)$ are minimal (and equal to $\mathrm{mod}(a,b)$).
It is clear that $h_{a,b}$ does
not change if we add to $\alpha$ or $\beta$ a multiple of $\gamma$. If we replace $\alpha_0$ by a multiple
$n\alpha_0$ then $F$ 
 is replaced by the union of $n$ shifts of $F$ 
and the claim reduces to:

\begin{lemma}
$h_3(z,\tau,\sigma)=\prod_{j=0}^{n-1}h_3(z+j\tau,n\tau,\sigma).$
\end{lemma}

\begin{proof}
This is the special case of a family of identities for multiple Bernoulli polynomials following from
the trivial identity of generating series
\begin{eqnarray*}
\sum_{j=0}^\infty B_{r,j}(\zeta,t_1,\dots,t_r)\frac {u^j}{j!}
&=&e^{\zeta u}\prod_{j=1}^r\frac u{e^{t_iu}-1}
\\
&=&\sum_{j=0}^{n-1}
e^{(\zeta+jt_1)u}\frac u{e^{nt_1u}-1}\prod_{j=2}^r\frac u{e^{t_iu}-1}.
\end{eqnarray*}
In our case $r=2$ and comparing the coefficients of $u^3$ gives
\[
R_3(\zeta,t,s)=B_{2,3}(\zeta,t,s)=\sum_{j=0}^{n-1}R_3(\zeta+jt,nt,s),
\]
which is equivalent to the claim.
\end{proof}

The case of a multiple of $\beta_0$ is treated in the same way. This shows that $h_{a,b}$
is well defined.

It is immediate to check that the functions $h_{a,b}$, $h_a$ are indeed smooth and positive on their
domain of definition and that $h_{a,b}h_{b,a}=1$. There remains to prove the conditions \eqref{e-hii}--\eqref{e-hiv}.

{\em Proof of \eqref{e-hii}\/:}
The first condition $h_{a,b}h_{b,c}h_{c,a}=|\phi_{a,b,c}|^{-2}$ are checked directly
 as in the proof of Theorem \ref{t-avedakedavra1}.

Consider first the case of linearly independent $a,b,c$. Then we proceed as in the proof of Theorem \ref{t-avedakedavra1} and use the direction vectors $\alpha,\beta,\gamma$ of the three pairs in the definition of $h_{a,b}$,
$h_{b,c}$ and $h_{c,a}$. The equation then reduces to
\[
R_3\left(\mathrm{Im}(w/x_3),\mathrm{Im}(x_1/x_3),\mathrm{Im}(x_2/x_3)\right)+\text{cycl.}=\mathrm{Im}\,P_3(w,x).
\]
This identity can be easily proved using the following formula.

\begin{lemma}\label{l-JohnBold}
Let $z_1,\dots,z_n,w_1,\dots, w_n$ be complex
numbers such that $\mathrm{Im}(w_i\bar w_j)\neq 0$,
$i\neq j$. Then
\[
\mathrm{Im}\,\frac{z_1\cdots z_n}{w_1\cdots w_n}
=\sum_{j=1}^n
\frac{\mathrm{Im}(z_1/w_j)\cdots\mathrm{Im}(z_n/w_j)}
{\mathrm{Im}(w_1/w_j)\cdots
\widehat{\mathrm{Im}(w_j/w_j)}\cdots
\mathrm{Im}(w_n/w_j)}.
\]
\end{lemma}

\begin{proof}
Apply the Cauchy
residue theorem to the meromorphic
differential
\[
\omega=\frac1{2i}\frac
{(\zeta z_1-\bar z_1)\cdots(\zeta z_n-\bar z_n)}
{(\zeta w_1-\bar w_1)\cdots(\zeta w_n-\bar w_n)}
\frac {d\zeta}\zeta.\]
The left-hand side of the claimed identity is
$-\mathrm{res}_{\zeta=\infty}\omega-\mathrm{res}_{\zeta=0}\omega$ and
the right-hand side is the sum of the residues
at the remaining poles $\zeta=\bar w_i/ w_i$.
\end{proof}

Consider now the case where the vectors $a,b,c$ span a plane.   We can assume that $a,b,c$
are pairwise linearly independent since otherwise the identity reduces to the inversion relation.
Again we follow the pattern of proof
of Theorem \ref{t-avedakedavra1}, for which
we need an analogue for $h_{a,b}$ of the infinite product representation for $\Gamma_{a,b}$.
This is again provided by ``Bernoulli calculus'':

\begin{lemma}\label{l-MrHardin} Let $a,b\in\Lambda_{\mathit{prim}}$ and $\gamma$ be the corresponding direction
vector. As in Section \ref{s-iii} introduce
$C_{+-}(a,b)=\{\delta\in\Lambda^\vee\,|\,\delta(a)>0,\delta(b)\leq 0\}$.
Then for $(w,x)\in V_a\cap V_b$ the series
\[
S_{a,b}(t)=-t^2\sum_{\delta\in C_{+-}(a,b)/\mathbb Z\gamma}\exp\left(t\, \mathrm{Im}\frac{w-\delta(x)}{\gamma(x)}\right),
\]
convergent for $\mathrm{Im}\,t>0$, is holomorphic around $t=0$ and
\[
h_{a,b}(w,x)=e^{-(2\pi/3)S_{a,b}'''(0)}.
\]
\end{lemma}

\begin{proof}
By \eqref{e-ObediahSlope} and the formula for the generating function of multiple Bernoulli polynomials,
$h_{a,b}$ is defined as $\exp(-(2\pi/3)\hat S'''_{a,b}(0))$ where
\begin{eqnarray*}
\hat S_{a,b}(t)
&=&
\sum_{\delta\in F/\mathbb Z\gamma}
\exp\left(t\,\mathrm{Im}\,\frac{w-\delta(x)}{\gamma(x)}\right)\frac{t^2}
{(e^{t\,\mathrm{Im}\,\tau}-1)
 (e^{t\,\mathrm{Im}\,\sigma}-1)}
\\
&=&\sum_{\delta\in F/\mathbb Z\gamma}
\exp\left(t\,\mathrm{Im}\,\frac{w-\delta(x)}{\gamma(x)}\right)\frac{t^2 e^{-t\,\mathrm{Im}\,\sigma}}
{(e^{t\,\mathrm{Im}\,\tau}-1)
 (1-e^{-t\,\mathrm{Im}\,\sigma})},
\end{eqnarray*}
where $\tau=\alpha(x)/\gamma(x)$ and $\sigma=\beta(x)/\gamma(x)$,
Note that $\mathrm{Im}\,\tau<0$ and $\mathrm{Im}\,\sigma>0$ on $V_a\cap V_b$.
So for $\mathrm{Re}\,t>0$ we may expand the geometric series and obtain, as in  the proof of
Proposition \ref{p-horcrux}, that $S_{a,b}(t)=\hat S_{a,b}(t)$. \end{proof}

Using this representation, we proceed as in the proof of Theorem \ref{t-avedakedavra1}: we
may assume that $(a,b)$, $(b,c)$, $(a,c)$ have the same direction vector $\gamma$ and deduce
from $e
C_{-+}(a,c)=C_{-+}(a,b)\sqcup C_{-+}(b,c),
$ that
\[
S_{a,c}(t)=S_{a,b}(t)+S_{b,c}(t).
\]
Taking the third derivative at $t=0$ implies that $h_{a,b}h_{b,c}=h_{a,c}$. This completes the
proof of \eqref{e-hii}.

{\it Proof of \eqref{e-hiii},\eqref{e-hiv}:\/}
Since the other components of the 2-cocycle are not so explicit and depend on the framing, it is
better to reformulate them in terms of norms of local sections.
Let $\|\Gamma_{a,b}\|^2=h_{a,b}|\Gamma_{a,b}|^2$ be the norm squared of the section $\Gamma_{a,b}$ of $L_{a,b}$
(a more precise notation for this norm would be $\|\cdot\|_{a,b}$). It is a smooth function
defined almost everywhere
on $V_a\cap V_b$. Similarly let $\|\Delta_a\|^2=h_a|\Delta_a|^2$. Then the remaining
conditions \eqref{e-hiii}, \eqref{e-hiv} read
\begin{gather}
\|\Gamma_{a,b}(w+\mu(x),x)\|^2\|\Delta_b(\mu;w,x)\|^2\label{equazione1}
=
\|\Gamma_{a,b}(w,x)\|^2\|\Delta_a(\mu;w,x)\|^2,
\\
\|\Delta_a(\mu+\nu;w,x)\|^2
=
\|\Delta_a(\mu;w,x)\|^2\|\Delta_a(\nu;w+\mu(x),x)\|^2,\label{equazione2}
\\
\|\Delta_a(\mu\circ g^{-1};w,x)\|^2
=
\|\Delta_{g^{-1}a}(\mu;w,g^{-1}x)\|^2,
\label{equazione3}
\end{gather}
for all $\mu,\nu\in\Lambda^\vee=\mathbb Z^3$, $g\in \mathrm{Aut}(\Lambda)=\SL_3(\mathbb Z)$.

We start with \eqref{equazione3}. For this it is better to temporarily make the dependence on the
framing explicit: let $\Delta_{a,f}$ denote $\Delta_{a}$ calculated with the framing
$f\colon\Lambda_{\mathit{prim}}\to\{$Bases of $\Lambda^\vee\}$.  Then
$\Delta_{g^{-1}a,f}(\mu\circ g;w,g^{-1}x)=\Delta_{a,g f}(\mu;w,x)$, for the natural
action of $SL_3(\mathbb Z)$ on the set of framings. Thus \eqref{equazione3} follows from
the following remark.

\begin{lemma}\label{l-Proudie}
The normed squared $\|\Delta_{a,f}(w,x)\|^2$ is independent of the framing $f$. Moreover,
$\|\Delta_{a,f}(w+\mu(x),x)\|^2=\|\Delta_{a,f}(w,x)\|^2$ if $\mu\in H(a)$.
\end{lemma}

\begin{proof}
Both $h_{a}$ and $\Delta_{a}$ are defined using the basis $f(a)=(\alpha_1,\alpha_2,\alpha_3)$
assigned to $a$ by the framing. By definition, $\alpha_2,\alpha_3$ are an oriented
basis of $H(a)$ and $\alpha_1(a)=1$. This basis can be changed in two ways:
$\alpha_2,\alpha_3$ can be replaced by another oriented basis of $H(a)$ or
one can add to $\alpha_1$  an integer linear combination of $\alpha_2$ and $\alpha_3$.
The normed squared is a product of factors
\[
h_2\left(\frac{w+j\alpha_1(x)}{\alpha_3(x)},\frac{\alpha_2(x)}{\alpha_3(x)}\right)
\left|\theta_0\left(\frac{w+j\alpha_1(x)}{\alpha_3(x)},\frac{\alpha_2(x)}{\alpha_3(x)}\right)\right|^2,
\]
and a change of basis amounts to an $SL_2(\mathrm Z)\ltimes \mathbb Z^2$-transformation
under which, by Proposition \ref{p-Bertie}, each factor in invariant. Similarly, a shift
of $w$ by a linear combination of $\alpha_1(x),\alpha_2(x)$ does not change the factors of the
product, proving the second statement.
\end{proof}

Let us turn to Equation \eqref{equazione2}. Clearly it is sufficient to prove it for $\nu$
belonging to a basis, which we take to be the one given by the framing. Note that
$\Delta_a(\mu+\nu;w,x)=\Delta_a(\mu;w,x)$ if $\nu\in H(a)$ and the same holds
for $h_a$. Also $h_a(0;w,x)=\Delta_a(0;w,x)=1$. This proves the claim for $\nu\in H(a)$.
It remains to check the identity for $\nu=\alpha_1$. In this case inserting the definitions
gives
\[
\|\Delta_a(\mu+\alpha_1;w,x)\|^2
=
\|\Delta_a(\mu;w,x)\|^2\|\Delta_a(\alpha_1;w+\mu(a)\alpha_1(x),x)\|^2.
\]
This implies the claim by Lemma \ref{l-Proudie} since $\mu(a)\alpha_1=\mu\mod H(a)$.

Finally, let us prove \eqref{equazione1}. Again, it is sufficient to prove it for $\mu$ belonging
to a basis of $\Lambda^\vee$. Let $\gamma=\gamma_{a,b}$ be the direction vector of $(a,b)$ and
let $\alpha,\beta$ be such that $\alpha(b)=\beta(a)=0$ and that $\alpha(a),\beta(b)$
are positive and minimal
(thus both equal to $\mathrm{mod}(a,b)$). If $\mu=\alpha$ or $\beta$, formulae are
simple: if $\mu=\alpha$, then $\Delta_b(\mu;w,x)=1$ and by Proposition \ref{p-horcrux}
\begin{equation}\label{e-DrGrantly}
\frac{\Gamma_{a,b}(w+\alpha(x),x)}{\Gamma_{a,b}(w,x)}=\prod_{\delta\in F}
\theta_0\left(\frac{w+\delta(x)}{\gamma(x)},\frac{\beta(x)}{\gamma(x)}\right).
\end{equation}
Suppose that the framing assigns to $a$ a basis of the form $(\alpha_1,\beta,\gamma)$. Then
the right-hand side differs from $\Delta_a(\alpha; w,x)$ by a shift of the first argument of
the theta function by a multiple of $\beta(x)/\gamma(x)$. On the other hand we have, thanks
to the identity \[h_3(z+\tau,\tau,\sigma)=h_2(z,\sigma)h_3(z,\tau,\sigma),\]
following from \eqref{e-difference},
\begin{equation}\label{e-MrsGrantly}
\frac{h_{a,b}(w+\alpha(x),x)}{h_{a,b}(w,x)}=\prod_{\delta\in F}
h_2\left(\frac{w+\delta(x)}{\gamma(x)},\frac{\beta(x)}{\gamma(x)}\right).
\end{equation}
Multiplying \eqref{e-MrsGrantly} with the absolute value squared of \eqref{e-DrGrantly}
gives
\[
\frac{\|\Gamma_{a,b}(w+\alpha(x),x)\|^2}{\|\Gamma_{a,b}(w,x)\|^2}
=
\|\Delta_a(\mu;w,x)\|^2,
\]
since by Proposition
 \ref{p-Bertie} the right-hand side is invariant under shifts of the first argument of $\theta_0$ by
$\beta(x)/\gamma(x)$ or by change of framing. Similarly,
\[
\frac{\|\Gamma_{a,b}(w+\beta(x),x)\|^2}{\|\Gamma_{a,b}(w,x)\|^2}
=
\|\Delta_b(\mu;w,x)\|^{-2}.
\]
Also, both $\Gamma_{a,b}$ and $h_{a,b}$ are invariant under a shift of $w$ by $\gamma(x)$.
This proves \eqref{equazione1} if $\mu$ is in the sublattice $M$ generated by $\alpha,\beta,\gamma$. Let $\mu_0\in\Lambda^\vee$ be such that $\mu_0(a)=1$
and $0\leq\mu_0(b)<\alpha(a)$. Such a vector is unique modulo $\mathbb Z\gamma$. Then
the quotient $\Lambda^\vee/M$ is cyclic of order $s=\mathrm{mod}(a,b)$ generated by
the image of $\mu_0$. All this can be best seen by going to the normal form,
where $a=e_1$, $b=re_1+se_2$, $\alpha=(s,-r,0)$, $\beta=(0,1,0)$, $\gamma=(0,0,1)$ and
$\mu_0=(1,0,0)$.

The proof is reduced to checking \eqref{equazione1} for $\mu=\mu_0$. Let
$s=\mathrm{mod}(a,b)=\alpha(a)=\beta(b)$ and $r=\mu_0(b)$ be the invariants of the
wedge $(a,b)$. In this
case the formula in the proof of Proposition \ref{p-felixfelicis} reduces to
\begin{equation}\label{e-Gwynne}
\frac{\Gamma_{a,b}(w+\mu_0(x),x)}{\Gamma_{a,b}(w,x)}=
\frac
{
\theta_0\left(\frac {w+\beta(x)}{\gamma(x)},\frac{\beta(x)}{\gamma(x)}\right)
}
{
\prod_{\delta\in F_3}\theta_0\left(\frac {w+\delta(x)}{-\gamma(x)},\frac{\alpha(x)}{-\gamma(x)}
\right)
}\,,
\end{equation}
where $F_3\subset\Lambda^\vee/\mathbb Z\gamma$ is defined by the inequalities
$0\leq\delta(a)<s,0\leq\delta(b)<r$.
On the other hand, shifting $\delta$ by $-\mu_0$ in the product defining $h_{a,b}$ gives
\[
h_{a,b}(w+\mu_0(x),x)=\prod_{\delta\in F'}
h_3\left(\frac{w+\delta(x)}{\gamma(x)},\frac{\alpha(x)}{\gamma(x)},\frac{\beta(x)}{\gamma(x)}\right).
\]
Here $F'$ is the set of $\delta\in\Lambda^\vee/\mathbb Z\gamma$ obeying
$0< \delta(a)\leq s$, $r\leq\delta(b)<s+r$.
For each value of $\delta(a)$ in this range there is a unique $\delta\in F'$. In particular $\delta(a)=s$ only
for $\delta=\alpha+\beta$. On the corresponding factor we apply the equation relating $h_3$
to $h_2$:
\[
\frac{
h_3\left(\frac{w+\beta(x)+\alpha(x)}{\gamma(x)},\frac{\alpha(x)}{\gamma(x)},\frac{\beta(x)}{\gamma(x)}\right)}{h_3\left(\frac{w+\beta(x)}{\gamma(x)},\frac{\alpha(x)}{\gamma(x)},\frac{\beta(x)}{\gamma(x)}\right)}
=h_2\left(\frac{w+\beta(x)}{\gamma(x)},\frac{\beta(x)}{\gamma(x)}\right).
\]
This allows us to replace the range $F'$ by the range $F''$ determined by $0\leq \delta(a)<s$, $r\leq\delta(b)<s+r$. Similarly, the factors corresponding to $\delta\in F''$ with $\delta(b)\geq s$ may
be related by a shift of $w$ by $\beta(x)$ to factors with $0\leq \delta(b)<r$. The result is
\begin{equation}\label{e-Arabin}
\frac{h_{a,b}(w+\mu_0(x),x)}{h_{a,b}(w,x)}=
\frac
{
h_2\left(\frac {w+\beta(x)}{\gamma(x)},\frac{\beta(x)}{\gamma(x)}\right)
}
{
\prod_{\delta\in F_3}h_2\left(\frac {w+\delta(x)}{-\gamma(x)},\frac{\alpha(x)}{-\gamma(x)}
\right)
}\,,
\end{equation}
As above, we obtain the claim by comparing \eqref{e-Gwynne} with \eqref{e-Arabin}.
\end{proof}
\subsection{The (1,1)-curvature}
Here is an explicit formula for the (1,1)-curvature $F=\bar\partial\partial\log\,h$.
The components $h^{0,1},h^{1,0}$ are given in terms of the functions $h_3$, $h_2$,
respectively.
Thus we first compute $\overline{\partial}\partial\log h_3$, $\overline{\partial}\partial\log h_2$. The latter was already computed in Subsection~\ref{ssec-thetabunherm}.
The former is given by the following expression:
\begin{align*}
\overline{\partial}\partial\log h_3&=\frac{2\pi}3\left[\left(\frac{6\zeta}{st}-3\left(\frac{1}s+\frac{1}t\right)\right)dz\wedge d\overline{z}-\right.\\
&\phantom{=}\left.-\left(-\frac{3\zeta^2}{t^2s}+\frac{3\zeta}{t^2}+\frac{1}{2s}-\frac{s}{2t^2}\right)\left(dz\wedge d\overline{\tau}-d\overline{z}\wedge d\tau\right)-\right.\\
&\phantom{=}\left.-\left(-\frac{3\zeta^2}{ts^2}+\frac{3\zeta}{s^2}+\frac{1}{2t}-\frac{t}{2s^2}\right)\left(dz\wedge d\overline{\sigma}-d\overline{z}\wedge d\sigma\right)-\right.\\
&\phantom{=}\left.-\left(\frac{2\zeta^3}{t^3s}-\frac{3\zeta^2}{t^3}+\frac{s\zeta}{t^3}\right)d\tau\wedge d\overline{\tau}-\right.\\
&\phantom{=}\left.-\left(\frac{\zeta^3}{t^2s^2}-\left(\frac{1}{2s^2}+\frac{1}{2t^2}\right)\zeta\right)\left(d\tau\wedge d\overline{\sigma}-d\overline{\tau}\wedge d\sigma\right)-\right.\\
&\phantom{=}\left.-\left(\frac{2\zeta^3}{ts^3}-\frac{3\zeta^2}{s^3}+\frac{t\zeta}{s^3}\right)d\sigma\wedge d\overline{\sigma}\right].
\end{align*}
By Theorem \ref{thm-bartycrouch}, the components $F_{a,b}$ of the (1,1)-curvature are
then
\[
F_{a,b}=\overline{\partial}\partial\log h_{a,b}=
\sum_{\delta\in F/\mathbb{Z}\gamma}\Phi_{a,b,\delta}^*(\overline{\partial}\partial\log h_3),
\]
where $\Phi_{a,b,\delta}$ is the map from $V_a\cap V_b$ to $\mathbb{C}\times H_{+}\times H_{+}$, defined via $(w,x)\mapsto \left(\frac{w-\delta(x)}{\gamma(x)},\frac{\alpha(x)}{\gamma(x)},\frac{\beta(x)}{\gamma(x)}\right)$.
As for the other components of $F$, we have
\[
F_{a;g}=\overline{\partial}\partial\log h_{a;g}=\sum_{j=0}^{\mu(a)}\Phi_{a,j}^*\overline{\partial}\partial\log h_2,
\]
where the notations are as in Theorem~\ref{thm-bartycrouch}, an explicit expression for the closed $(1,1)$-form $\overline{\partial}\partial\log h_2$ can be found in Subsection~\ref{ssec-thetabunherm}, and the map $\Phi_{a,j}$ from $V_a$ to $C\times H_{+}$ is defined via $(w,x)\mapsto \left(\frac{w+j\alpha_1(x)}{\alpha_3(x)},\frac{\alpha_2(x)}{\alpha_3(x)}\right)$, for any choice of a framing $f$.
We notice that, despite being present in the Definition of $\Phi_{a,j}$, the differential form $\overline{\partial}\partial\log h_{a;g}$ does not depend on the choice of a framing $f$, since a change of framing corresponds to multiplying $h_{a;g}$ by
the absolute value squared of a holomorphic function.

\appendix
\section{Gerbes on stacks}\label{s-gerbes}

Gerbes have been extensively studied in algebraic
geometry. Gerbes on manifolds have also been studied
for example in \cite{Bry93} \cite{hitchin:gerbe}
\cite{murray:gerbe}. They were recently studied over differentiable stacks \cite{bx1}.
In this section, we  summarize and prove
certain necessary facts used in the main part of the paper about
stacks and gerbes on stacks in the differentiable and in the analytic category.

\subsection{Definitions}
First of all, we recall that a stack $\cX$ is a category fibred in
groupoids over the category of smooth manifolds $\cC$  satisfying
two conditions: ``isomorphisms form a sheaf'' and ``effective data
descends''. (See detailed treatments for example in
\cite{bx1} \cite{metzler} \cite{pronk}). A manifold $M$ is a stack by viewing
it as the category of all morphisms from some manifold to it
$\{U\to M\}$. Morphisms between stacks are functors between fibred
categories. A stack morphism $f: \cX \to \cY$ is called a {\em
representable submersion} if for any morphism $Z\to \cY$ where $Z$
is a manifold, the fibre product $\cX \times_{\cY} Z$ is again a
manifold and $\cX\times_{\cY} Z \to Z$ is a submersion. We call
$f$ a {\em representable surjective submersion} if it is
furthermore an epimorphism of fibred categories. We call $f$ {\em
\'etale} if when $\cX\times_{\cY} Z$ is a manifold, the morphism $\cX\times_{\cY} Z \to Z$ is an \'etale map,
namely a local diffeomorphism.

A {\em differentiable stack} $\cX$ is a stack with a representable
surjective submersion from a manifold $X$ to it. We call such a
manifold $X$ a {\em chart} of $\cX$. Then $X\times_{\cX} X \rra X$
is a Lie groupoid whose classifying space is isomorphic to $\cX$
as a stack. We call such a Lie groupoid a {\em groupoid
presentation} of $\cX$.
Different charts give rise to Morita equivalent Lie groupoids.
Morita equivalent Lie groupoids have isomorphic classifying spaces
viewed as stacks. Therefore one can go back and forth between
differentiable stacks and Lie groupoids. Further $\cX$ is an
{\em \'etale stack} if it admits a chart $X$ such that the map $X\to \cX$
is \'etale.

\begin{ep}
Let a group $G$ act on a manifold $M$. The action groupoid
$G\ltimes M \rra M$ is defined by $\bt(g, x) =x$ and $\bs(g, x)=
g^{-1} x$ as  target and source maps, and $(g, x)\cdot (h, y)=
(gh, x)$  when $\bs(g, x)=\bt(h, y)$.  If the action is free and
proper, then the quotient $M/G$ is again a manifold whose
presenting groupoid $M/G\rra M/G$ is Morita equivalent to the
action groupoid $G\ltimes M \rra M$.  If the action is not free
and proper, then the quotient might not be a manifold. But one
could still understand it as a differentiable stack $[M/G]$
presented by the action groupoid. In this language, both manifolds
and orbifolds are differentiable stacks.
\end{ep}

Given an abelian Lie group $A$, an $A$-central extension of a
groupoid $X_1\rra X_0$ is a groupoid $R\rra X_0$ fitting in the
following short exact sequence of groupoids:
\begin{equation}\label{eq:cent-ex}
\xymatrix{
1 \ar[r] &A\times X_0 \ar[d] \ar@<-1ex>[d] \ar[r] & R \ar[d] \ar@<-1ex>[d] \ar[r] & X_1  \ar[d] \ar@<-1ex>[d] \ar[r] & 1 \\
& X_0 \ar[r]^{id} & X_0 \ar[r]^{id} &  X_0
}
\end{equation}
in such a way that $R\to X_1$ is a principal $A$-bundle
over $X_1$. It is equivalent to the following three conditions:
$R|_{X_0}$ is the trivial $A$-bundle; $pr_1^{*} R \otimes
m^{*} R^{-1} \otimes pr_2^* R \cong 1$ via an isomorphism $\phi$
with $pr_i: X_1 \times_{X_0} X_1 \to X_1$ the projections and $m$
the multiplication of $X_1 \rra X_0$; $\phi$ satisfies a higher
coherence condition. Two central extensions $R\rra X_0$ and
$R'\rra X_0$ are isomorphic if all the involved groupoids in
\eqref{eq:cent-ex} are isomorphic correspondingly and all the
expected diagrams commute.

An $A$-{\em gerbe} over $\cX$ is a stack
presented by an
$A$-central extension of a groupoid presentation $X_1\rra X_0$ of
$\cX$. Two gerbes $\cG$ and $\cG'$ over $\cX$ are {\em isomorphic}
to each other if their respective presentations $R\rra X_0$ and
$R'\rra X_0$ are isomorphic as $A$-central extensions of some
groupoid presentation $X_1\rra X_0$ of $\cX$. This goes under the name
of ``stable isomorphism'' in the framework of bundle gerbes of
Murray \cite{murray:gerbe}.

We can generalize the above to the holomorphic setting too. A
complex groupoid is a groupoid in the category of complex
manifolds, namely a groupoid $G_1\rra G_0$ with both $G_1$ and
$G_0$ complex manifolds and structure maps holomorphic maps.
Morita equivalence between two complex groupoids is given by a
complex manifold which is a Morita bibundle in the usual sense but
with holomorphic action from  both groupoids. Then an {\em \'etale
complex stack} $\cX$ is presented by Morita equivalent \'etale
complex groupoids. Here \'etale means that the source map is
\'etale. By definition, an \'etale complex stack has charts $X$'s
which are complex manifolds and the projections $X\to \cX$ are
\'etale. A {\em $\C^\times-$ (holomorphic) gerbe} $\cG$ over an \'etale
complex stack $\cX$ then is a stack presented by a $\C^\times$-central
extension of complex groupoids of a groupoid presentation $X_1\rra
X_0$ of $\cX$. Isomorphisms of $\C^\times$-gerbes are similarly
defined as for gerbes (note that the isomorphisms used are now in
the category of complex groupoids). To relate this to the categorical
interpretation of gerbes, one takes the category of all the torsors
(principal bundles) of the $\C^\times$-central extension. We also
relate
this to the equivariant gerbe described in the introduction in Proposition \ref{prop:relate}.

With the viewpoint of groupoid central extension, some definitions
become short. A {\em hermitian structure} of a $\C^\times$-gerbe
$\cG$ over $\cX$ is a hermitian structure $h$ on the line bundle
$L\to X_1$ associated to $R\to X_1$ such that the isomorphism
$\phi$ sends the tensored hermitian structures from  $pr_1^* L
\otimes m^* L^{-1} \otimes pr_2^*L$ to the trivial hermitian
structure on $\C \times X_1$.  A hermitian structure of $\cG$
comes with a {\em canonical connective structure} which is the canonical
connection 1-form $\theta$ on $R$ associated to $h$, and this is
exactly the unique connective structure constructed in Theorem
\ref{thm-eqgerbesconnstr}.

\subsection{Sheaf and \v Cech cohomology of
stacks}\label{sec:coho}
To describe the Dixmier--Douady class of a gerbe over a
stack, we need to introduce cohomology of stacks. Here we  first
summarize the results in \cite{Del} \cite{Fri82}, then define the \v Cech and sheaf
cohomology  for stacks.

Given a differentiable stack $\cX$ presented by a groupoid
$X_1\rra X_0$, the nerve of the groupoid is a simplicial manifold
$X_\bullet$ given by
\[
\xymatrix{
 {\dots} X_2 \ar[r]\ar@<1ex>[r]\ar@<-1ex>[r] & X_1 \ar[r]\ar@<-1ex>[r] & X_0,} \] where $X_n=X_1\times_{X_0} X_1\cdots
\times_{X_0} X_1$ with $n$-copies of $X_1$ and the
face maps $d_k: X_n \to X_{n-1}$ are given by $d_k(\gamma_1,{\dots},
\gamma_n)$ $=$ $ ( \gamma_1,{\dots}, \gamma_k \gamma_{k+1}, {\dots},
\gamma_n)$ for $k\in [1, n-1]$, $d_0(\gamma_1, {\dots},
\gamma_n)=(\gamma_2, {\dots}, \gamma_n)$ and $d_n(\gamma_1, {\dots},
\gamma_n)=(\gamma_1, {\dots}, \gamma_{n-1})$. There are also
degeneracy maps in the structure of a simplicial manifold but we
will omit them here since they are used only implicitly in our context.

A covering $\cV$ of a simplicial manifold $X_\bullet$ is made up
by open coverings $\cV_n=\{ V_{n, \alpha} \}$ for every manifold
$X_n$ compatible with the structure maps of the simplicial
manifold. More precisely, we
have $n+1$ maps $\partial_k$ from the set of indices of $\cV_n$ to
that of $\cV_{n-1}$ satisfying compatibility conditions as that of
$d_k$'s. Moreover if $x\in V_{n, \alpha}$ then $d_k(x)\in V_{n-1,
\partial_k(\alpha)}$.  To better understand it, we give two
examples.

\begin{ep}\label{1cover}[pull-back coverings]
In the case when $X_\bullet$ is the nerve of  some groupoid $X_1
\rra X_0$, a covering $\cV$ of $X_\bullet$ can be induced by a
covering of $X_0$. Let $\cV_0 = \{ V_{i}\}$ be a covering of
$X_0$.  We define the covering $\cV_n$ of $X_n$ to be the
collection of the following open sets,
 \begin{equation}\label{vn-1}
 V_{i_0 i_1{\dots}i_n}:=\{( \gamma_1, \gamma_2, {\dots}, \gamma_n): \bt(\gamma_1)\in V_{ i_0},
  \bs(\gamma_1) \in V_{i_1}, {\dots}, \bs(\gamma_n) \in V_{i_n}\}
\end{equation}
\entrymodifiers={+++[o][F-]} \SelectTips{cm}{}
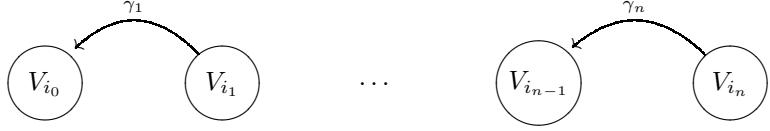
\begin{center}
\xymatrix @+1pc {
 V_{i_0}  &  V_{i_1} \ar@/_2pc/[l]_{\gamma_1} & *{{\dots}} & V_{i_{n-1}} & V_{i_n} \ar@/_2pc/[l]_{\gamma_n} }
\end{center}
It is not hard to see that $\cV_n$ is a covering of $X_n$ and the
map between  indices $\partial_k(i_0{\dots}i_n)=
i_0{\dots}\hat{i}_k{\dots}i_n$ makes $\cV$ all together a covering of
$X_\bullet$. Moreover, $Y_n:=\sqcup_{i_0{\dots}i_n} V_{n, i_0{\dots}i_n}$
is the nerve of the groupoid $Y_1\rra Y_0$. The structure maps of
$Y_\bullet$ are naturally inherited from $X_\bullet$. For
example, $\bt: (g, x) \mapsto x$ and $\bs: (g, x)\mapsto g^{-1}x$.
Moreover $Y_1\rra Y_0$ is Morita equivalent to $X_1\rra X_0$.(see
also the proof of Proposition \ref{prop:h2}).
\end{ep}

\begin{ep}\label{2cover}[invariant coverings]
In the case when $X_\bullet$ is the nerve of  some action groupoid
$G\ltimes X_0 \rra X_0$, there is another covering of $X_\bullet$
induced  by an invariant covering $\cV_0$ of $X_0$. By invariance,
we mean that there is an action of $G$ on the set of indices
$\{a\}$, such that $gV_{a}:=\{ gx: x\in V_{a}\}$ is equal to
$V_{ga}$. A covering $\cV$ of $X_\bullet$ can be induced by a
covering of $X_0$. An element in $X_n$ can be written as $(g_1,
g_2, {\dots}, g_n, x)$ representing $((g_1, x), (g_2, g_1^{-1} x),
{\dots}, (g_n, g_{n-1}^{-1}\cdots g_1^{-1} x )) $    in $(G\ltimes
X_0)\times_{X_0}\cdots\times_{X_0} (G\ltimes X_0)$. Then we define
$\cV_n$ to be the collection of the following open sets,
\[ V_{ g_1, {\dots} ,g_n ,a} := \{ (g_1, {\dots}, g_n, x): x\in  V_{a} \}. \] 
It is easy to see that $\cV_n$ is a covering of $X_n$
and the map
\[
\partial_k(g_1,{\dots},g_n ,a)=
  \begin{cases}
  (g_2, {\dots}, g_n, g_1^{-1} a) & k=0 \\
  (g_1, {\dots}, g_k g_{k+1}, {\dots}, g_n, a) & 0<k<n, \\
  (g_1, {\dots}, g_{n-1} , a)  & k=n,
  \end{cases}
\]
altogether makes $\cV$ a covering of $X_\bullet$. We notice that
the intersection $V_{ g_1, {\dots} ,g_n ,a} \cap V_{h_1, {\dots}, h_n, b}$
is non-empty only if $(g_1,{\dots}, g_n)=(h_1,{\dots},h_n)$ and then the
intersection as a set equals to $V_a\cap V_b$. We also notice that
this notion of covering  works well in the case when $G$ is a discrete
group, but not for a general Lie
group. Moreover, in the first case, such an invariant covering
always exists. For example, we can take a covering of $V_a$ of
$X$, then translate it around and define $V_{ga}$ as $gV_a$.

\end{ep}

Given such a simplicial covering $\cV$ of $X_\bullet$, we can
define a double complex \begin{equation} \label{double-cx}
\begin{CD}
&  & @A\check{\delta} AA  @A \check{\delta} AA @A \check{\delta} AA
{\dots}\\
V_{..}\cap V_{..} \cap V_{..} & \quad  & C^{0,2} @>\delta>> C^{1,2} @>\delta>> C^{2, 2} @>\delta>>{\dots}\\
&  & @A\check{\delta} AA  @A \check{\delta} AA @A \check{\delta} AA
{\dots}\\
V_{..}\cap V_{..} & & C^{0,1} @>\delta>> C^{1,1} @>\delta>> C^{ 2,1} @>\delta>>{\dots} \\
&  & @A \check{\delta} AA  @A \check{\delta} AA @A \check{\delta} AA
{\dots}\\
V_{..} & & C^{0,0} @>\delta>> C^{1,0} @>\delta>> C^{2,0} @>\delta>>{\dots}\\
  & & X_0 & & X_1 &  & X_2
\end{CD}
\end{equation}
where
\begin{equation} \label{eq:c-in-double-cx}
\begin{split}
C^{n, k}:= & \{ \text{global invertible holomorphic functions on
the} \\
& (k+1) \text{-intersections of} \; V_{n, \alpha} \}.
\end{split}
\end{equation}
The columns are \v{C}ech cochains on $X_n$ with covering
$\cV_n$ and
 \v{C}ech differential $\check{\delta}$, and the rows are simplicial cochains, so that the horizontal differential
is defined by
\begin{equation} \label{eq:delta}
(\delta f)_{\cap_{i}V_{n, \alpha_i}}=\prod_{j=0}^{n}
(d_j^*(f_{\cap_i V_{n-1, \partial_j (\alpha_i)}}))^{(-1)^{j+n}}.
\end{equation}
As usual, the cohomology $H^\bullet_{\cV}(X_\bullet, \cO^\times)$
of the double complex is defined as the cohomology of the total
complex $C^{N}= \oplus_{n+k=N} C^{n, k}$ with the total
differential $D=\sum_{n+ k=N}(\delta^{n, k} +(-1)^n
(\check{\delta}^{n, k})): C^{N} \to C^{N+1}$ and $\cO^\times$ is
interpreted as a sheaf of additive groups with $+$ the
multiplication. In fact, $\cO^\times$ can be replaced by any sheaf
$\cF$ of abelian groups over the stack $\cX$. In fact, $\cF$
induces sheaves $\cF^n$ on $X_n$ for all $n$ that automatically
satisfies the expected compatibility conditions \cite{Del} 5.1.6, so that $\cF^\bullet=
(\cF^n){n\geq 0}$ is a simplicial sheaf on the simplicial manifold $X_\bullet$. Then one
replaces invertible holomorphic functions by sections of
$\cF^n$ in the definition of $C^{n, k}$ and everything else
remains the same.

The {\em sheaf cohomology} $H^i(X_\bullet ,\cF^\bullet)$ is
defined as the derived functor of the invariant section functor
$\Gamma_{inv}(X_\bullet, \cF^\bullet):= \ker(\Gamma(X_0, \cF^0) \rra
\Gamma(X_1, \cF^1))$ (\cite{Del} 5.2). Here $\Gamma()$ denotes the set of global sections.
More precisely, we take an injective resolution $I^{\bullet, i}$ of the
simplicial sheaf $\cF^\bullet$,
then $H^i(X_\bullet, \cF^\bullet)$ is the $i$-th cohomology of the
complex
\[ 0\to \Gamma_{inv}(X_\bullet, I^{\bullet, 0}) \to  \Gamma_{inv}(X_\bullet, I^{\bullet, 1})
\to  \Gamma_{inv}(X_\bullet, I^{\bullet, 2}) \to {\dots} \] It can be
also calculated via an acyclic resolution $K^{p, i}$ of $\cF^{p}$,
namely a resolution so that $H^{\geq 1}(X_p, K^{p, i})=0$. Then $H^*(X_\bullet,
\cF^\bullet)$ is the total cohomology of the complex $\Gamma(X_p,
K^{p, q})$. As shown in \cite{haefliger-gpdme}, the sheaf cohomology is
invariant under Morita equivalence. We define $H^i(\cX, \cF):=
H^i(X_\bullet, \cF^\bullet)$ for an \'etale groupoid presentation
$X_1\rra X_0$. For face maps $d^n_{j}:$ $X_{n}$$\to$$X_{n-1}$ for
$j=0, {\dots}, n$ and degeneracy maps $s^n_i:$ $X_{n}$$ \to $$X_{n+1}$
for $i=0, {\dots}, n$, we have $d_j^{n+1}
s^n_{i}$=$s_i^{n-1} d^n_{j-1}$ for $0\leq i \leq j-2$, $d^{n+1}_{i+1}
s^{n}_i=id$=$d^{n+1}_i  s^n_i$ for $0\leq i \leq n$, and
$d^{n+1}_j s_i^n=s^{n-1}_{i-1}d^n_j$ for $j+1 \leq i \leq n$.
Therefore $\delta^n  s^{n-1}+s^{n}\delta^{n+1}=0$ for $n\geq 1$, where
$s^n $$:= $$\sum_i (-1)^i s^n_i$. So the cohomology groups $H^i$
of the complex ${\dots}\overset{\delta}{\to}$ $\Gamma(X_n, \cF^n)$
$\overset{\delta}{\to}$ $\Gamma(X_{n+1}, \cF^{n+1})${\dots} vanish
when $i\leq 1$ for any $\cF^\bullet$. Hence, $H^i(X_\bullet,
\cF^\bullet)$ can be also calculated as the total cohomology of the double
complex $\Gamma(X_p, I^{p, q})$.

Given a covering $\cV$ of $\cX$, $C^q(\cV_p, \cF^p)$ is a sheaf,
where the sections of $C^q(\cV_p, \cF^p)$ on some open set $U$ are
sections of $\cF^p$ on $U$ restricted on $q+1$-fold intersections
of $V_{p,j}$'s. We call $C^q(\cV_p, \cF^p)$ the {\em \v Cech resolution}
associated to a covering $\cV$ of $X_\bullet$ with values in  a
sheaf $\cF^\bullet$ of abelian groups. In fact, we have the following complex of
sheaves:
\[ 0\to \cF^p \to C^0(\cV_p, \cF^p) \to C^1(\cV_p, \cF^p) \to {\dots} \]
 Notice that $H^i(X_p, C^{q}(\cV,
\cF^p))=H^i(\sqcup \cap_{j=0,{\dots},q} V_{p, j}, \cF^p)$.  The
covering $\cV$ of $X_\bullet$ is acyclic, that is $H^{\geq 1} (\cap_{j=0, {\dots}, q}
V_{p,j}, \cF^p)=0$, for all $p$, iff  $C^\bullet(\cV,
\cF^\bullet)$ is an acyclic resolution of $\cF^\bullet$,
 that is, $H^{\geq 1}(X_p, C^{q}(\cV_p, \cF^p))=0$ for all $p$. For constant sheaves such as $\ZZ$ or $\C$, if on each level $n$ the intersections $\bigcap_\alpha
V_{n, \alpha}$ are contractible, then $\cV$ is acyclic. For $\cO$
(resp. $\cO^\times$), one needs the intersections to be Stein
(resp. Stein and contractible) to guarantee being acyclic.

Using the \v Cech resolution $C^\bullet(\cV_p, \cF^p)$ for
$\cF^p$, the total cohomology of the double complex
$\Gamma(X_p,C^q(\cV_p, \cF^p))$ is exactly $H^i_{\cV}(X_\bullet,
\cF^\bullet)$. Since every resolution of $\cF^\bullet$ maps to an
injective resolution of $\cF^\bullet$, we have the induced map
$H^i_{\cV}(X_\bullet, \cF^\bullet)\to H^i(X_\bullet,
\cF^\bullet)$. Moreover, this map is an isomorphism, i.e. $H^i_\cV(X_\bullet,
\cF^\bullet)$ $\cong$$H^i(X_\bullet, \cF^\bullet)$ when $\cV$ is acyclic since $C^q(\cV,
\cF^p)$ is an acyclic resolution of $\cF^p$. The coverings $\cV$'s
of $X_\bullet$ form a direct system by defining $\cV  \prec \cU$
($\cV$ is finer than $\cU$) if $\cV_n \prec
\cU_n$(that is for all $V_{n, \alpha} \in \cV_n$ there is a $U_{n, \beta}
\in \cU_n$ such that $V_{n, \alpha} \subset U_{n, \beta}$) for all
$n$ and the maps $\alpha \to \beta$ are compatible with the facial
maps. Then the map $\underset{\lra}{\lim} H^i_{\cV}(X_\bullet,
\cF^\bullet)$$\to$$ H^i(X_\bullet, \cF^\bullet)$ is an isomorphism
when an acyclic covering exists\footnote{One could drop the
condition of existence of an acyclic coverings, see \cite{tu}, but this condition
always holds in the cases we
consider in the main part of the paper.}. Therefore, we
can use \v Cech cohomology to interpret $H^i(\cX, \cF)$.

\subsection{Geometric objects for $H^1$ and $H^2$ and Dixmier--Douady class}
A {\em holomorphic $\C^\times$-principal bundle on an \'etale
complex stack} $\cX$ is defined as an epimorphism
of \'etale complex stacks $\cL\to\cX$
such that for any \'etale chart $X\to \cX$ the pull-back of $\cL$
on $X$, namely $X\times_{\cX} \cL$, is a holomorphic
$\C^\times$-principal bundle on $X$. 
Let $X_1=X\times_{\cX} X$ and $X_0=X$.
The holomorphic
$\C^\times$-principal bundles on $\cX$ are in 1-1 correspondence
with  the $X_1$ equivariant holomorphic $\C^\times$-principal bundles on $X_0$. See \cite{bx1}
for details. We will prove a similar statement for gerbes in
Proposition \ref{prop:relate}.

A holomorphic $\C^\times$-principal bundle on a stack $\cX$ is
determined by a class in $H^1(\cX, \cO^\times)$. We choose a
groupoid presentation of $\cX$ and a covering of the nerve
$X_\bullet$ of the groupoid. Suppose we have a cocycle $(\eta,
\xi) \in C^{0, 1} \oplus C^{1, 0}$ in the double complex. Then
$\check{\delta}(\eta)=1$, $\cd(\xi)=\delta(\eta)$ and
$\delta(\xi)=1$. So $\eta$ corresponds to a $\C^\times$-principal
bundle $L_0$ on $X_0$.  We define an $X_1$ action on it by
\[ \underbrace{(x, \lambda)}_{ V_{0, i}} \cdot \underbrace{\gamma}_{V_{1,\alpha}} = \underbrace{(\bs(\gamma),
\lambda\xi_{\alpha}(\gamma))}_{V_{0,\partial_1(\alpha)}}, \] where
$\partial_0{\alpha}=i$. One can always choose such $\alpha$ for a
$\gamma$ since $d_0(\gamma)=x$. Then $\cd(\xi)=\delta(\eta)$
ensures that the action is independent of choice of the charts of
$(x, \lambda)$ and of $\gamma$. Moreover $\delta(\xi)=1$
ensures that what we define is indeed an action. Moreover the action
obviously commutes with $\C^\times$ multiplication. Therefore
$L_0$ is an $X_1$ equivariant $\C^\times$-bundle over $X_0$ which
corresponds to a $\C^\times$-bundle $\cL$ on $\cX$. This
construction works for a general covering, in particular for the two
coverings mentioned in the examples.

On the other hand, a $\C^\times$-bundle $\cL$ on
$\cX$ gives an $X_1$-equivariant line bundle $L_0$ on
$X_0$ for any groupoid presentation $X_1\rra X_0$ of $\cX$. We
choose a covering on $X_0$, where $L_0$ trivializes. This covering
generates a covering $\cV$ of the nerve of the groupoid
$X_\bullet$ as in Example \ref{1cover}. Then the transition
functions of $L_0$ give $\eta \in C^{0,1}$ and the $X_1$ action
gives $\xi\in C^{1,0}$ as described above. Then the cohomology
class of $(\eta, \xi)$ in $H^1_{\cV}(X_\bullet, \cO^\times)$
converges to an element in $H^1(\cX, \cO^\times)$ and that is the
class corresponding to this equivariant $\C^\times$-bundle.

Recall that a local section of any $\C^\times$-bundle $L$ over an
open set $U$ of the base manifold $X$ is a function $\theta_0$
defined on $\sqcup (V_{ i} \cap U)$ such that $\cd \theta_0 =
\eta$, for a covering $\{V_i\}$ of $X$. A {\em local section} of
$\cL$ is a function $\theta$ defined on an open subset of $\sqcup
V_{0, i}$ such that $\cd \theta = \eta$ and $\delta \theta =\xi$,
namely $D(\theta) = (\eta, \xi)$.

It is not hard to see that the entire construction applies to
holomorphic line bundles over stacks. In fact, holomorphic line
bundles are equivalent to holomorphic $\C^\times$-principal
bundles.

\begin{proposition}\label{prop:h2}
A holomorphic $\C^\times$-gerbe on a stack $\cX$ is characterized by a class
in $H^2(\cX, \cO^\times)$.
\end{proposition}
\begin{proof}
Take a groupoid presentation $X_1\rra X_0$ of $\cX$ and a covering
$\cV_\bullet$ in the form of Example \ref{1cover} of the nerve
$X_\bullet$. To be more consistent with other part of the paper
and reduce the usage of notations, we assume $X_1\rra X_0$ is an
action groupoid, but the whole proof works for general stacks. Let
$(\Phi, \Theta, \Psi)\in C^{0, 2} \oplus C^{1, 1} \oplus C^{2, 0}$
be a 2-cocycle representing the class in $H^2(\cX, \cO^\times)$.
We have,
\begin{equation} \label{2cocycle}
\cd\Phi =1, \quad \delta \Phi =\cd \Theta, \quad \delta \Theta \cdot \cd
\Psi=1 ,\delta \Psi=1.
\end{equation} Let $Z_1:= \sqcup (V_{i} \cap V_{j})$ and
$Y_1:= \sqcup V_{ij}$ and $U_0:= \sqcup V_{i}$. Then $ Y_1 \rra
U_0$ is a groupoid Morita equivalent to $X_1 \rra X_0$ via the
Morita bibundle $U_0 \times_{pr, X_0, \bt }X_1$. In fact,
$(Y_1=U_0 \times_{X_0} X_1 \times_{X_0} U_0) \rra U_0$ is the
pull-back groupoid by the covering map $U_0\to X_0$. Moreover,
$\Phi$ and $\Theta$ are determined by $\Psi$ via the following
equations,
\begin{equation} \label{put-1} \Phi_{i,j,k}(y)=\Psi_{ijk}(1,1;y)^{-1},
\quad \Theta_{ij, i'j'} (g; y) =\Psi_{ijj'}(g, 1; y)
\Psi_{ii'j'}(1, g; y)^{-1}.\end{equation} There is a trivial
$\C^\times$-bundle $L_\Theta$ on $Y_1$ given by the cocycle
$\Theta_{ij, i'j'}^{-1} \Theta_{ij, i'' j''}$ (namely with
$\Theta_{ij, i'j'}$ as the trivialisation function) with respect
to the covering $V_{ ij} \cap V_{i'j'}$ of $V_{ij}$, where
$\Theta_{ij, i'j'}$ is defined on $V_{ij} \cap V_{i'j'}$. Then the
third and fourth equation of \eqref{2cocycle} ensure that
$L_\Theta \rra U_0$ is a $\C^\times$-groupoid central extension of
$Y_1 \rra U_0$ with the multiplication on $L_\Theta$
\begin{equation}\label{eq:mul-cent}
(g, y, \lambda) \cdot (h, g^{-1} y, \mu) = (gh, y, \lambda \mu
\Psi_{ijk}(g, h; y)),
\end{equation}
for global sections $(g, y, \lambda)$ (appearing as
$(g, y, \lambda \Theta_{ij, i'j'})$ on $V_{ij}\cap V_{i'j'}$) over $V_{ij}$ and $(h,
g^{-1} y, \mu ) $ over $V_{jk}$. Therefore this central extension
gives us the $\C^\times$-gerbe $\cG$ on $\cX$.

It is not surprising by \eqref{put-1} that the first and second
conditions in \eqref{2cocycle} are not used. More precisely,  $\cd
\Phi =1$ implies that $\Phi$ gives a $\C^\times$-gerbe $\cG_\Phi$
on $X_0$. In particular, $\Phi$ gives a $\C^\times$-bundle
$L_{\Phi}$ on $Z_1$ by the 1-cocycle $\Phi_{ikk'}
\Phi^{-1}_{jkk'}=\Phi_{ijk'} \Phi_{ijk}^{-1}$ (namely with
$\Phi_{ijk}$ as the trivialisation function) with respect to the
covering $(V_{i} \cap V_{j}) \cap V_{k}$ on each $V_{i} \cap
V_{j}$. Notice that $i: Z_1 \to Y_1$ by $y \in V_i\cap V_j$
mapping to $(1; y) \in V_{ij}$ is an embedding. The second
equation of \eqref{2cocycle} tells us that $i^*L_\Theta =L_\Phi$.

Now we show  that different cocycles representing the same
cohomology class give the same gerbe. If we have another cocycle
representing the same cohomology class as $(\Phi, \Theta, \Psi)$
possibly on a different covering of a different groupoid, then the
corresponding $Y'_1\rra U'_0$ however is  still Morita equivalent
to $Y_1\rra U_0$ via some Morita bibundle which we call $U_0''$.
The pull-back groupoid $Y''_1\rra U''_0$ via $U''_0\to U_0$
viewed as the common refinement of the two coverings is Morita
equivalent to both $Y_1\rra U_0$ and $Y'\rra U'_0$ via groupoid
morphisms. (The conditions on Morita bibundle ensure that the
pull-backs via $U_0''\to U_0$ and $U_0' \to U_0$ give the same
groupoid). Notice that the pull-back groupoid of $R\rra U_0$ on
$Y''\rra U''_0$ is Morita equivalent to $R\rra U_0$, and $(\Phi,
\Theta, \Psi)$ converges to the same cohomology class as its
pull-back cocycle on the covering generated by $U''_0$. We
conclude that we only have to show the statement for two classes
on the same covering of the same groupoid, namely if $(\Phi,
\Theta, \Psi)$ and $(\Phi', \Theta', \Psi')$ differ by a
coboundary, which says
\begin{equation}\label{eq:differ}
\Phi=\Phi' \cdot \cd \eta, \; \Theta =\Theta'\cdot \delta
\eta \cdot (\cd \xi)^{-1}, \; \Psi =\Psi' \cdot \delta \xi,
\end{equation}
then their corresponding central extensions are isomorphic. We
construct an isomorphism $f: L_\Theta \to L_\Theta'$ by defining
it on the local charts $V_{ij}\cap V_{i'j'} \times \C^\times \to
V_{ij}\cap V_{i'j'} \times \C^\times$, with
\[ (g,y, \lambda)\mapsto (g,y,
\lambda \cdot \frac{\eta_{ii'}(y)} {\eta_{jj'}(g^{-1}x)
\xi_{i'j'}(g;x)}).\] The second equation in \eqref{eq:differ}
guarantees that they glue to a global map.

In local charts, the trivialisation of $L_\Theta$ on $V_{ij}\cap
V_{i'j'}$ is given by function $\Theta_{ij, i'j'}$. Using the
third condition in \eqref{2cocycle} and \eqref{eq:mul-cent}, the multiplication of
$L_\Theta$ written in local trivialisation is
\[
(g, y, \lambda ) \cdot (h, g^{-1}y, \mu) = (gh, y, \lambda \mu
\Psi_{i'j'k'}(g,h; y) ),
\]
for $(g, y)\in V_{ij}\cap V_{i'j'}$ and $(h, g^{-1}y) \in
V_{jk}\cap V_{j'k'}$. From this it is easy to see that $f$ is a
groupoid homomorphism. It is not hard to conclude that $f$ is an
isomorphisms of central extensions.

For the converse, a $\C^\times$-gerbe $\cG$ over the stack
$\cX$ can be viewed as a $\C^\times$-central extension of a
groupoid $X_1\rra X_0$ presenting $\cX$. Take a covering $V_i$ of
$X_0$. Then it generates a covering $\cV$ of the simplicial
manifold $X_\bullet$ as in Example \ref{1cover}. By using the same notation as
before, $Y_1\rra U_0$ is the pull-back groupoid also
presenting $\cX$, so the central extension can be pulled back to
this new groupoid and we call it $L\rra U_0$. Then $d_0^*L \otimes
d_1^* L^{-1} \otimes d_2^* L \cong 1$ via a trivialisation
function $\Psi_{ijk} $ over $Y_1\times_{U_0} Y_1$, recalling that
$Y_1\times_{U_0} Y_1 =\sqcup V_{ijk}$ with $V_{ijk}= \{(g,h,y):
y\in V_i, g^{-1}y \in V_j, (gh)^{-1}y \in V_k\}$. By
\eqref{put-1}, $\Psi_{ijk}$ determines a class $[(\Phi, \Theta,
\Psi)]$ in $H^2_\cV(Y_\bullet, \cO^\times)$ and it converges to a
class in $H^2(\cX, \cO^\times)$. Moreover, if we have two
isomorphic $\C^\times$-central extension $R$ and $R'$ of $X_1\rra
X_0$ given by isomorphisms $f_{ij, i'j'}$ locally defined on the
piece $V_{ij}\cap V_{i'j'}$, then their corresponding cocycles
differ by a coboundary $D(\eta, \xi)$ with $\eta_{ij}(y)= f_{ij,
ij}(1, y)$ and  $\xi_{ij}(g, y)= f_{ij, ij}(g, y)$.
\end{proof}

As the proof above shows, $L_\Theta|_{Z_0} \cong L_\Phi$.
Therefore the gerbe given by $L_\Theta$ can be viewed as giving a
$G$-{\em equivariant structure} to the gerbe $\cG_\Phi$ on $X_0$
given by $L_\Phi$.

Using the short exact sequence of sheaves,
\[ \ZZ \hookrightarrow \cO \overset{\exp 2 \pi i} {\rightarrowtail}
\cO^\times, \] we will have a long exact sequence of cohomology
groups. In particular, we have $H^2(\cX, \cO^\times)\to H^3(\cX,
\ZZ)$. According to the above proposition, a gerbe $\cG$ over
$\cX$ corresponds to a class in $H^2(\cX, \cO^\times)$. The image
of this class in $H^3(\cX, \ZZ)$ is called the {\em Dixmier--Douady
class} of $\cG$.

\subsection{Equivariant gerbes and gerbes over stacks}

In this subsection, we relate the two sorts of gerbes: the one
described in the introduction and the one in this section. Let a
group $G$ act on a manifold $X$. Then we have the action groupoid
$G\times X \rra X$. The {\em equivariant cohomology} $H^i_G(X,\cF)$ is
defined as the cohomology of the nerve of the action groupoid
(an equivariant sheaf $\cF$ on $X$ induces a simplicial
sheaf $\cF^\bullet$ on the nerve of the action groupoid).
Therefore we have $H^i_G(X, \cF)=H^i([X/G], \cF)$, where $[X/G]$
is the differentiable stack presented by $G\times X\rra X$.
Moreover, an invariant covering $\cV$ of $X$ induces
 a covering (still denoted by $\cV$) of the nerve
$G^\bullet \times X$ as in Example \ref{2cover}. Then as in Section \ref{sec:coho}, we have a map
$H^i_\cV(G^\bullet\times X, \cF)$$\to$$H^i([X/G], \cF)$. It is an
isomorphism when $\cV$ is an acyclic covering, which is true in our
cases. Recall that a 2-cocycle  $(\phi_{a,b,c}, \phi_{a,b},
\phi_a)$ representing a class in $H^2_\cV(G^\bullet\times X,
\cO^\times)$ corresponds to an equivariant holomorphic bundle
gerbe as described in the introduction. Under
$H^\bullet_\cV(G^\bullet\times X, \cO^\times)$$ \to $$H^\bullet
([X/G], \cO^\times)$, the image of $(\phi_{a,b,c}, \phi_{a,b},
\phi_a)$ corresponds to a $\C^\times$-gerbe on $[X/G]$. In the
following proposition, we construct this gerbe on $[X/G]$
in a concrete geometrical way from $(\phi_{a,b,c}, \phi_{a,b}, \phi_a)$.

\begin{proposition}\label{prop:relate} If we use an invariant covering as in
Example \ref{2cover} for an action groupoid $X_1:=G\times X$ over
$X_0:=X$, then an equivariant holomorphic bundle
gerbe on $X$ defined by a cocycle
$(\phi_{a,b,c}, \phi_{a,b}, \phi_a)$ corresponds to a holomorphic $\C^\times$-gerbe on $[X/G]$ given as a groupoid central extension.
\end{proposition}
\begin{proof}
Take a 2-cocycle  $(\phi_{a,b,c}, \phi_{a,b}, \phi_a)$
representing a class in $H^2_\cV(G^\bullet\times X, \cO^\times)$.
Then $(\phi_{a,b,c}, \phi_{a,b}, \phi_a)$ satisfies the conditions
\begin{equation} \label{2cocycle-equi}
\cd (\phi_{.,.,.}) =1, \;\delta (\phi_{.,.,.})=\cd
(\phi_{.,.}), \; \delta(\phi_{.,.})=\cd(\phi_{.})^{-1},\;
\delta(\phi_{.})=1.
\end{equation}
Denote by $(V_{a})$  the covering of $X$ that generates the
invariant covering $\cV$ as in Example \ref{2cover}. Let
$U_0=\sqcup V_a$, $Y_1:=\sqcup V_{g, a}$ and $Z_1:=\sqcup (V_a\cap
V_b)$. The groupoid $Y_1\rra U_0$ is not Morita equivalent to
$X_1\rra X_0$ any more. But the groupoid ``generated'' by $Y_1$
and $Z_1$ is the pull-back groupoid $U_1:= U_0\times_{X_0} X_1
\times_{X_0} U_0$, therefore Morita equivalent to $X_1\rra X_0$.
More precisely, $Y_1$ and $Z_1$ are subgroupoids of $U_1$ over the
same base $U_0$. For an element $y \in V_{a}\cap V_{b}$ in $Z_1$
and an element $(g, y)\in V_{g, b} $ in $Y_1$, the multiplication
is $y \cdot (g, y)= (g, y)\in W_{ g, a, g^{-1}b}$, where $W_{g, a,
g^{-1}b}=\{ (g, y), y\in V_a, g^{-1} y \in V_{g^{-1}b}\}$. As in
the case of equivariant gerbes, there is a trivial
$\C^\times$-bundle $L_{11}$ on $Y_1$ given by the trivialisation
function $\phi_{a,a'}(g,-)$ with the covering $V_{g, a}\cap V_{g,
a'}$ of $V_{g, a}$,  and a trivial $\C^\times$-bundle $L_{20}$ on
$Z_1$ given by the trivialisation function $\phi_{a,b,c}$  with
the covering $(V_a\cap V_b) \cap V_c$ of $V_a \cap V_b$. The
corresponding holomorphic line bundles of $L_{11}$ and $L_{20}$
are exactly what we called $L_a(g)$ and $L_{ab}$ before. Notice
that $U_1=\sqcup_{g, a, g^{-1}b} W_{g, a, g^{-1}b}$. We define a
$\C^\times$-bundle on $U_1$ piece by piece, namely by $L|_{W_{g,
a, g^{-1}b}}:= \iota_{a, b}^* L_{20} \otimes \iota_{g, b}^*
L^{-1}_{11}$, where $\iota_{a,b}:W_{g, a, g^{-1}b}\to V_{a}\cap
V_{ b}$, by $(g,y)\mapsto y$ and $\iota_{g, b}:W_{g, a,
g^{-1}b}\to V_{ g,b}$, by $(g,y)\mapsto (g,y)$, where the second
piece has to be understood as an element of $V_{g,b}$. Then
equation \eqref{2cocycle-equi} ensures that $L\rra U_0$ is a
$\C^\times$-central extension of $U_1\rra U_0$. More precisely,
$L$ being a $\C^\times$-central extension is equivalent to that
$L|_{U_0}\cong \C^\times\times U_0$,
\begin{equation}\label{eq:cent-ext}
d_0^* L \otimes d_1^* L^{-1} \otimes d_2^* L \cong 1
\end{equation}
on $U_2:=U_1 \times_{U_0} U_1 $, and the isomorphisms satisfy
further coherence condition, where $d_i: U_2\to U_1$ are the
facial maps. Below, we verify these conditions one by one. Since
$\phi_{a,a,c}(y)=1$ and $\phi_{a,a'}(1; y)=1$ by requirements of
the equivariant cocycles, we have $L|_{U_0}=1$. For
\eqref{eq:cent-ext}, we notice that the second (resp. the first
and third) condition in \eqref{2cocycle-equi} says that
$\phi_{b,c}(g;y)$ (resp. $\phi_{a,b,c}^{-1}(y)\phi_{c}(g,h;y)$)
gives the first (resp. second) isomorphism ``$\cong$'' in the
following chain of isomorphisms,
\[
\begin{split}
 d_0^* L \otimes d_2^* L &=  (\iota_{a,b}^*L_{20}\otimes \iota_{g,
b}^*L_{11}^{-1}) \otimes (\iota_{g^{-1} b ,
g^{-1}c}^*L_{20}\otimes \iota_{h, g^{-1}c}^* L_{11}^{-1}) \\
&=\iota_{a,b}^*L_{20}\otimes (\iota_{g, b}^*L_{11}^{-1} \otimes
\iota_{g^{-1} b , g^{-1}c}^*L_{20})\otimes \iota_{h, g^{-1}c}^*
L_{11}^{-1} \\
&\cong \iota_{a,b}^* L_{20} \otimes (\iota_{b,c}^* L_{20} \otimes
 \iota_{g, c}^* L_{11}^{-1}) \otimes \iota_{h,
 g^{-1}c}^*L_{11}^{-1} \\
&=(\iota_{a,b}^* L_{20} \otimes \iota_{b,c}^* L_{20}) \otimes
 (\iota_{g, c}^* L_{11}^{-1} \otimes \iota_{h,
 g^{-1}c}^*L_{11}^{-1}) \\
&\cong \iota_{a, c}^* L_{20} \otimes \iota_{gh, c}^*
 L_{11}^{-1}=d_1^* L,
 \end{split}
\]where we omit certain pull-backs for simplicity.
Therefore we only have to check that the higher coherence, namely
\[\Psi_{a,g^{-1}b,(gh)^{-1}c}(g,h;y):= (\phi_{a,b,c}^{-1}(y)
\phi_{b,c}(g;y) \phi_a(g,h;y))\] satisfies $\delta(\Psi)_{a,
g^{-1}b, (gh)^{-1}c,(ghf)^{-1} d}(g, h, f; y) =1$ on the triple
fold $W_{g, a, g^{-1}b}$ $ \cap $ $W_{h, g^{-1}b,(gh)^{-1} c }$$
\cap $ $W_{f,(gh)^{-1} c, (ghf)^{-1}d }$. It follows from
\eqref{2cocycle-equi} and the following calculation:
\[
\begin{split}
&\delta (\Psi)_{a,
g^{-1}b, (gh)^{-1}c,(ghf)^{-1} d}(g,h,f;y) \\
=& (\delta \phi_{.,.,.})_{b,c,d}^{-1}(g; y) (\cd
\phi_{.})_{c,d}(g, h; y) (\delta \phi_{.,.})^{-1}_{c,d}(g, h; y)
(\cd \phi_{.,.})_{b, c, d} (g; y).
\end{split}
\] Therefore $L$ is a $\C^\times$-central extension and gives a $\C^\times$-gerbe $\cG$ on $\cX$.

Finally, we notice that $V_{a}$, $W_{g, a, g^{-1}b}$ etc. actually
gives a covering as in Example \ref{1cover}, which is  a
refinement of the invariant covering. From the above
construction, the cocycle corresponds to the constructed gerbe
$\cG$ is $(\Phi, \Theta,\Psi)$ with $\Phi$ and $\Theta$ determined
by $\Psi$ as in \eqref{put-1}. It is easy to see that this new
cocycle restricting on the invariant covering is exactly
$(\phi_{a,b,c}, \phi_{a,b}, \phi_{a})$. Therefore the
$\C^\times$-gerbe $\cG$ corresponds to the cohomology class in
$H^2(\cX, \cO^\times)$ represented by $(\phi_{a,b,c}, \phi_{a,b},
\phi_{a})$.

\end{proof}

\subsection{Equivariant cohomology and group cohomology}
\label{sect:coho-group-stack} Since we use the definition of
equivariant cohomology with action groupoid instead of with a
classifying space $EG$, we supply the proof in this setting of a
basic fact of equivariant cohomology we have used.
\begin{thm} \label{thm:coho-group-stack}
Let $G$ be  a discrete group acting (not necessarily freely) on a
manifold $X$. Then for a sheaf $\cF$ on the stack $[X/G]$, namely
an equivariant sheaf $\cF$ on $X$, there is a spectral sequence
\[ E_2^{p, q}=H^p(G, H^q(X, \cF)) \Rightarrow H^{p+q}_G(X, \cF)=H^{p+q}([X/G], \cF).\]
\end{thm}
\begin{proof}
An equivariant sheaf $\cF$ on $X$ is a sheaf with isomorphisms of
sheaves $\phi_g: g_*\cF\to \cF$ such that
$\phi_g\phi_h=\phi_{gh}$. We take an acyclic equivariant
resolution $K^{q}$ of $\cF$. There always exists one (for example
one could use Godement's resolution). Then $K^q$ generate a
simplicial sheaf $K^{\bullet,q}$ on $G^\bullet\times X$. We form a
spectral sequence associated with the double complex
$\Gamma(G^p\times X, K^{p,q})$ of maps $G^p \to \Gamma(X, K^{q})$,
namely  $E_0^{p, q}$$:=$$\Gamma(G^p\times X, K^{p,q})$, filtered
by $\oplus_{p\geq n} E_0^{p, q}$ (in the direction of $p$). On one
hand, the double complex calculates the sheaf cohomology
$H^i([X/G], \cF)$. On the other hand, there is a $G$ action on
$H^i(X, \cF)$ induced by the $G$ action on $X$. In fact, $E_1^{p,
q}=C^p(G, H^i(X, \cF))$ since $K^q$ is acyclic. Hence the first
page $(E_1^{p, q}, d_1)$ is exactly the complex $(C^p(G,  H^i(X,
\cF)), \delta_G)$ to calculate the group cohomology. Therefore
$E_2^{p,q}=H^p(G, H^q(X, \cF))$.
\end{proof}

\subsection{Hypercohomology on stacks}
Hypercohomology of complexes of sheaves on manifolds and action groupoids
is explicitly studied in \cite{Bry93},
\cite{gomi-deligne}. Here we extend it to the
category of complex stacks and relate it to our construction using
invariant coverings. Our treatment is similar to the one in the
references.

We first discuss the concept of hypercohomology of a complex of sheaves on
a complex stack. Take a stack $\cX$, and a complex of sheaves
of abelian groups
\begin{equation}
\cF_0 \to \cF_1 \to {\dots}\to \cF_m
\end{equation}
over $\cX$. Take a presentation $X_1\rra X_0$ of $\cX$ and form
the nerve $X_\bullet$. Then we naturally have a complex of
simplicial sheaves
\begin{equation}\label{eq:sheaf-cx}
\cF_0^\bullet \to \cF_1^\bullet \to {\dots} \to \cF_m^\bullet,
\end{equation}
in the sense of Section 3.2 in \cite{gomi-deligne}. Following
this reference\footnote{Notice that there are enough injective
objects in the category of sheaves of abelian groups over stacks,
as explained in \cite{grothendieck-tohoku}.}, we take an injective resolution of \eqref{eq:sheaf-cx},
which looks like a triple complex with one direction a complex of
injective sheaves $I_{0}^{p, q} \to I_{1}^{p, q} \to  {\dots} \to
I_{n}^{p, q}$ and one direction an injective resolution $I_r^{p,
\bullet}$ of $\cF_r^p$. Furthermore, the maps satisfy similar
compatibility conditions as  in \cite{Bry93} 1.2.6., where the
author defines an injective resolution for a complex of sheaves
over a manifold.  We define the {\em
hypercohomology} $H^i (X_\bullet, \cF_0^\bullet \to {\dots} \to
\cF_m^\bullet)$ as the total cohomology of the double complex
\begin{equation} \label{hyper-cx}
C^{q, r}:=\Gamma_{inv}(X_\bullet, I^{\bullet, q}_r).
\end{equation}
In fact it is also the total cohomology of the triple complex
\begin{equation} \label{hyper-tri-cx}
C^{p, q, r}:= \Gamma(X_p, I_r^{p, q}).
\end{equation}
This is so  as in the case of sheaf cohomology: we filter
$C^{p,q,r}$ by $\oplus_{q+r\geq n}C^{p,q,r}$,
then we get the $E_1$ page with a lot of zeros except $C^{q,r}$.
We take the spectral sequence with  $E^{q,r}_0:=$
$C^{q,r}$=$\Gamma_{inv}(X_\bullet, I^{\bullet, q}_r)$ and filtered
by $\oplus_{r\geq n} C^{q, r}$. Then the first page
$E^{q,r}_1$=$H^q(X_\bullet, \cF_r^\bullet)$=$H^q(\cX, \cF_r)$
which is independent of the choice of groupoid presentation of
$\cX$. Hence the hypercohomology $H^i (X_\bullet,$$
\cF_0^\bullet$$ \to$$ {\dots}$$ \to$$ \cF_m^\bullet)$ is
independent of the choice of the groupoid presentation $X_1 \rra
X_0$ of $\cX$. Therefore, we define the {\em hypercohomology of a
stack} $\cX$ as  $H^i(\cX,$ $ \cF_0$ $\to $ $\cF_1{\dots}$ $\to$ $
\cF_m)$ $:=$ $H^i(X_\bullet,\cF_0^\bullet$ $ \to$ $ {\dots} $
$\to$ $ \cF_m^\bullet)$ for some groupoid presentation $X_1\rra
X_0$ of $\cX$.

Take a simplicial open covering $\cV$ of $X_\bullet$---the nerve of the
groupoid $X_1 \rra X_0$ presenting $\cX$.  Then
we let
\begin{equation}\label{hyper-check-tri-cx}
 {\check{C}}^{p, q, r}= \Gamma(\cap_{i=0}^{q} V_{p, \alpha_i},
\cF_{r}^p).
\end{equation}
Then we have $\delta: \check{C}^{p,q,r}\to
\check{C}^{p+1, q, r}$ defined as in \eqref{eq:delta} and $\cd:
\check{C}^{p,q,r} \to \check{C}^{p, q+1, r}$ the \v Cech
differential. They are defined exactly as the $\delta$ and $\cd$
in the double complex \eqref{double-cx}. Moreover, there is
another differential $d: \check{C}^{p,q,r} \to
\check{C}^{p,q,r+1}$, induced by the differential $\cF_r \to
\cF_{r+1}$ in the sheaf complex. Then $(\check{C}^{p,q,r}, \delta,
\cd, d)$ is a triple complex and the total complex is
$\check{C}^{N} :=\oplus_{N=p+q+r} \check{C}^{p,q,r}$ with the
total differential $D_3= \delta + (-1)^{q} \cd + (-1)^{q+r} d$. We
define the \v Cech hypercohomology ${\check{H}}^\bullet_\cV(X_\bullet,
\cF_0^\bullet \to {\dots} \to \cF_m^\bullet)$ with respect to the covering $\cV$ as
the total cohomology of this triple complex.

Then we have the following proposition which generalizes the case
of the relation between \v Cech cohomology and sheaf cohomology.

\begin{proposition} \label{prop:deligne-cech}
In the above setup, there is a map
\[{\check{H}}^n_\cV(X_\bullet, \cF_0^\bullet \to {\dots} \to
\cF_m^\bullet) \to H^n(\cX, \cF_0\to {\dots} \to \cF_m). \] In particular, when
$\cV$ is an acyclic covering of $\cX$, ${\check{H}}^n_\cV(\cX, \cF_0
\to {\dots} \to \cF_m) = H^n(\cX, \cF_0\to {\dots} \to \cF_m)$.
\end{proposition}
\begin{proof}
Take the \v Cech resolution $C^q(\cV_p, \cF_r^p)$ of a covering $\cV$ of $X_\bullet$.
Then the triple complex in \eqref{hyper-check-tri-cx}, $\check{C}^{p,q,r}$ is $\Gamma(X_p,$$C^q(\cV_p, \cF^p))$. Since
every resolution maps to an injective resolution, we arrive at a
morphism of triple complexes \eqref{hyper-check-tri-cx} and \eqref{hyper-tri-cx}: $\check{C}^{p,q,r}$$\to$$C^{p,q,r}$.
Therefore there is a map on the level of cohomology.

Notice that the cohomology of the triple complex
$\check{C}^{p,q,r}$ can be also calculated as the cohomology of
the double complex $\check{C}^{n,r}:=\oplus_{p+q=n}
{\check{C}}^{p, q, r}$ with one direction $D: \check{C}^{n, r} \to
\check{C}^{n+1, r}$ the total differential in \eqref{double-cx} with
sheaf $\cF_r$ and the other direction $d: \cF_r\to \cF_{r+1}$.
Moreover,  the hypercohomology $H^n(\cX, \cF_r)$ of stack $\cX$ can be
calculated by the double complex $C^{q, r}$ in \eqref{hyper-cx}. The first pages
$E_1^{n, r}$ of $\check{C}^{., .}$ and $C^{., .}$ are $H^n_\cV(\cX, \cF_r)$=$H^n (\cX, \cF_r)$
respectively when $\cV$ is a acyclic covering of $\cX$.
Therefore we immediately have ${\check{H}}^n_\cV(\cX, \cF_0
\to {\dots} \to \cF_m) = H^n(\cX, \cF_0\to {\dots} \to\cF_m)$.
\end{proof}

For the stack $\cX=[X/G]$, we have a result relating hypercohomology of
stacks and group cohomology as in Section
\ref{sect:coho-group-stack}.

\begin{proposition}\label{prop:group-deligne}
The hypercohomology $H^{\bullet}([X/G], \cF_0\to {\dots}\to\cF_m)$
can be calculated via a spectral sequence with $E_2$ page
\[ E_2^{p, q} := H^p(G, H^q(X, \cF_0 \to {\dots} \to \cF_m)) .\]
\end{proposition}
\begin{proof}
We also denote $\cF_r$ the equivariant sheaf on $X$. Then it has
an acyclic equivariant resolution $K_r^q$ (for example Godement's
construction). Then $K_r^q$ generates the simplicial sheaves
$K^{p,q}_r$ on $G^p\times X $. Consider the triple complex
$C^{p,q,r}:=\Gamma(G^p\times X, K_r^{p,q})$=the maps $G^p
\to$$\Gamma(X, K^q)$. Then if we filter by $r$, namely, we form
the double complex $C^{n, r}$$:=$$\oplus_{p+q=n}C^{p,q,r}$ and
filtered in the direction of $r$. Then we arrive at the first page
$H^n(G^\bullet\times X, \cF_r^\bullet)$ which is the same as the
first page of \eqref{hyper-cx}. Therefore this triple complex has
the total cohomology $H^{\bullet}([X/G], \cF_0\to
{\dots}\to\cF_m)$. In particular, the double complex $\Gamma(X,
K^q_r)$ has total cohomology $H^i(X, \cF_0\to{\dots}$ $\to$
$\cF_m)$. This tells us that if we filter $C^{p,q,r}$ by $p$
first, we arrive at the first page $E^{p, n}_1$=$ C^p(G,$$
H^n(X,$$\cF_0\to$${\dots}$ $\to$ $\cF_m))$. Then $(E^{p,n}_1,
d_1)$ is exactly the complex $( C^p(G,$ $ H^n(X,$
$\cF_0\to{\dots}\to\cF_m),$ $\delta_G)$. Therefore,
$E_2^{p,n}=H^p(G, H^n(X,$$\cF_0\to{\dots}$ $\to$ $\cF_m))$ and it
converges to $H^{\bullet}([X/G], \cF_0\to {\dots}$ $\to$ $\cF_m)$.
\end{proof}


\end{document}